\documentclass[a4paper,11pt]{amsart}
\usepackage{amsmath}
\usepackage{amssymb}
\usepackage{amsfonts}
\usepackage{amsthm}
\usepackage{mathtools}
\usepackage{mathrsfs}
\usepackage{xcolor}
\usepackage{url}
\usepackage{tikz}
\usepackage[hidelinks]{hyperref}
\hypersetup{pdftitle={Archimedean Ichino--Ikeda constants and period relations},
  pdfauthor={Fabian Januszewski}}
\theoremstyle{plain}
\newtheorem{introtheorem}{Theorem}

\newtheorem{theorem}{Theorem}[section]
\newtheorem{lemma}[theorem]{Lemma}
\newtheorem{corollary}[theorem]{Corollary}
\newtheorem{proposition}[theorem]{Proposition}

\newtheorem{hypothesis}[theorem]{Hypothesis}

\theoremstyle{remark}
\newtheorem{remark}[theorem]{Remark}

\numberwithin{equation}{section}
 \newcommand{\CC}{{\rm\bf C}}
\newcommand{\RR}{{\rm\bf R}}
\newcommand{\QQ}{{\rm\bf Q}}
\newcommand{\ZZ}{{\rm\bf Z}}
\newcommand{\NN}{{\rm\bf N}}
\newcommand{\bfA}{{\rm\bf A}}
\DeclareMathOperator{\U}{\mathrm{U}}
\DeclareMathOperator{\SU}{\mathrm{SU}}
\DeclareMathOperator{\GL}{\mathrm{GL}}
\DeclareMathOperator{\Hom}{\mathrm{Hom}}
\DeclareMathOperator{\vol}{\mathrm{vol}}
\DeclareMathOperator{\Res}{\mathrm{Res}}

\DeclareMathOperator{\diag}{\mathrm{diag}}

\newcommand{\lieg}{\mathfrak g}
\newcommand{\lieh}{\mathfrak h}
\newcommand{\liek}{\mathfrak k}
\newcommand{\liep}{\mathfrak p}

\newcommand{\lieu}{\mathfrak u}
\newcommand{\gl}{\mathfrak{gl}}
\newcommand{\slie}{\mathfrak{sl}}
\newcommand{\one}{{\bf 1}}
\newcommand{\vz}{v_0}
\newcommand{\Sinf}{S_\infty}
\newcommand{\hatpi}{\widehat\pi}
\newcommand{\hatsigma}{\widehat\sigma}
\newcommand{\volc}{\vol^{\rm comp}}
\newcommand{\Rq}{R_q}
\newcommand{\cP}{\mathcal P}
\newcommand{\cI}{\mathcal I}
\newcommand{\cN}{\mathcal N}
\newcommand{\cM}{\mathcal M}

\newcommand{\xit}{\tilde\xi}
\newcommand{\Qhw}{Q_{\rm hw}}

\DeclareMathOperator{\SL}{\mathrm{SL}}
\DeclareMathOperator{\Aut}{\mathrm{Aut}}
\DeclareRobustCommand{\hypref}[2]{\hyperref[#1]{(#2)}}
 
\begin{document}

\title[Archimedean Ichino--Ikeda constants and period relations]{The archimedean Ichino--Ikeda constant for\\
$\U(n-1,1)\supset\U(n-2,1)$ and\\period relations of Grobner, Harris and Lin}
\author{Fabian Januszewski}
\address{Institut f\"ur Mathematik, Fakult\"at EIM, Paderborn University, Warburger Str.\ 100, 33098 Paderborn, Germany}
\email{fabian.januszewski@math.uni-paderborn.de}
\date{}

\subjclass[2020]{Primary: 11F67, 22E45; Secondary: 11F70, 22E46, 22E47, 11G18, 14G35}
\keywords{Ichino--Ikeda conjecture, Gan--Gross--Prasad conjecture, discrete series, unitary groups, Gelfand--Tsetlin bases, archimedean periods, Deligne's conjecture, Shimura varieties}

\begin{abstract}
  Grobner, Harris and Lin established a form of Deligne's conjecture on critical values of $L$-functions for a class of automorphic motives over CM fields under a certain regularity condition, a conjecture on the non-vanishing of twisted central $L$-values, and a hypothesis of rationality of archimedean integrals of Ichino--Ikeda type. Let $\pi$ be a discrete series representation of $\U(n-1,1)$ and $\sigma$ a discrete series representation of $\U(n-2,1)$ in the position considered by Grobner, Harris and Lin. We compute the local Ichino--Ikeda integral of $\pi\otimes\sigma^\vee$ on the minimal $K$-types. The integral equals the volume of the compact dual, divided by the dimension of a finite-dimensional representation and by an explicit product of rational functions of the infinitesimal character. As a consequence the product over the archimedean places of the normalized local integrals is algebraic, for the rational structures of Harris at a CM point and for the normalization of the measures of Grobner, Harris and Lin. This resolves the product form of the archimedean hypothesis of Grobner, Harris and Lin for suitable test vectors. Over imaginary quadratic fields we deduce, without the archimedean hypothesis, the relation between automorphic periods on which their result rests, at a given embedding and with a larger number field, assuming their non-vanishing conjecture and a condition on the de~Rham structures along the induction of their proof. Over general CM fields the Petersson norms of de~Rham rational forms differ by CM periods which we expect to be transcendental, and we obtain the relation for the periods of forms chosen along the induction, up to a factor determined by its base, under the same assumptions and a further condition there.
\end{abstract}

\maketitle
\tableofcontents

\section{Introduction}\label{ghl:sec:introduction}

Let $F^+$ be a totally real number field, $F/F^+$ a CM extension, and $\Pi$ a cohomological conjugate self-dual cuspidal automorphic representation of $\GL_n(\bfA_F)$. In \cite{ghl}, Grobner, Harris and Lin relate critical values of Rankin--Selberg $L$-functions of $\Pi$ to Petersson norms through the refined Gan--Gross--Prasad conjecture of Ichino--Ikeda and N.~Harris \cite{nharris} for a pair of unitary groups $H=\U(V)\supset H'=\U(V')$, $\dim V=n$, $\dim V'=n-1$, which are definite at all archimedean places but one place $\vz$, where they have signatures $(n-1,1)$ and $(n-2,1)$. At $\vz$ the representations are discrete series $\pi(q)_{\vz}$ and $\pi'(q)_{\vz}$, indexed by $0\leq q\leq n-2$. The archimedean local Ichino--Ikeda integrals are not computed in \cite{ghl}. Their algebraicity on the minimal $K$-types is \cite[Conjecture 4.15]{ghl} and \cite[Conjecture 4.16]{ghl2}, on which \cite[Theorem 6.1]{ghl} and \cite[Theorem 5.8]{ghl2} are conditional.

In this paper we compute the integral at $\vz$ in closed form, for every $q$ (Theorem \ref{ghl:thm:A}), and deduce the product form of the conjecture (Theorem \ref{ghl:thm:B}) for the de~Rham structures at a CM point, which we describe (cf.\ Proposition \ref{ghl:prop:DR}), and for the normalization of the measures of \cite{ghl} recalled below. As an application we obtain, without the archimedean hypothesis of \cite{ghl}, the relation between automorphic periods of \cite[Theorem 6.1]{ghl}, on which their result rests. For $F^+=\QQ$ it holds at a given embedding and with a larger number field, assuming \cite[Conjecture 2.10]{ghl} and the condition \hypref{ghl:hyp:dR}{dR} on the de~Rham structures (Theorem \ref{ghl:thm:Cintro}). Over general CM fields we obtain it, again at a given embedding and with a larger number field, for the periods of forms chosen along the induction, up to a factor determined by its base, under a further condition \hypref{ghl:hyp:B0}{B0} (Theorem \ref{ghl:thm:Cprime}).

\subsection{The local result}\label{ghl:sec:intronotation}\label{ghl:sec:introlocal}

We use the notation of \cite[\S\S3.3, 4.2]{ghl}. We describe the more detailed setting and a dictionary in section \ref{ghl:sec:setting}, Table \ref{ghl:tab:dictionary}. Let $A=(A_1>\dots>A_n)\in(\ZZ+\frac{n-1}2)^n$ be the infinitesimal character of $\Pi$ at $\vz$. For $0\leq q\leq n-2$ put
\begin{equation}\label{ghl:eq:Aq}
\begin{aligned}
  A(q)&=(A_1,\dots,\widehat{A_{q+1}},\dots,A_n\,;\,A_{q+1}),\\
  A'(q)&=(A_1-\tfrac12,\dots,\widehat{A_{q+1}-\tfrac12},\dots,A_{n-1}-\tfrac12\,;\,A_{q+1}+\tfrac12),
\end{aligned}
\end{equation}
with the compact block before the semicolon \cite[p.~31, (4.10)]{ghl} and the hat denoting omission of the entry. Let $\pi=\pi(q)_{\vz}$ be the discrete series of $\U(n-1,1)$ with Harish-Chandra parameter $A(q)$, and $\sigma$ the discrete series of $\U(n-2,1)$ with Harish-Chandra parameter $A'(q)$. By \cite[Theorems 4.11, 4.14]{ghl}, $\pi'(q)_{\vz}=\sigma^\vee$. Let $\tau$ be the minimal $K$-type of $\pi$, $\tau_\sigma$ the minimal $K'$-type of $\sigma$, and let $dh'$ be a Haar measure on $H'_{\vz}$. For $f\in\tau$ and $f'$ in the minimal $K'$-type of $\pi'(q)_{\vz}$ let
$$
  I^*(f,f')=\frac{1}{\langle f,f\rangle\langle f',f'\rangle}\int_{H'_{\vz}}
  \langle\pi(h)f,f\rangle\,\langle\pi'(q)_{\vz}(h)f',f'\rangle\,dh'
$$
be the normalized local integral of \cite[p.~30]{ghl}. We write $d(\sigma)$ for the formal degree of
$\sigma$ with respect to $dh'$, $\volc$ for the volume of the compact dual $\U(n-1)$ of $H'_{\vz}$
for the measure matched with $dh'$ (section \ref{ghl:sec:measures}), and $F_{A'(q)}$ for the
irreducible representation of $\U(n-1)$ with infinitesimal character $A'(q)$, so that
$$
  \dim F_{A'(q)}=\prod_{i<j}\frac{A'_i-A'_j}{j-i},
$$
where the entries of $A'(q)$ are sorted in decreasing order.

The argument rests on the following non-vanishing statement, which we prove for $q\geq1$ and
$A_q-A_{q+1}\geq2$ in section \ref{ghl:sec:converse} (cf.\ Theorem \ref{ghl:thm:J}), and for $q=0$ (the
holomorphic case of \cite[Lemma 4.9]{ghl}) in Proposition \ref{ghl:prop:Jq0}:
\begin{equation}\tag{J}\label{ghl:hyp:J}
  \Hom_{H'_{\vz}}(\hatsigma,\hatpi|_{H'_{\vz}})\neq0,
\end{equation}
where $\hatpi$, $\hatsigma$ are the Hilbert spaces of $\pi$ and $\sigma$ and $\Hom$ denotes bounded
intertwining operators. The statement \hypref{ghl:hyp:J}{J} is a case of the local Gan--Gross--Prasad conjecture
for discrete series of real unitary groups, which is a theorem of He \cite[Theorem 1.1]{he} (cf.\ Remark
\ref{ghl:rmk:Jholds}). Our proof gives in addition the components of the image of $\tau_\sigma$ in the $K$-types
of $\pi$, from which we compute the constant $s(q)$ of Theorem \ref{ghl:thm:A}.

\begin{introtheorem}\label{ghl:thm:A}
  Let $n\geq3$, $0\leq q\leq n-2$, and $A\in(\ZZ+\frac{n-1}2)^n$ strictly decreasing, with
  $A_q-A_{q+1}\geq2$ if $q\geq1$. Then there is an isometric
  $H'_{\vz}$-embedding $j\colon\hatsigma\to\hatpi$, unique up to a scalar of absolute value one,
  and the following holds.
  \begin{itemize}
    \item[(i)] For $P_\tau$ the orthogonal projection onto $\tau$ and $0\neq u\in\tau_\sigma$, the ratio
      $$
        s(q):=\frac{\|ju\|^2}{\|P_\tau ju\|^2}
      $$
      is finite and independent of $u$, and
      \begin{equation}\label{ghl:eq:P}
        s(q)=\prod_{i=1}^{q}\frac{(A_i-A_{q+1})(A_i-A_n)}{(A_i-A_{q+1}-1)(A_i-A_n+1)}
        \;\in\;\QQ_{>0}.
      \end{equation}
    \item[(ii)] For all $f\in\tau$ and $f'$ in the minimal $K'$-type of $\pi'(q)_{\vz}$,
      $$
        I^*(f,f')=C_\infty(q)\cdot R^\circ(f,f'),\qquad
        C_\infty(q)=\frac{1}{d(\sigma)\,s(q)}=\frac{\volc}{\dim F_{A'(q)}\cdot s(q)},
      $$
      where
      $$
        R^\circ(f,f')=\frac{|\langle f,w\rangle|^2}{\|f\|^2\|w\|^2},
      $$
      $w=P_\tau j\bar f'$, is the normalized branching coefficient on the minimal $K$-types (section \ref{ghl:sec:R}). In
      particular, $C_\infty(q)/\volc\in\QQ_{>0}$.
  \end{itemize}
\end{introtheorem}

The second expression for $C_\infty(q)$ in (ii) is Harish-Chandra's known formula for the formal degree \eqref{ghl:eq:fdeg}. For $q=0$ the product \eqref{ghl:eq:P} is empty, $s(0)=1$, and Theorem \ref{ghl:thm:A} is the expected formal degree identity. For $n=3$, $q=1$ it reads
$$
s(1)=\frac{a(a+b)}{(a-1)(a+b+1)}
$$
with $a=A_1-A_2$, $b=A_2-A_3$. For $n=4$ and $A=(\frac92,\frac52,\frac12,-\frac32)$ one finds for $q=0,1,2$:
$$
s(q)=1,\frac{12}7,\frac{64}{35}\quad\text{and}\quad C_\infty(q)/\volc=\frac1{15},\frac7{72},\frac{35}{192}.
$$
In general the value $s(q)$ does not depend on $A_{q+2},\dots,A_{n-1}$, and only the gap $A_q-A_{q+1}$ is restricted. The latter condition cannot be dropped (cf.\ Remark \ref{ghl:rmk:sharp}). For $A_q-A_{q+1}=1$ the parameter $A'(q)$ is singular, since $A_q-\frac12=A_{q+1}+\frac12$.

\subsection{The product form of the conjecture}\label{ghl:sec:introglobal}

The local condition on which the main results in \cite{ghl} rest is \cite[Conjecture 4.15]{ghl}. It asserts that
$$
I^*_v(f_v,f'_v)\in E(\pi(q))\cdot E(\pi'(q))
$$
for all $v\in\Sinf$, for $f_v$, $f'_v$ in the minimal $K$-types and rational for the de~Rham rational structures of \cite[\S3.2.2]{ghl}. The proof of \cite[Theorem 6.1]{ghl} uses only the product $\prod_{v\in\Sinf}I^*_v$. As in \cite[\S1.1]{ghl} all number fields are subfields of $\bar\QQ\subset\CC$, and the statements below are understood relative to this embedding.

We take the de~Rham structures of \cite[\S2.3]{harris13a} on $\pi(q)$ and $\pi'(q)$ at a common CM point, as in \cite[pp.~2045--2046]{harris13a} and \cite[\S3.1, p.~87]{harris14}. On the minimal $K$-types they are the rational structures given by the fibres at this point of models of the automorphic vector bundles over the Shimura varieties of $H$ and $H'$. The point is attached to an orthogonal decomposition of $V$ into lines, $n-1$ of which span $V'$, and exactly one of the lines, $V_-$ say, is negative at $\vz$. It lies in $V'$. We describe these structures in section \ref{ghl:sec:DR} (Proposition \ref{ghl:prop:DR}), and we extend their scalars to the number field $E$ defined there, which contains $E(\pi(q))E(\pi'(q))$ and is stable under complex conjugation. In general $E$ is larger than $E(\pi(q))E(\pi'(q))$, since the fibres at the point are vector spaces over a field which contains the residue field of the point, and $E$ need not be CM (Remark \ref{ghl:rmk:nonCM}).

The measures on $H(\bfA_{F^+})$ and $H'(\bfA_{F^+})$ are those of \cite[p.~28]{ghl}, normalized as in \cite[\S5.2, p.~2043]{harris13a}, with the convention
$$
\vol(\U(V'_1)(F^+)\backslash\U(V'_1)(\bfA_{F^+}))=1
$$
for the group $\U(V'_1)\cong\U(V_1)$ of \cite[\S1.2]{ghl}. Hence the measure $dh'=\prod_vdh'_v$ on $H'(\bfA_{F^+})$ which enters $I^{\rm can}$ is a rational multiple of the Tamagawa measure, and its local factors at the finite places take rational values on compact open subgroups and give volume $1$ to hyperspecial maximal compact subgroups at the places where $H'$ is unramified (section \ref{ghl:sec:globalmeasures}). The archimedean factors are not constrained.

\begin{introtheorem}\label{ghl:thm:B}
  Let $n\geq3$ and $0\leq q\leq n-2$, and let $\Pi$, $\pi(q)$, $\pi'(q)$, $f=\otimes f_v$, $f'=\otimes f'_v$ be as in \cite[\S4.2, Theorem 4.14]{ghl}, with $f_v$, $f'_v$ rational for the de~Rham structures over $E$ of section \ref{ghl:sec:DR} (Proposition \ref{ghl:prop:DR}) and, for every $v\in\Sinf$, eigenvectors of $\U(V_-)(F^+_v)$, and let $dh'$ be normalized as above. Then
  $$
    \prod_{v\in\Sinf}I^*_v(f_v,f'_v)=\mathrm{Vol}_\infty\cdot\frac{1}{s(q)}\cdot
    \prod_{v\in\Sinf}\frac{1}{\dim F_{\Lambda'_v}}\cdot\prod_{v\in\Sinf}R^\circ_v(f_v,f'_v)
    \;\in\;E,
  $$
  with $\mathrm{Vol}_\infty\in\QQ^\times$ independent of $q$ (the product of the volumes of the compact duals, Lemma \ref{ghl:lem:M}), $\Lambda'_{\vz}=A'(q)$, $\Lambda'_v$ the infinitesimal character of $\pi'(q)_v^\vee$ for $v\neq\vz$, and $R^\circ_v\in E$. The identity holds for every Haar measure on $H'(F^+_\infty)$. The normalization of $dh'$ is required for $\mathrm{Vol}_\infty\in\QQ^\times$, and hence for the $E$-rationality of the expression.
\end{introtheorem}

At $\vz$ the eigenvector condition is empty, since $\U(V_-)(F^+_{\vz})$ is a central factor of $K_{\vz}$ and of $K'_{\vz}$. For $F^+=\QQ$ there is no other place. If $f$, $f'$ as in \cite[Theorem 4.14]{ghl} exist, then such $f$, $f'$ with $I^{\rm can}(f,f')\neq0$ exist (section \ref{ghl:sec:proofB}). This is the product form of \cite[Conjecture 4.15]{ghl} for these vectors. For $F^+=\QQ$ it is the conjecture as stated, for the field $E$. For $[F^+:\QQ]\geq2$ the per-place statement is meaningful only after a splitting of the archimedean measure, which the normalization does not fix, and it holds for every splitting with rational volumes at the definite places (Corollary \ref{ghl:cor:perplace}). At the definite places the eigenvector condition cannot be removed without an assumption on CM periods which we expect to fail (Remark \ref{ghl:rmk:definiteperiods}). For $n=2$ the only index is $q=0$, and the product form holds as well (Remark \ref{ghl:rmk:n2}).

Theorem \ref{ghl:thm:B} does not use \cite[Proposition 3.19]{harris13a}, which asserts the rationality of the normalized inner products $Q(v)^{-1}\langle v_1,v_2\rangle$ for all de~Rham rational vectors $v\neq0$, $v_1$, $v_2$, and on which the independence of the automorphic periods $Q(\pi)$ of \cite[\S4.3]{ghl} of the choice of the de~Rham rational vector rests. Its proof, as we understand it, needs the de~Rham structure $U(\lieg)_{{\rm DR},x}$ of the enveloping algebra of the Lie algebra $\lieg$ of the group at a CM point $x$ \cite[p.~2011]{harris13a} to be stable under complex conjugation. With the description given there this is not the case already for $\mathrm{GU}(1,1)$, and there the statement fails: for a de~Rham rational lowest weight vector $w$ of a holomorphic discrete series and the non-compact root $a$ with $X_{-a}w=0$, the ratio $\langle p(a)X_aw,p(a)X_aw\rangle/\langle w,w\rangle$ is a positive integer multiple of $|p(a)|^2$, where $p(a)$ is the CM period by which the root vector $X_a$ is multiplied, and $|p(a)|^2$ is transcendental (Remark \ref{ghl:rmk:definiteperiods} (3)). We use the rationality of inner products only on the minimal $K$-types, for eigenvectors of $\U(V_-)$ (Lemma \ref{ghl:lem:Rrational}), and for Petersson norms of forms whose archimedean components lie in the minimal $K$-types, which for $[F^+:\QQ]\geq2$ depend on their patterns (Lemma \ref{ghl:lem:pattern}). Theorem \ref{ghl:thm:B} is stated for the given embedding. The version equivariant under $\Aut(\CC/F^{\rm Gal})$, which is used in \cite[Theorem 6.1]{ghl}, requires in addition that the de~Rham structures are transported by $\Aut(\CC/F^{\rm Gal})$, and that $E$ be CM or totally real. For $F^+=\QQ$ it is obtained under the condition \hypref{ghl:hyp:Gal}{Gal} of section \ref{ghl:sec:global} in the proof of Theorem \ref{ghl:thm:C} (ii), and \hypref{ghl:hyp:Gal}{Gal} cannot hold in the example of Remark \ref{ghl:rmk:nonCM}.

In section \ref{ghl:sec:global} we determine what Theorem \ref{ghl:thm:B} gives for \cite[Theorem 6.1]{ghl}. The
proof of \cite[Theorem 6.1]{ghl} is an induction on $n$ \cite[\S6.2]{ghl}, which passes from $\Pi$ and $q$ to a
representation $\Pi'$ of $\GL_{n-1}(\bfA_F)$ and $n-q-2$. Besides \cite[Conjecture 4.15]{ghl} it uses at every step \cite[Conjecture
2.10]{ghl}, and the de~Rham structures of \cite[\S3.2.2]{ghl} on the representations of unitary groups which occur.
In section \ref{ghl:sec:globalinputs} we formulate the condition \hypref{ghl:hyp:dR}{dR} (Hypothesis \ref{ghl:hyp:dR}), which requires, inter alia, that the latter
are those of \cite{harris13a} at the CM points attached to one decomposition of $V$ into lines, as in section
\ref{ghl:sec:DR}, and, for the Galois-equivariant version, the condition \hypref{ghl:hyp:Gal}{Gal} (Hypothesis \ref{ghl:hyp:Gal}), which requires in particular that the
conjugate data under $\Aut(\CC)$ satisfy the hypotheses of Theorem \ref{ghl:thm:B} and \hypref{ghl:hyp:dR}{dR}, that $\Aut(\CC/F^{\rm Gal})$ transports the rational
data, and that $E$ is CM or totally real. We use the notation of \cite{ghl}: $\check\xi_\Pi$ is the conjugate
inverse of the central character of $\Pi$ \cite[\S1.1]{ghl}, $p(\check\xi_\Pi,\Sigma)$ its CM period
\cite[\S2.1.1]{ghl}, $P^{(q)}(\Pi,\imath_{\vz})$ the local arithmetic automorphic period \cite[Theorem 2.6]{ghl} at
the embedding $\imath_{\vz}$ of $F$ attached to $\vz$, which belongs to the CM type $\Sigma$ of \cite[\S1.1]{ghl}, and $E_q(\Pi)=E(\Pi)E(\pi(q))$ \cite[p.~46]{ghl}. For $F^+=\QQ$ we
obtain the following.

\begin{introtheorem}\label{ghl:thm:Cintro}
  Let $F^+=\QQ$, so that $F$ is imaginary quadratic, let $n\geq2$, and let $\Pi$ be as in \cite[Theorem 6.1]{ghl}: it satisfies \cite[Hypothesis
  2.4]{ghl}, and $\Pi_\infty$ is $(n-1)$-regular. For $0\leq q\leq n-2$ let $\pi(q)$ be as in \cite[\S6.1]{ghl}. Assume \cite[Conjecture 2.10]{ghl} for the representations to
  which the proof of \cite[Theorem 6.1]{ghl} applies it (Hypothesis \ref{ghl:hyp:nonvanishing}), and \hypref{ghl:hyp:dR}{dR}. Let
  the fields $E(\pi)$ of the representations $\pi$ of unitary groups which occur in the induction be enlarged, as
  in section \ref{ghl:sec:DR}, so that for each pair of the induction $E(\pi(q))=E(\pi'(q))$ is the field $E$
  defined there. Let $Q(\pi(q))$ be the class modulo $E^\times$ of the Petersson norm of a non-zero de~Rham rational form in $\pi(q)$
  whose archimedean component lies in the minimal $K$-type. By Lemma \ref{ghl:lem:pattern} it does not depend on
  the form. Then for $0\leq q\leq n-2$
  $$
    Q(\pi(q))\in p(\check\xi_\Pi,\Sigma)^{-1}\,\frac{P^{(q+1)}(\Pi,\imath_{\vz})}{P^{(q)}(\Pi,\imath_{\vz})}
    \cdot E^\natural(\Pi,q)^\times
  $$
  at the given embedding, for a number field $E^\natural(\Pi,q)\supset E_q(\Pi)$ determined by the induction. If
  in addition \hypref{ghl:hyp:Gal}{Gal} holds, the relation holds with $E_q(\Pi)$ in place of $E^\natural(\Pi,q)$, as a relation
  equivariant under $\Aut(\CC/F^{\rm Gal})$, and this is \cite[Theorem 6.1]{ghl} for $F^+=\QQ$, for
  these fields and forms, without \cite[Conjecture 4.15]{ghl}.
\end{introtheorem}

The precise statement is Theorem \ref{ghl:thm:C}. It also gives (6.4) of \cite[Theorem 6.3]{ghl} at the given
embedding and with a larger number field. The relation (6.5) there, with the motivic periods of $M(\Pi)$, uses further inputs
\cite[pp.~46--47]{ghl} and is not treated here. For $[F^+:\QQ]\geq2$ the Petersson norms of the de~Rham rational
forms which enter the induction depend on their patterns, through CM periods which we expect to be transcendental
(Lemma \ref{ghl:lem:pattern}). Under the same assumptions and a further condition \hypref{ghl:hyp:B0}{B0} (Hypothesis
\ref{ghl:hyp:B0}) on the base of the induction we obtain the relation for $Q(\pi(q))$ taken on forms of a pattern
fixed by the induction, up to a factor determined by the base of the induction (Theorem \ref{ghl:thm:Cprime}). If
the CM periods are transcendental, \cite[Theorem 6.1]{ghl}, in the form proved in \cite[\S6.2]{ghl}, can hold only
for the forms of at most one pattern, and we interpret Theorem \ref{ghl:thm:Cprime} as a statement of this kind
(Remark \ref{ghl:rmk:GHL43}).

\subsection{The argument}\label{ghl:sec:introargument}

The proof of Theorem \ref{ghl:thm:A} has five steps. The notation is that of sections \ref{ghl:sec:setting} to
\ref{ghl:sec:proofA}.

\noindent {\rm(1)} {\em The orthogonality relation} (section \ref{ghl:sec:lemmaP}). Matrix coefficients of the discrete series $\pi=\pi(q)_{\vz}$ of $\U(n-1,1)$ with
    Harish-Chandra parameter $A(q)$, restricted to $H'_{\vz}$, are square integrable. By an equivariant closed operator into $L^2(H'_{\vz})$ the
    part of $\hatpi|_{H'_{\vz}}$ orthogonal to the $\sigma$-isotypic part contributes nothing, and
    multiplicity one \cite[Theorem B]{sunzhu} gives $I^*=(d(\sigma)s(q))^{-1}R^\circ$, once \hypref{ghl:hyp:J}{J} holds and
    $P_\tau ju\neq0$. No direct integral decomposition of $\hatpi|_{H'_{\vz}}$ enters our argument.

\noindent {\rm(2)} {\em The Gelfand--Tsetlin model} (section \ref{ghl:sec:gt}). Gelfand--Tsetlin formulas
    with top row $A-\frac{n-1}2$ on an explicit region $\Rq$ of $K$-types (cf.\ Figure \ref{ghl:fig:Rq}) define an irreducible
    unitarizable $(\gl_n,K)$-module $V_{\Rq}$, which is $\pi_K$ by Schmid's characterization of the discrete
    series by its $K$-types \cite[Theorem (1.3)]{schmid75}. We fix an isomorphism $\Phi\colon V_{\Rq}\to\pi_K$
    (Theorem \ref{ghl:thm:identification}). In this way $\pi_K$ is given explicitly: its $K$-types,
    the top rows in $\Rq$, each with multiplicity one; in each $K$-type a basis adapted to the branching to
    $K'$; the action of every $E_{ab}$ on this basis by rational coefficients; and the squared norms of the basis
    vectors for the invariant form \eqref{ghl:eq:NP}, which are positive rational numbers.

\noindent {\rm(3)} {\em The Schmid vector} (section \ref{ghl:sec:schmid}). Let $\tau_\sigma=(\mu;\tau_0)$,
    $\mu=(\tau_1,\dots,\tau_{n-2})$ \eqref{ghl:eq:tausigma}, and let $u$ be a highest weight vector of
    $\tau_\sigma$. Assume first that $j$ exists. For a $K$-type $\delta$ of $\pi$ the component $P_\delta ju$ is $0$
    or a $K'$-highest weight vector of weight $(\mu;\tau_0)$ in $\pi_K(\delta)$. By the multiplicity one of step (2)
    and the multiplicity free branching from $K$ to $K'$, the latter occurs only for $\delta=(m;\tau_0)$ with $m\in\Rq$
    interlacing $\mu$, and then $P_\delta ju$ is a multiple of $\Phi(\xi_\delta)$, where $\xi_\delta=\xi(m,\mu)$ is the
    basis vector of the model with top row $m$ whose lower rows form the highest pattern of weight $\mu$. We identify
    $\delta$ with $m$. These $m$ form an explicit set $S$, a box in the coordinates $d_r=m_r-\tau_r$, $1\leq r\leq
    q$, the other entries of $m$ being determined by them; $\delta=\tau$ is $d=0$ (Lemma \ref{ghl:lem:support}, cf.\ Figure \ref{ghl:fig:S}). Hence
    $$
      P_\delta ju=c(\delta)\,\Phi(\xi_\delta)\quad(\delta\in S),\qquad P_\delta ju=0\quad(\delta\notin S),
    $$
    for a function $c\colon S\to\CC$ \eqref{ghl:eq:wdelta}. The Schmid conditions for $ju$ say that $\pi(\liep')ju$
    has no component in the $K'$-types $\gamma^{(i)}=\tau_\sigma-\beta_i$, $\beta_i$ the non-compact positive roots
    of $\sigma$, which do not occur in $\sigma$ (Lemma \ref{ghl:lem:Kprime}, \eqref{ghl:eq:schmidcond}). Through the
    explicit action of the $E_{n,i}$, $i\leq q$, on the $\xi_\delta$ (Lemma \ref{ghl:lem:matrixelement}) they give
    linear equations for $c$, the reduced equations, which are Cauchy systems (section \ref{ghl:sec:uniqueness}). Their
    solutions are the multiples of
    $$
      c_\delta=\frac{\Delta_I(x)}{\prod_{r=1}^qd_r!},
    $$
    where $x_i=m_i-i+1$ and $\Delta_I(x)=\prod_{a<b\in I}(x_a-x_b)$, $I=\{1,\dots,q,n-1\}$ (Proposition
    \ref{ghl:prop:unique}, Lemma \ref{ghl:lem:WS} (c)). The adjointness $E_{n-1,n}^*=-E_{n,n-1}$
    \eqref{ghl:eq:adjoint} and the explicit coefficients of these two operators give
    $\|\Phi(\xi_\delta)\|^2=\kappa\,\cN(\xi_\delta)$ on $S$, with $\cN$ explicit \eqref{ghl:eq:NS} and an unknown
    $\kappa>0$ (Lemma \ref{ghl:lem:normS}). Put
    $$
      t_\delta:=\frac{|c_\delta|^2\cN(\xi_\delta)}{|c_\tau|^2\cN(\xi_\tau)}\;\in\;\QQ_{>0},
    $$
    given in closed form by \eqref{ghl:eq:tdelta} (Theorem \ref{ghl:thm:term}).
    The set $S$, the reduced equations,
    $c_\delta$, $\cN$ and $t_\delta$ do not involve $j$.
    If $j$ exists, then $c(\delta)=\lambda c_\delta$ with
    $\lambda\neq0$, as $ju\neq0$, so $\|P_\delta ju\|^2=|\lambda|^2\kappa|c_\delta|^2\cN(\xi_\delta)$. As
    $c_\tau\neq0$, $P_\tau ju\neq0$, and $t_\delta=\|P_\delta ju\|^2/\|P_\tau ju\|^2$, in which $\lambda$ and
    $\kappa$ cancel. Parseval gives
    $$
      s(q)=\frac{\|ju\|^2}{\|P_\tau ju\|^2}=\sum_{\delta\in S}\frac{\|P_\delta ju\|^2}{\|P_\tau ju\|^2}
      =\sum_{\delta\in S}t_\delta
    $$
    (Theorem \ref{ghl:thm:schmid}, Remark \ref{ghl:rmk:schmidroute}).

\noindent {\rm(4)} {\em The summation} (section \ref{ghl:sec:summation}). It remains to show that
    $\sum_{\delta\in S}t_\delta$ equals \eqref{ghl:eq:P}. By \eqref{ghl:eq:tdelta}, $t_\delta$ is a function of
    $d_1,\dots,d_q$, of $a_r=A_r-A_{q+1}$ and of $B=A_{q+1}-A_n$, through $u_r=a_r+d_r$ and
    $d_0=d_1+\dots+d_q$. Multiplying by a factor depending only on $a$ and $B$ turns the series into $S^{(q)}(a;B)$
    \eqref{ghl:eq:Sq}, with summands free of rising factorials, and the claim into the closed form
    \eqref{ghl:eq:Sqclosed} (section \ref{ghl:sec:normalization}).
    The series $S^{(q)}(a;B)=\cI_{q,0}(a;B,0)$ is a member of a family
    $\cI_{p,m}(a;B,D;z)$, a sum over $d_1,\dots,d_p$, in which the summation variables of the outer sums enter only
    through auxiliary variables $z_1,\dots,z_m$ and a shift $D$ of $d_0$: splitting off the sum over $d_p$, the value
    $u_p$ becomes a new variable $z_{m+1}$, and $d_p$ is added to $D$ \eqref{ghl:eq:split}. For this family
    $$
      \cI_{p,m}(a;B,D;z)=\cM_p(a;B)\,D!\,(B+D)!\,\Omega_{p,m}(a;B,D;z)
    $$
    with explicit $\cM_p$, $\Omega_{p,m}$ (Theorem
    \ref{ghl:thm:induction}). For $p=1$ the sum over $d_1$ telescopes by
    the recurrence \eqref{ghl:eq:recurrence}, with vanishing limit at infinity (Lemma \ref{ghl:lem:base}). For
    $p\geq2$ the closed form for $p-1$, specialized at $z_{m+1}=u_p$, makes the sum over $d_p$ telescope again, with
    vanishing upper boundary term. Only the boundary terms at $d_p=0$ remain. For $p=q$, $m=0$, $D=0$ one has
    $\Omega_{q,0}=\prod_r(a_r-1)^{-1}$, and undoing the normalization gives \eqref{ghl:eq:P}, with the factor
    $a_r(a_r+B)/((a_r-1)(a_r+B+1))$ for each $r\leq q$ (Corollary \ref{ghl:cor:summation}).

\noindent {\rm(5)} {\em The Schmid vector generates $\sigma$} (section \ref{ghl:sec:converse}).
    The formal Schmid vector
    $$
      v=\sum_{\delta\in S}c_\delta\,\Phi(\xi_\delta)
    $$
    has orthogonal terms, and $\|v\|^2=|c_\tau|^2\|\Phi(\xi_\tau)\|^2\sum_\delta t_\delta$, which is finite by
    step (4) (Lemma \ref{ghl:lem:X}). It lies in the closed subspace $\mathcal X$ of vectors of $K'$-type $\tau_\sigma$
    and weight $(\mu;\tau_0)$, the Hilbert sum of the lines $\CC\Phi(\xi_\delta)$, $\delta\in S$. On $\mathcal X$ the
    weak Schmid conditions \eqref{ghl:eq:WS} are equivalent to the reduced equations of step (3) for the coefficients
    (Lemma \ref{ghl:lem:WS}), so their solutions in $\mathcal X$ form the line $\CC v$ (Lemma \ref{ghl:lem:sol}). The
    bounded operators commuting with $H'_{\vz}$ preserve this line, so the closed $H'_{\vz}$-span $V$ of $v$ is irreducible
    (Proposition \ref{ghl:prop:Virred}). It is a discrete series (Lemma \ref{ghl:lem:Vds}), and the Schmid conditions
    for $v$ and the $K'$-types of the discrete series of $H'_{\vz}$ and the transitions between them identify it
    with $\hatsigma$ (Proposition \ref{ghl:prop:Hprime}, Lemmas \ref{ghl:lem:classif}, \ref{ghl:lem:Vschmid},
    \ref{ghl:lem:Id}). This proves \hypref{ghl:hyp:J}{J} for $q\geq1$ and $A_q-A_{q+1}\geq2$, with
    $v=ju$ for an isometric embedding $j$ and a highest weight vector $u$ of $\tau_\sigma$ (Theorem \ref{ghl:thm:J}).

\begin{figure}[ht]
  \centering
  \begin{tikzpicture}[>=stealth,x=1.1cm,y=1cm,
      box/.style={draw,rounded corners=1pt,align=center,inner sep=4pt,font=\small}]
    \node[box] (s4) at (0,4)     {(2) Section \ref{ghl:sec:gt}: $V_{\Rq}\cong\pi_K$\\Theorem \ref{ghl:thm:identification}};
    \node[box] (s5) at (4.7,4)   {(3) Section \ref{ghl:sec:schmid}: $S$, $c_\delta$, $\cN$, $t_\delta$\\Lemma \ref{ghl:lem:support} to Theorem \ref{ghl:thm:term}};
    \node[box] (s6) at (9.4,4)   {(4) Section \ref{ghl:sec:summation}: $\sum_\delta t_\delta=$ \eqref{ghl:eq:P}\\Corollary \ref{ghl:cor:summation}};
    \node[box] (q0) at (0,2)     {$q=0$: \hypref{ghl:hyp:J}{J}, $s(0)=1$\\Proposition \ref{ghl:prop:Jq0}};
    \node[box] (s3) at (4.7,2)   {(1) Section \ref{ghl:sec:lemmaP}\\Lemmas \ref{ghl:lem:L2}, \ref{ghl:lem:closed}, \ref{ghl:lem:P}};
    \node[box] (s7) at (9.4,2)   {(5) Section \ref{ghl:sec:converse}: \hypref{ghl:hyp:J}{J}, $q\geq1$\\Theorem \ref{ghl:thm:J}};
    \node[box] (cf) at (4.7,0)   {$I^*=(d(\sigma)s(q))^{-1}R^\circ$\\Corollary \ref{ghl:cor:Cinf}};
    \node[box] (sq) at (9.4,0)   {$s(q)=\sum_\delta t_\delta$\\Theorem \ref{ghl:thm:schmid}};
    \node[box] (tA) at (7.05,-1.6) {Theorem \ref{ghl:thm:A}};
    \draw[->] (s4) -- (s5);
    \draw[->] (s5) -- (s6);
    \draw[->] (s4) -- (q0);
    \draw[->] (s5) -- (q0) node[pos=0.45,below,sloped,font=\scriptsize] {Lemma \ref{ghl:lem:Kprime}};
    \draw[->] (s5) -- (s7);
    \draw[->] (s6) -- (s7);
    \draw[->] (s3) -- (s7);
    \draw[->] (q0) -- (cf);
    \draw[->] (s3) -- (cf);
    \draw[->] (s7) -- (cf);
    \draw[->] (s7) -- (sq);
    \draw[->] (cf) -- (tA);
    \draw[->] (sq) -- (tA);
  \end{tikzpicture}
  \caption{The proof of Theorem \ref{ghl:thm:A}. An arrow points from a result to one whose proof uses it. The numbers (1) to (5) are those of the steps in section
    \ref{ghl:sec:introargument}, and the boxes without a number contain results of section \ref{ghl:sec:proofA}.
    Theorem \ref{ghl:thm:J} proves \hypref{ghl:hyp:J}{J} without assuming it, and no result of sections \ref{ghl:sec:lemmaP} to
    \ref{ghl:sec:converse} uses \hypref{ghl:hyp:J}{J}. Section \ref{ghl:sec:conditions} is
    not needed for Theorem \ref{ghl:thm:A}.}\label{ghl:fig:steps}
\end{figure}
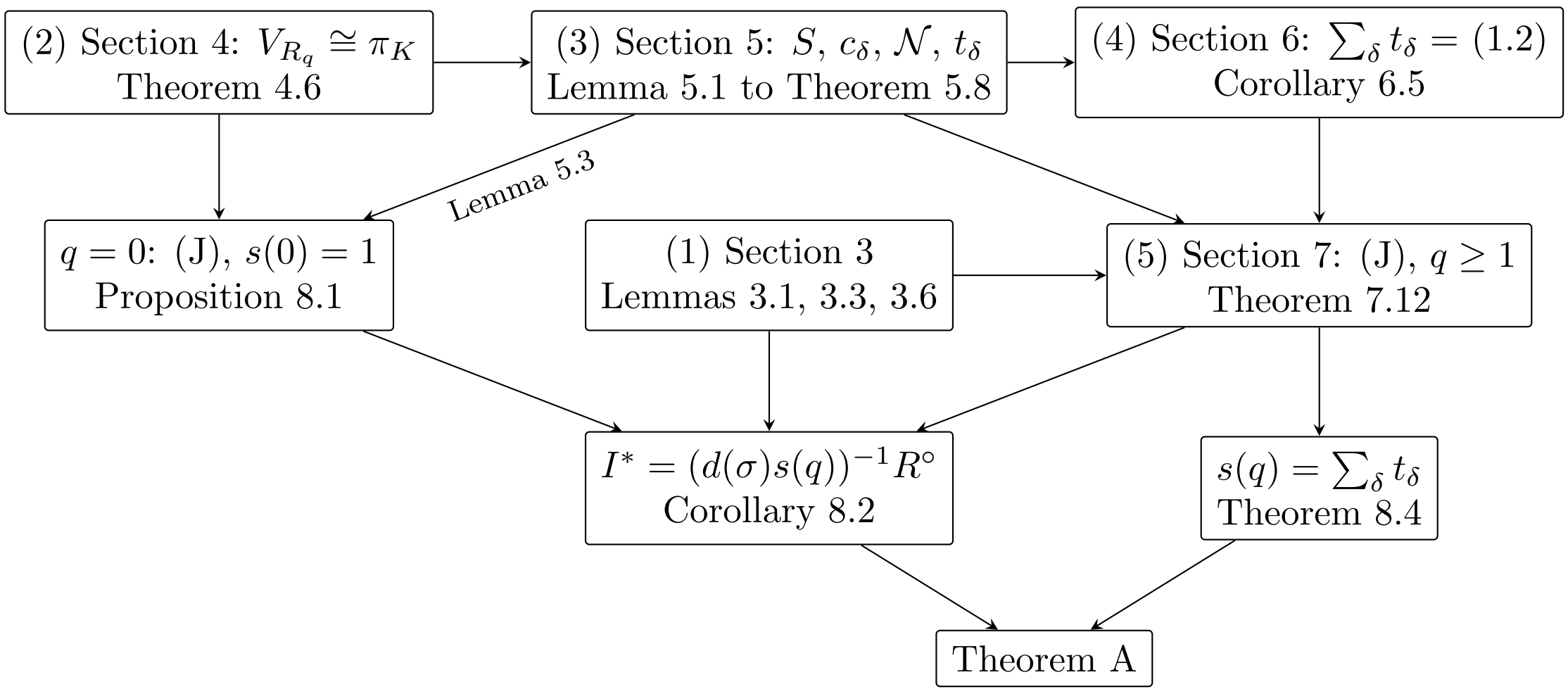

\medskip
\noindent Steps (1) to (5) are carried out in sections \ref{ghl:sec:lemmaP} to \ref{ghl:sec:converse}, in this
order, as far as they do not use \hypref{ghl:hyp:J}{J}. Theorem \ref{ghl:thm:A} is proved in section \ref{ghl:sec:proofA}. For
$q\geq1$, by step (5) the vector $ju$ is a multiple of $v$, so $s(q)=\sum_\delta t_\delta$ (Theorem
\ref{ghl:thm:schmid}), which is \eqref{ghl:eq:P} by step (4), and step (1) gives the factorization of $I^*$ (Corollary \ref{ghl:cor:Cinf}). The argument
of step (3) for a given $j$, from the Schmid conditions to the reduced equations, is in section
\ref{ghl:sec:conditions}. With \hypref{ghl:hyp:J}{J} taken from \cite[Theorem 1.1]{he}, it gives Theorem
\ref{ghl:thm:A} without step (5) (Remark \ref{ghl:rmk:Jholds} (3)).
For $q=0$ steps (3) to (5) are
replaced by Proposition \ref{ghl:prop:Jq0}, which gives \hypref{ghl:hyp:J}{J} and $s(0)=1$ directly. Figure \ref{ghl:fig:steps}
shows the dependencies.

Theorem \ref{ghl:thm:B} follows in section \ref{ghl:sec:consequences}, with the treatment of the
definite places by Schur orthogonality and the rationality of the product of archimedean volumes
(Lemma \ref{ghl:lem:M}).

\subsection{Prior work}\label{ghl:sec:prior}

For $n=3$, $q=0$, the local integral is computed in \cite[Appendix B, Proposition B.3]{hhy}, by
theta lifts and a seesaw. Theorem \ref{ghl:thm:A} specializes to a case of that result, up to the comparison of normalizations. For $n=3$,
$q=1$, the $\tau$-spherical function of $\pi$ is in \cite[Theorem 4.10]{liu}.
The totally definite case of
the product rationality is \cite[Corollaries 4.1.6, 4.1.7, p.~145]{harris13b} and \cite[Proposition
6.16]{grobnerlin}. The existence of the embedding $j$ of Theorem \ref{ghl:thm:A}, i.\,e.\ \hypref{ghl:hyp:J}{J}, is
contained in the branching law of He \cite[Theorem 1.1]{he} for the discrete spectrum of the restriction of
discrete series from $\U(p,q)$ to $\U(p-1,q)$, and the uniqueness of the distinguished pair in the Vogan
packet is a case of the theorem of Beuzart-Plessis \cite[Theorem 12.4.1]{beuzartplessis}
(Remark \ref{ghl:rmk:Jholds}). We have not found the formula \eqref{ghl:eq:P}, or the rationality of $C_\infty(q)$
for $q\geq1$, in the literature.

\subsection*{Acknowledgements}

The author acknowledges support by the Deutsche Forschungsgemeinschaft (DFG, German Research
Foundation) -- SFB-TRR 358/1 2023 -- 491392403.
 \section{Setting and the local Ichino--Ikeda integral}\label{ghl:sec:setting}

\subsection{Groups}\label{ghl:sec:groups}

Fix $n\geq3$. Let $e_1,\dots,e_n$ be the standard basis of $\CC^n$ with the Hermitian form of matrix
$J=\diag({\bf 1}_{n-1},-1)$, and $H=\U(n-1,1)$ its unitary group. Let $H'\cong\U(n-2,1)$ be the
stabilizer of the positive vector $e_{n-1}$. It acts on $\langle e_1,\dots,e_{n-2},e_n\rangle$. Put
$K=\U(n-1)\times\U(1)_{e_n}$ and $K'=K\cap H'=\U(n-2)\times\U(1)_{e_n}$. All four groups are
connected. This is the place $\vz$ of \cite[\S3.3]{ghl}: $V=V'\oplus\CC e_{n-1}$ with $V'$ of
signature $(n-2,1)$, so that $H=H_{\vz}$ and $H'=H'_{\vz}$. We drop the index $\vz$ from now on until
section \ref{ghl:sec:consequences}.

We identify $\lieh_\CC=\gl_n(\CC)$ with basis $E_{ab}$. Then $\lieh'_\CC$ is spanned by the $E_{ab}$
with $a,b\in\{1,\dots,n-2,n\}$, and $\liep'_\CC=\langle E_{in},E_{ni}:i\leq n-2\rangle$. On a unitary
$(\lieh,K)$-module
\begin{equation}\label{ghl:eq:adjoint}
  E_{ab}^*=E_{ba}\ (a,b\leq n-1),\qquad E_{an}^*=-E_{na}\ (a\leq n-1),\qquad E_{nn}^*=E_{nn}.
\end{equation}
A $K$-type is written $(m;m_0)=(m_1\geq\dots\geq m_{n-1};m_0)$, where $m$ is the highest weight of
the $\U(n-1)$-factor and $m_0$ the character of $\U(1)_{e_n}$, and likewise $(\mu;\mu_0)$ for $K'$. The
branching from $\U(n-1)\times\U(1)$ to $\U(n-2)\times\U(1)$ is multiplicity free, the
$\U(n-2)$-types in $F_m$ being the $\mu$ which interlace $m$, $m_i\geq\mu_i\geq m_{i+1}$.

Let $a_t=\exp t(E_{1n}+E_{n1})\in H'$. Both $H$ and $H'$ have real rank one, $H'=K'A^+K'$ with
$A^+=\{a_t:t\geq0\}$, and the restricted roots of $\lieh$ on $\log a_t$ are $\pm t$ with
multiplicity $2(n-2)$ and $\pm2t$ with multiplicity one. Hence $\rho_H(\log a_t)=(n-1)t$ and
$\rho_{H'}(\log a_t)=(n-2)t$.

\subsection{Parameters}\label{ghl:sec:parameters}

Let $A=(A_1>\dots>A_n)\in(\ZZ+\frac{n-1}2)^n$ and put $L_j:=A_j-\frac{n-1}2\in\ZZ$, so that
$L_1>\dots>L_n$. Fix $0\leq q\leq n-2$, and let $A(q)$, $A'(q)$ be as in \eqref{ghl:eq:Aq}. We
formulate the following standing assumption, which is weaker than the regularity
$A_i-A_{i+1}\geq2$ of \cite[p.~31]{ghl}:
\begin{equation}\label{ghl:eq:gap}
  A_q-A_{q+1}\geq2\qquad\text{if }q\geq1.
\end{equation}

{\em The representation $\pi$.} Let $\pi=\pi_{A(q)}$ be the discrete series of $H$ with
Harish-Chandra parameter $A(q)$, on the Hilbert space $\hatpi$, with Harish-Chandra module
$\pi_K$. The positive system $\Psi=\{\alpha:\langle A(q),\alpha\rangle>0\}$ has the standard compact
part $\{e_i-e_j:i<j\leq n-1\}$ and the non-compact part
$$
  \beta_i=e_i-e_n\ (i\leq q),\qquad\beta_i=e_n-e_i\ (q<i\leq n-1).
$$
With $\rho_{c,i}=\frac n2-i$, the Blattner parameter $\tau=A(q)+\rho_n-\rho_c$ is
\begin{equation}\label{ghl:eq:tau}
  \tau_i=L_i+i\ (i\leq q),\qquad\tau_i=L_{i+1}+i-1\ (q<i\leq n-1),\qquad\tau_0=L_{q+1}+n-1-q.
\end{equation}
It is the minimal $K$-type of $\pi$ and has multiplicity one. We put
$$
  \mathrm{tot}:=\sum_{i\leq n-1}\tau_i+\tau_0=\sum_jL_j+\frac{n(n-1)}2.
$$

{\em The representation $\sigma$.} Let $\sigma=\pi_{A'(q)}$ be the discrete series of $H'$ with
Harish-Chandra parameter $A'(q)$. By \eqref{ghl:eq:gap},
$$
  A_q-\frac12>A_{q+1}+\frac12>A_{q+2}-\frac12.
$$
Let $\Delta^+(\sigma)$ be the set of roots $\alpha$ of
$\lieh'_\CC$ with $\langle A'(q),\alpha\rangle>0$, where the roots of
$\lieh'_\CC$ are those involving only $e_1,\dots,e_{n-2},e_n$, with the index $n$ of $\gl_n(\CC)$. Its non-compact part, i.\,e.\ the set of positive
non-compact roots of $\sigma$, is
$$
  \Delta_n^+(\sigma)=\{\beta_i:1\leq i\leq n-2\},\qquad\beta_i=e_i-e_n\ (i\leq q),\quad
  \beta_i=e_n-e_i\ (q<i\leq n-2),
$$
and its Blattner parameter $\tau_\sigma=(\mu;\tau_{\sigma,0})$ is, by the same computation as for
\eqref{ghl:eq:tau},
\begin{equation}\label{ghl:eq:tausigma}
  \mu_i=\tau_i\quad(1\leq i\leq n-2),\qquad\tau_{\sigma,0}=\tau_0.
\end{equation}
We put $\gamma^{(i)}:=\tau_\sigma-\beta_i$, i.\,e.\ $\gamma^{(i)}=(\mu-e_i;\tau_0+1)$ for $i\leq q$ and
$\gamma^{(i)}=(\mu+e_i;\tau_0-1)$ for $q<i\leq n-2$.

{\em Orientation.} By \cite[Theorems 4.11, 4.14]{ghl} we have $\pi'(q)_{\vz}=\pi^\vee_{A'(q)}=\sigma^\vee$.
The definite places confirm this: there $H'_v$ is compact and the integrand of $I_v$ carries no complex
conjugation \cite[p.~30]{ghl}, so by Schur orthogonality $I_v\not\equiv0$ if and only if
$\pi_v'^\vee\subset\pi_v|_{H'_v}$, which holds by interlacing for $\pi'_v$ the dual of the representation
with infinitesimal character $A_v-\frac12$.

\subsection{The integral}\label{ghl:sec:integral}

The local integral of \cite[p.~30]{ghl} at $\vz$ is
$$
  I(f,f')=\int_{H'}\langle\pi(h)f,f\rangle\langle\pi'(h)f',f'\rangle\,dh',\qquad\pi'=\pi'(q)_{\vz},
$$
without complex conjugation. Let $\bar\pi'$ be the conjugate representation on the conjugate Hilbert
space, where $\langle\bar x,\bar y\rangle:=\overline{\langle x,y\rangle}$. Then $\bar\pi'\cong\pi'^\vee\cong
\sigma$, and with $u:=\bar f'\in\hatsigma$ we have
$$
  \langle\pi'(h)f',f'\rangle=\overline{\langle\sigma(h)u,u\rangle},
$$
so
\begin{equation}\label{ghl:eq:I}
  I(f,f')=\int_{H'}\langle\pi(h)f,f\rangle\,\overline{\langle\sigma(h)u,u\rangle}\,dh',\qquad
  I^*(f,f')=\frac{I(f,f')}{\|f\|^2\|u\|^2}.
\end{equation}
The minimal $K'$-type of $\pi'$ is $\overline{\tau_\sigma}$, so $u$ ranges over $\tau_\sigma$. The
formal degree $d(\sigma)$ is defined by Schur orthogonality,
\begin{equation}\label{ghl:eq:schur}
  \int_{H'}\langle\sigma(h)a,b\rangle\,\overline{\langle\sigma(h)c,e\rangle}\,dh'
  =d(\sigma)^{-1}\langle a,c\rangle\overline{\langle b,e\rangle}\qquad(a,b,c,e\in\hatsigma).
\end{equation}

\subsection{The branching coefficient}\label{ghl:sec:R}

Let $j\colon\hatsigma\to\hatpi$ be an isometric $H'$-embedding, and let $P_\tau$ be the orthogonal projection of
$\hatpi$ onto $\tau$. Since $j$ is $K'$-equivariant, $P_\tau\circ j|_{\tau_\sigma}\colon\tau_\sigma\to\tau$
is a $K'$-map, and $\Hom_{K'}(\tau_\sigma,\tau)$ is one-dimensional by the multiplicity free branching from $K$ to
$K'$, since $\mu$ interlaces $(\tau_1,\dots,\tau_{n-1})$ and $\tau_{\sigma,0}=\tau_0$. Hence
$P_\tau\circ j|_{\tau_\sigma}$ is a scalar multiple of an isometry, and if it is non-zero, then
$$
  s(q):=\frac{\|ju\|^2}{\|P_\tau ju\|^2}=\frac{\|u\|^2}{\|P_\tau ju\|^2}\qquad(0\neq u\in\tau_\sigma)
$$
does not depend on $u$. For $f\in\tau$ and $u\in\tau_\sigma$ with $w:=P_\tau ju\neq0$ we put
$$
  R^\circ(f,f'):=\frac{|\langle f,w\rangle|^2}{\|f\|^2\|w\|^2},\qquad u=\bar f'.
$$
It depends only on the lines $\CC f$ and $\CC u$, and not on the choice of $j$, which is unique up to
a phase by Lemma \ref{ghl:lem:P} (i).

\subsection{Measures}\label{ghl:sec:measures}

On $H'=K'A^+K'$ we write
$$
  dh'=dk_1\,J'(t)\,dt\,dk_2,\qquad J'(t)=c\,\sinh^{2n-5}t\cosh t,\qquad c>0,
$$
where $dk$ is a Haar measure on $K'$. This is the integration formula for the Cartan decomposition. The
measure $dh'$ determines an invariant top form on $\lieh'$, hence on the compact form $\liek'\oplus
i\liep'$ of $\lieh'_\CC$, and so a Haar measure on the compact dual $\U(n-1)$ of $H'$. We write
$\volc=\volc_{dh'}(\U(n-1))$ for its total mass. Harish-Chandra's formula for the formal degree reads
\begin{equation}\label{ghl:eq:fdeg}
  d(\sigma)=\frac{\dim F_{A'(q)}}{\volc}.
\end{equation}

\begin{table}[ht]
  \caption{Dictionary between \cite{ghl} and this paper (at $\vz$).}\label{ghl:tab:dictionary}
  \begin{tabular}{ll}
    \cite{ghl} & this paper\\\hline
    $H_{\vz}=\U(V)_{\vz}$, $H'_{\vz}=\U(V')_{\vz}$ & $H=\U(n-1,1)\supset H'=\U(n-2,1)$, $H'$ fixing $e_{n-1}$\\
    $A_{\lambda,\vz}=(A_1>\dots>A_n)$ & $A$; $L_j=A_j-\frac{n-1}2$\\
    $\pi(q)_{\vz}=\pi_{\lambda,q}$, parameter $A(q)$ & $\pi=\pi_{A(q)}$, minimal $K$-type $\tau$ \eqref{ghl:eq:tau}\\
    $\pi'(q)_{\vz}=\pi^\vee_{A'(q)}$ & $\sigma^\vee$, $\sigma=\pi_{A'(q)}\cong\bar\pi'(q)_{\vz}$\\
    $f_{\vz}$, $f'_{\vz}$ (minimal $K$-types) & $f\in\tau$, $u=\bar f'\in\tau_\sigma$ \eqref{ghl:eq:tausigma}\\
    $c_{f}c_{f'}$, $I_{\vz}$, $I^*_{\vz}$ & \eqref{ghl:eq:I}\\
    $dh'_{\vz}$ & $dh'$; $\volc$, $d(\sigma)$ \eqref{ghl:eq:fdeg}\\
    --- & $j$, $s(q)$, $R^\circ$ (section \ref{ghl:sec:R}); $C_\infty(q)=I^*/R^\circ$
  \end{tabular}
\end{table}
 \section{Matrix coefficients restricted to \texorpdfstring{$H'$}{H'}}\label{ghl:sec:lemmaP}

In this section $\pi$ is any tempered irreducible unitary representation of $H$ and $\sigma$ any
discrete series representation of $H'$. We write $\rho:=\hatpi|_{H'}$ and $\lambda$ for the left
regular representation of $H'$ on $L^2(H')$, $(\lambda(g)F)(h)=F(g^{-1}h)$.

We prove the orthogonality relation (Lemma \ref{ghl:lem:P}).

\subsection{The $L^2$ estimate}\label{ghl:sec:L2}

\begin{lemma}\label{ghl:lem:L2}
  There is $X\in U(\liek_\CC)^K$ such that for $f_1,f_2\in\pi^\infty$ the function
  $h\mapsto\langle\pi(h)f_1,f_2\rangle$ lies in $L^2(H')$, with norm at most
  $C\|\pi(X)f_1\|\|\pi(X)f_2\|$, $C$ independent of $f_1,f_2$.
\end{lemma}

\begin{proof}
  The representation $\pi$ is weakly contained in the regular representation of $H$. By
  \cite[Theorem 1.2, (1.2)]{sun09} there is $X\in U(\liek_\CC)^K$, depending on $K$ only, with
  $$
    |\langle\pi(g)f_1,f_2\rangle|\leq\|\pi(X)f_1\|\|\pi(X)f_2\|\,\Xi_H(g),
  $$
  where $\Xi_H$ is Harish-Chandra's spherical function. Since $H$ has real rank one with root
  multiplicities $2(n-2)$ and $1$,
  $\Xi_H$ is the Jacobi function with $\alpha=n-2$, $\beta=0$, $\lambda=0$:
  $$
    \Xi_H(a_t)={}_2F_1\big(\tfrac{n-1}2,\tfrac{n-1}2;n-1;-\sinh^2t\big)
    =\cosh^{-(n-1)}t\cdot{}_2F_1\big(\tfrac{n-1}2,\tfrac{n-1}2;n-1;\tanh^2t\big),
  $$
  with the second expression given by Pfaff's transformation. The last series has positive coefficients,
  constant term $1$ and parameter excess $0$, hence behaves like a positive multiple of $-\log(1-x)$
  as $x\to1$, and $-\log(1-\tanh^2t)=2\log\cosh t$. Therefore
  $$
    \Xi_H(a_t)\asymp(1+t)e^{-(n-1)t}
  $$
  for $t\geq0$. Since $\Xi_H$ is bi-invariant under $K\supset K'$, the integration formula of section
  \ref{ghl:sec:measures} gives
  $$
    \int_{H'}\Xi_H(h)^2\,dh'\leq C'\int_0^\infty(1+t)^2e^{-2(n-1)t}e^{2(n-2)t}\,dt<\infty.\qedhere
  $$
\end{proof}

\begin{remark}\label{ghl:rmk:L1}
  By the lower bound, $\Xi_H|_{H'}\notin L^1(H')$ for $n\geq3$, so \cite[Corollary 1.3]{sun09} does not
  apply directly. The argument below goes through $L^2$. Absolute convergence of $I$ is also
  \cite[Proposition 2.1]{nharris}.
\end{remark}

By Lemma \ref{ghl:lem:L2} and Cauchy--Schwarz the integrals \eqref{ghl:eq:I} converge absolutely for
$f\in\pi^\infty$, $u\in\hatsigma$, and define a continuous form on $\pi^\infty\times\sigma^\infty$.

\subsection{The $\sigma$-free part}\label{ghl:sec:free}

To isolate the $\sigma$-contribution to \eqref{ghl:eq:I} without decomposing $\rho$, we show that a vector
orthogonal to all copies of $\hatsigma$ whose matrix coefficient is square integrable on $H'$ yields zero integral (Lemma \ref{ghl:lem:free}).

Let $\rho[\sigma]\subset\rho$ be the closed span of the images of all bounded $H'$-maps
$\hatsigma\to\rho$, and $\rho':=\rho[\sigma]^\perp$. Both are closed and $H'$-stable, and
$\Hom_{H'}(\hatsigma,\rho')=0$, since the image of such a map lies in $\rho[\sigma]\cap\rho'=0$.

\begin{lemma}\label{ghl:lem:closed}
  Let $x\in\rho$ with $h\mapsto\langle\rho(h)x,x\rangle$ in $L^2(H')$, and let $V$ be the closed
  $H'$-span of $x$. Then there is an isometric $H'$-embedding $U\colon V\to L^2(H')$, and
  $Tz(h):=\langle z,\rho(h)x\rangle$ satisfies $T=U|T|$ on its domain $D(T)=\{z\in V:Tz\in L^2(H')\}$.
\end{lemma}

\begin{proof}
  Since $T(\rho(g)z)=\lambda(g)Tz$, the domain $D(T)$ is $H'$-stable. It contains $x$, hence is dense
  in $V$. The operator $T$ is closed: if $z_k\to z$ and $Tz_k\to F$ in $L^2$, then $Tz_k\to Tz$
  pointwise and a subsequence converges to $F$ almost everywhere. Furthermore $\ker T=0$, since $Tz=0$
  means $z\perp\rho(H')x$. Let $T=U|T|$
  be the polar decomposition. Then $U$ is a partial isometry with initial space $(\ker T)^\perp=V$. For
  $g\in H'$, $\lambda(g)T\rho(g)^{-1}=T$, and
  $\lambda(g)T\rho(g)^{-1}=(\lambda(g)U\rho(g)^{-1})(\rho(g)|T|\rho(g)^{-1})$ is again a polar
  decomposition: $\rho(g)|T|\rho(g)^{-1}$ is positive self-adjoint with square
  $\rho(g)T^*T\rho(g)^{-1}=T^*T$ (as $T\rho(g)^{-1}=\lambda(g)^{-1}T$, domains included, and
  $(T\rho(g)^{-1})^*=\rho(g)T^*$, $(\lambda(g)^{-1}T)^*=T^*\lambda(g)$, the operators $\rho(g)$ and $\lambda(g)$
  being unitary), and
  $\lambda(g)U\rho(g)^{-1}$ is a partial isometry with initial space $V$. By uniqueness of the polar
  decomposition, $U$ intertwines $\rho|_V$ and $\lambda$.
\end{proof}

\begin{lemma}\label{ghl:lem:free}
  Let $x\in\rho'$ with $h\mapsto\langle\rho(h)x,x\rangle$ in $L^2(H')$. Then
  $$
    \int_{H'}\langle\rho(h)x,x\rangle\overline{\langle\sigma(h)y,y'\rangle}\,dh'=0
  $$
  for all $y,y'\in\hatsigma$.
\end{lemma}

\begin{proof}
  Let $T$, $U$, $V\subset\rho'$ be as in Lemma \ref{ghl:lem:closed}. Put
  $\Phi_{y'}w(h):=\langle w,\sigma(h)y'\rangle$. By \eqref{ghl:eq:schur} this is a bounded $H'$-map
  $\hatsigma\to L^2(H')$. Since $\langle\rho(h)x,x\rangle=\overline{Tx(h)}$ and
  $\overline{\langle\sigma(h)y,y'\rangle}=\Phi_yy'(h)$, the integral is $\langle\Phi_yy',Tx\rangle_{L^2}$.
  Let $E$ be the orthogonal projection onto the closed span $L^2[\sigma]$ of the images of all
  $H'$-maps $\hatsigma\to L^2(H')$. It commutes with $\lambda$, and $\Phi_yy'\in L^2[\sigma]$. For every
  $H'$-map $\iota\colon\hatsigma\to L^2[\sigma]$, $(EU)^*\circ\iota\in\Hom_{H'}(\hatsigma,V)=0$. These
  images span a dense subspace of $L^2[\sigma]$, and the bounded operator $(EU)^*=U^*E$ vanishes on
  $L^2[\sigma]^\perp$, so $(EU)^*=0$ and $EU=0$. Hence $Tx=U|T|x\in\ker E$, and
  $$
    \langle\Phi_yy',Tx\rangle=\langle E\Phi_yy',Tx\rangle=\langle\Phi_yy',ETx\rangle=0.
  $$
\end{proof}

\begin{remark}\label{ghl:rmk:discrete}
  The same argument shows that every irreducible closed subrepresentation $\vartheta$ of $\rho$ is a discrete
  series of $H'$: if $P_\vartheta x\neq0$ for a $K$-finite $x$, $P_\vartheta$ the orthogonal projection onto
  $\vartheta$, then $P_\vartheta|_V$ composed with the adjoint and
  with $U$ gives a non-zero $H'$-map $\vartheta\to L^2(H')$.
  Theorem \ref{ghl:thm:J} will show that for $\pi=\pi_{A(q)}$ the representation $\sigma=\pi_{A'(q)}$ occurs in $\rho$ if $A_q-A_{q+1}\geq2$.
\end{remark}

\subsection{The orthogonality relation}\label{ghl:sec:lemmaPstatement}

\begin{lemma}[Orthogonality relation]\label{ghl:lem:P}
  Let $\pi$ be a tempered irreducible unitary representation of $H$ and $\sigma$ a discrete series of $H'$.
  \begin{itemize}
    \item[(i)] $\dim\Hom_{H'}(\hatsigma,\hatpi|_{H'})\leq1$.
    \item[(ii)] If this space is $0$, then $I\equiv0$ on $\pi^\infty\times\hatsigma$.
    \item[(iii)] If it is spanned by an isometric embedding $j$, then for $f_1,f_2\in\pi^\infty$ and
      $u_1,u_2\in\hatsigma$
      $$
        \int_{H'}\langle\pi(h)f_1,f_2\rangle\overline{\langle\sigma(h)u_1,u_2\rangle}\,dh'
        =d(\sigma)^{-1}\langle j^*f_1,u_1\rangle\overline{\langle j^*f_2,u_2\rangle};
      $$
      in particular
      $$
        I(f,f')=d(\sigma)^{-1}|\langle f,ju\rangle|^2
      $$
      for $u=\bar f'$.
  \end{itemize}
  All integrals converge absolutely.
\end{lemma}

\begin{proof}
  For (i), let $j$ be a bounded $H'$-map. Its adjoint $j^*$ maps $\pi^\infty$ continuously to
  $\sigma^\infty$, since for $Y\in U(\lieh'_\CC)$ we have
  $$
    \|\sigma(Y)j^*x\|=\|j^*\pi(Y)x\|\leq\|j^*\|\|\pi(Y)x\|.
  $$
  Therefore $(x,\bar w)\mapsto\langle j^*x,w\rangle$ is a continuous $H'$-invariant bilinear form on
  $\pi^\infty\times\overline{\sigma^\infty}=\pi^\infty\times\pi'^\infty$, and $j$ is determined by it.
  Both $\pi^\infty$ and $\pi'^\infty$ are irreducible Casselman--Wallach representations, and by
  \cite[Theorem B]{sunzhu} for $\U(1,n-1)\supset\U(1,n-2)$ (replace the Hermitian form by its negative)
  the space of such forms has dimension at most one.

  For (ii) and (iii), write $f=f_\sigma+f_\perp$ with $f_\sigma\in\rho[\sigma]$ and
  $f_\perp\in\rho'$. Both summands are
  $H'$-stable and orthogonal, so the cross terms $\langle\rho(h)f_\sigma,f_\perp\rangle$ and
  $\langle\rho(h)f_\perp,f_\sigma\rangle$ vanish, and $\langle\pi(h)f,f\rangle=F_1(h)+F_2(h)$ with
  $F_1(h)=\langle\rho(h)f_\sigma,f_\sigma\rangle$, $F_2(h)=\langle\rho(h)f_\perp,f_\perp\rangle$. By
  (i), $\rho[\sigma]=j(\hatsigma)$ and $f_\sigma=jj^*f$ if $j$ exists, and $f_\sigma=0$ otherwise, so $F_1$ is a
  matrix coefficient of $\sigma$, in $L^2(H')$. By Lemma \ref{ghl:lem:L2}, $F_1+F_2\in L^2$, hence
  $F_2\in L^2$, and Lemma \ref{ghl:lem:free} gives $\int F_2\overline{\langle\sigma(h)u,u\rangle}=0$.
  With \eqref{ghl:eq:schur} we obtain
  $$
    \int F_1\overline{\langle\sigma(h)u,u\rangle}=d(\sigma)^{-1}|\langle j^*f,u\rangle|^2.
  $$
  This gives the diagonal case, and the general case follows by polarization in $f$ and in $u$,
  since both sides are sesquilinear in $(f_1,f_2)$ and in $(u_1,u_2)$.
\end{proof}
 \section{The Gelfand--Tsetlin model}\label{ghl:sec:gt}

In this section no gap condition is assumed: $0\leq q\leq n-2$ and $L_1>\dots>L_n$ are arbitrary.

We construct an explicit $(\gl_n,K)$-module $V_{\Rq}$ on Gelfand--Tsetlin patterns over a region $\Rq$ and
identify it with $\pi_K$ (Theorem \ref{ghl:thm:identification}), using the parameters and \eqref{ghl:eq:adjoint}
of section \ref{ghl:sec:setting}.

\subsection{The region $\Rq$ and the patterns}\label{ghl:sec:region}

\begin{figure}[ht]
  \centering
  \begin{tikzpicture}[x=-0.72cm,y=1cm,>=stealth,
      pt/.style={circle,fill=black,inner sep=0pt,minimum size=3.2pt},
      bl/.style={circle,draw=black,inner sep=0pt,minimum size=7pt}]
    \fill[black!12] (3.5,1.45) rectangle (-1.5,2.15);
    \foreach \l in {6,3,0,-2,-4}
      \draw[black!45,dashed] (\l+0.5,-0.95) -- (\l+0.5,2.35);
    \draw (10.6,0) -- (-7.6,0);
    \foreach \v in {-7,...,10} \draw (\v,0.06) -- (\v,-0.06);
    \foreach \l/\j in {6/1,3/2,0/3,-2/4,-4/5} {
      \draw[thick] (\l,0.12) -- (\l,-0.12);
      \node[below] at (\l,-0.1) {$\scriptstyle L_{\j}$};}
    \foreach \c/\r in {8.6/0,5/1,2/2,-0.5/3,-2.5/4,-5.8/5}
      \node at (\c,-0.75) {$\scriptstyle I_{\r}$};
    \node[anchor=east] at (11.0,0.8) {$I_1\times\dots\times I_4$};
    \foreach \a/\b/\i in {6/4/1,3/1/2,0/-1/3,-2/-3/4} {
      \draw[thick] (\a,0.8) -- (\b,0.8);
      \foreach \v in {\a,...,\b} \node[pt] at (\v,0.8) {};
      \node[above] at ({(\a+\b)/2},0.82) {$\scriptstyle x_{\i}$};}
    \node[anchor=east] at (11.0,1.8) {$\Rq$};
    \draw[thick,->] (7,1.8) -- (10.4,1.8);
    \foreach \v in {7,...,10} \node[pt] at (\v,1.8) {};
    \node[above] at (8.5,1.82) {$\scriptstyle x_1$};
    \foreach \a/\b/\i in {6/4/2,-2/-3/3} {
      \draw[thick] (\a,1.8) -- (\b,1.8);
      \foreach \v in {\a,...,\b} \node[pt] at (\v,1.8) {};
      \node[above] at ({(\a+\b)/2},1.82) {$\scriptstyle x_{\i}$};}
    \draw[thick,->] (-4,1.8) -- (-7.4,1.8);
    \foreach \v in {-4,...,-7} \node[pt] at (\v,1.8) {};
    \node[above] at (-5.5,1.82) {$\scriptstyle x_4$};
    \foreach \v in {7,4,-2,-4} \node[bl] at (\v,1.8) {};
  \end{tikzpicture}
  \caption{The region $\Rq$ \eqref{ghl:eq:Rq} for $n=5$, $q=2$ and $L=(6,3,0,-2,-4)$, i.\,e.\
    $A=(8,5,2,0,-2)$, with the integers decreasing from left to right. The lower row shows the top rows
    $x_i\in I_i$ of the finite-dimensional representation with highest weight $(L_j+j-1)_j$, the upper row
    the region $\Rq$, which omits $I_q\cup I_{q+1}=[L_4+1,L_2]$ (shaded). The circled points are
    $x(\tau)=(7,4,-2,-4)$.}\label{ghl:fig:Rq}
\end{figure}
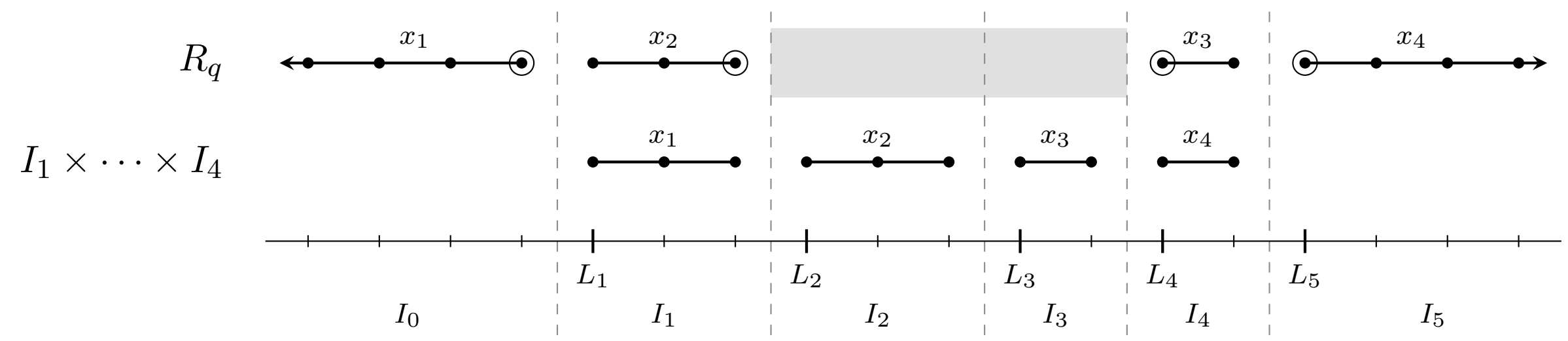

For a Gelfand--Tsetlin pattern $P=(m_{k,i})_{1\leq i\leq k\leq n-1}$ put $l_{k,i}:=m_{k,i}-i+1$ and
$l_{n,j}:=L_j$. On the top row write $m_i:=m_{n-1,i}$, $x_i:=l_{n-1,i}$. The integers $L_1>\dots>L_n$
cut $\ZZ$ into the non-empty intervals $I_0=[L_1+1,\infty)$, $I_r=[L_{r+1}+1,L_r]$ $(1\leq r\leq n-1)$,
$I_n=(-\infty,L_n]$. Let
\begin{equation}\label{ghl:eq:Rq}
  \Rq:\quad x_i\in I_{i-1}=[L_i+1,L_{i-1}]\ (i\leq q),\qquad
  x_i\in I_{i+1}=[L_{i+2}+1,L_{i+1}]\ (q<i\leq n-1),
\end{equation}
where $L_0:=+\infty$ and $L_{n+1}:=-\infty$. It omits the two intervals $I_q$, $I_{q+1}$ adjacent to $L_{q+1}$ (the
top rows of the finite-dimensional representation with highest weight $(L_j+j-1)_j$ fill
$I_1\times\dots\times I_{n-1}$). The region $\Rq$ is a non-empty product of integer intervals,
$x_1>\dots>x_{n-1}$ on $\Rq$, and $\tau\in\Rq$ by \eqref{ghl:eq:tau}:
$x_i(\tau)=L_i+1$ is the lower end for $i\leq q$, and $x_i(\tau)=L_{i+1}$ is the upper
end for $i>q$ (Figure \ref{ghl:fig:Rq}). Let $\cP(\Rq)$ be the set of patterns with top row in
$\Rq$ and $m_{k+1,i}\geq m_{k,i}\geq m_{k+1,i+1}$ $(k\leq n-2)$, i.\,e.\
\begin{equation}\label{ghl:eq:interlace}
  l_{k+1,i+1}+1\leq l_{k,i}\leq l_{k+1,i}.
\end{equation}
In every row, $l_{k,1}>l_{k,2}>\dots>l_{k,k}$.

\subsection{The operators}\label{ghl:sec:operators}

Put $g_i:=i$ $(i\leq q)$, $g_i:=i+2$ $(q<i\leq n-2)$, $\gamma_i(t):=t-1-L_{g_i}$ $(i\leq n-2)$,
$\gamma_{n-1}(t):=1$, and $s_k(P):=\sum_{i\leq k}m_{k,i}$, $s_0:=0$, $s_n:=\mathrm{tot}$. No $g_{n-1}$ is
defined: for $i=n-1$ a product over $j\neq g_i$ is over all $j=1,\dots,n$. On the vector space
$V_{\Rq}$ with basis $\xit_P$, $P\in\cP(\Rq)$, let
$$
  H_a\xit_P=(s_a-s_{a-1})(P)\,\xit_P,\qquad
  E_k\xit_P=\sum_i\tilde A^+_{k,i}(P)\,\xit_{P+\delta_{k,i}},\qquad
  F_k\xit_P=\sum_i\tilde A^-_{k,i}(P)\,\xit_{P-\delta_{k,i}},
$$
with the sums taken over those $i$ for which $P\pm\delta_{k,i}\in\cP(\Rq)$, where for $k\leq n-2$
\begin{equation}\label{ghl:eq:GTcoeff}
  \tilde A^+_{k,i}=-\frac{\prod_{j=1}^{k+1}(l_{k,i}-l_{k+1,j})}{\prod_{s\neq i}(l_{k,i}-l_{k,s})},\qquad
  \tilde A^-_{k,i}=\frac{\prod_{j=1}^{k-1}(l_{k,i}-l_{k-1,j})}{\prod_{s\neq i}(l_{k,i}-l_{k,s})},
\end{equation}
and for $k=n-1$
\begin{equation}\label{ghl:eq:GTtop}
  \tilde A^+_{n-1,i}=-\frac{\prod_{j\neq g_i}(x_i-L_j)}{\prod_{s\neq i}(x_i-x_s)},\qquad
  \tilde A^-_{n-1,i}=\frac{\gamma_i(x_i)\prod_{j=1}^{n-2}(x_i-l_{n-2,j})}{\prod_{s\neq i}(x_i-x_s)}.
\end{equation}
These are the non-normalized Gelfand--Tsetlin coefficients, with the top row $L$ and with the
boundary factor $x_i-1-L_{g_i}$ moved from the raising to the lowering coefficient, so that every wall
of $\Rq$ is a zero of the coefficient that would leave it. Put $E_{aa}:=H_a$, $E_{k,k+1}:=E_k$,
$E_{k+1,k}:=F_k$, and recursively $E_{ab}:=[E_{a,b-1},E_{b-1,b}]$ $(b>a+1)$,
$E_{ab}:=[E_{a,a-1},E_{a-1,b}]$ $(a>b+1)$. All denominators are non-zero on $\cP(\Rq)$.

\subsection{The invariant form}\label{ghl:sec:form}

We define a diagonal Hermitian form on $V_{\Rq}$ and establish its positivity and adjointness for the
generators of section \ref{ghl:sec:operators} before verifying the Lie algebra relations.

Let
\begin{gather*}
  \Delta_k(P):=\prod_{a<b\leq k}(l_{k,a}-l_{k,b}),\\
  G_{k-1,k}(P):=\frac{\prod_{1\leq a\leq b\leq k-1}(l_{k,a}-l_{k-1,b})!}
  {\prod_{1\leq a<b\leq k}(l_{k-1,a}-l_{k,b}-1)!}\quad(2\leq k\leq n-1),
\end{gather*}
$W_q(x):=\prod_{i\leq n-1}w_i(x_i)$ with
$$
  w_i=\frac{(x_i-L_i-1)!\prod_{j<i}(L_j-x_i)!}{\prod_{j>i}(x_i-L_j-1)!}\ (i\leq q),\qquad
  w_i=\frac{(x_i-L_{i+2}-1)!\prod_{j\leq i+1}(L_j-x_i)!}{\prod_{j\geq i+3}(x_i-L_j-1)!}
$$
for $q<i\leq n-2$,
$$
  w_{n-1}=\prod_{j=1}^n(L_j-x_{n-1})!,
$$
and
\begin{equation}\label{ghl:eq:NP}
  N(P):=W_q(x)\prod_{k=1}^{n-1}\Delta_k(P)^{-1}\prod_{k=2}^{n-1}G_{k-1,k}(P),\qquad
  \langle\xit_P,\xit_Q\rangle:=\delta_{PQ}N(P).
\end{equation}

\begin{theorem}\label{ghl:thm:form}
  \noindent {\rm(1)} All factorial arguments are non-negative integers and $\Delta_k(P)>0$, so $N(P)\in\QQ_{>0}$
  and the form is positive definite.

  \noindent {\rm(2)} $H_a^*=H_a$ and
  $$
    \langle E_k\xit_P,\xit_Q\rangle=\epsilon_k\langle\xit_P,F_k\xit_Q\rangle
  $$
  with $\epsilon_k=+1$
  $(k\leq n-2)$, $\epsilon_{n-1}=-1$. On every edge $P,P+\delta_{k,i}\in\cP(\Rq)$ $(1\leq k\leq n-1)$ one has
  $$
    \epsilon_k\tilde A^+_{k,i}(P)>0,\qquad\tilde A^-_{k,i}(P+\delta_{k,i})>0.
  $$
\end{theorem}

\begin{proof}
  (1) Since the rows decrease strictly, $\Delta_k>0$. In $G_{k-1,k}$, \eqref{ghl:eq:interlace} gives
  $l_{k,a}\geq l_{k,b}\geq l_{k-1,b}$ $(a\leq b\leq k-1)$ and $l_{k-1,a}\geq l_{k,a+1}+1\geq l_{k,b}+1$
  $(a<b)$. In $w_i$, for $i\leq q$, $x_i\in[L_i+1,L_{i-1}]$ gives $x_i-L_i-1\geq0$, $L_j-x_i\geq0$ $(j<i)$
  and $x_i-L_j-1\geq L_i-L_{i+1}\geq1$ $(j>i)$. The case $q<i\leq n-2$ is analogous, and $x_{n-1}\leq L_n$.

  (2) It suffices to show, for every edge, with $y:=l_{k,i}$,
  \begin{equation}\label{ghl:eq:edgeratio}
    \frac{N(P+\delta_{k,i})}{N(P)}=\epsilon_k\frac{\tilde A^-_{k,i}(P+\delta_{k,i})}{\tilde A^+_{k,i}(P)}.
  \end{equation}
  Validity of $P+\delta_{k,i}$ gives $y+1\leq l_{k+1,i}$ $(k\leq n-2)$, $y+1\leq l_{k-1,i-1}-1$ $(i\geq2)$ and
  $l_{k,s}\geq y+2$ $(s<i)$. The three factors change as follows:
  \begin{gather*}
    \frac{\Delta_k(P)}{\Delta_k(P+\delta_{k,i})}=\prod_{s\neq i}\frac{y-l_{k,s}}{y+1-l_{k,s}}>0,\\
    \frac{G_{k-1,k}(P+\delta_{k,i})}{G_{k-1,k}(P)}=(-1)^{i-1}\prod_{j\leq k-1}(y+1-l_{k-1,j})>0,
  \end{gather*}
  and for $k\leq n-2$,
  $$
    \frac{G_{k,k+1}(P+\delta_{k,i})}{G_{k,k+1}(P)}=\frac{(-1)^i}{\prod_{j\leq k+1}(y-l_{k+1,j})}>0.
  $$
  For $k\leq n-2$ the product of the three ratios is $\tilde A^-_{k,i}(P+\delta_{k,i})/\tilde A^+_{k,i}(P)$,
  and counting negative factors gives the two signs. For $k=n-1$, $x_i$ and $x_i+1$ lie in the same
  interval of \eqref{ghl:eq:Rq}. The number of $j\neq g_i$ with $x_i<L_j$ is $i-1$, $i+1$ or $n$ (for
  $i\leq q$, $q<i\leq n-2$, $i=n-1$), always $\equiv i-1\bmod 2$, and $\gamma_i(x_i+1)>0$. Under
  $x_i\mapsto x_i+1$ the factor $(x_i-L_{g_i}-1)!$ of $w_i$ $(i\leq n-2)$ is multiplied by $\gamma_i(x_i+1)$,
  and for $j\neq g_i$ each factor $(L_j-x_i)!$ is divided by $L_j-x_i$ and each factor $(x_i-L_j-1)!$ of the
  denominator is multiplied by $x_i-L_j$, so
  $$
    \frac{w_i(x_i+1)}{w_i(x_i)}=\frac{\gamma_i(x_i+1)}{\prod_{j\neq g_i}|x_i-L_j|}
    =\frac{(-1)^{i-1}\gamma_i(x_i+1)}{\prod_{j\neq g_i}(x_i-L_j)}>0.
  $$
  With the ratios of $\Delta_{n-1}$ and $G_{n-2,n-1}$ the signs $(-1)^{i-1}$ cancel and
  \eqref{ghl:eq:edgeratio} holds with $\epsilon_{n-1}=-1$.
\end{proof}

\subsection{The relations}\label{ghl:sec:relations}

\begin{lemma}\label{ghl:lem:serre}
  Let $H_a$ $(a\leq n)$, $E_k$, $F_k$ $(k\leq n-1)$ be operators on a complex vector space with
  
  \noindent {\rm(CS1)} $[H_a,H_b]=0$;

  \noindent {\rm(CS2)} $[H_a,E_k]=(\delta_{a,k}-\delta_{a,k+1})E_k$,
  $[H_a,F_k]=-(\delta_{a,k}-\delta_{a,k+1})F_k$;

  \noindent {\rm(CS3)} $[E_k,F_l]=\delta_{kl}(H_k-H_{k+1})$;
  
  \noindent {\rm(CS4)} $[E_k,E_l]=[F_k,F_l]=0$ for $|k-l|\geq2$;

  \noindent {\rm(CS5)} $[E_k,[E_k,E_{k\pm1}]]=[F_k,[F_k,F_{k\pm1}]]=0$.

  \noindent Then

  \noindent {\rm(1)} the $E_{ab}$ satisfy $[E_{ab},E_{cd}]=\delta_{bc}E_{ad}-\delta_{ad}E_{cb}$;
    
  \noindent {\rm(2)} if there is a Hermitian form with $H_a^*=H_a$ and $E_k^*=\epsilon_kF_k$, $\epsilon$ as in
  Theorem \ref{ghl:thm:form}, then \eqref{ghl:eq:adjoint} holds.
\end{lemma}

\begin{proof}
  \noindent (1) (CS1)--(CS5) are the relations of Serre's presentation of $\slie_n$ with $h_k=H_k-H_{k+1}$
  \cite[\S18.1, \S18.3, Theorem (Serre), p.~99]{humphreys}. That theorem is a statement about the Lie
  algebra presented by these relations, so the operators define a representation of $\slie_n$, with no
  condition on the space. The operator $Z:=\sum_aH_a$ commutes with all $E_k$, $F_k$ by (CS2), and
  $\gl_n=\slie_n\oplus\CC\one$ acts with $E_{aa}\mapsto H_a$.

  \noindent (2) We proceed by induction on $|a-b|$ with $[X,Y]^*=[Y^*,X^*]$, obtaining
  $[E_{a,b-1},E_{b-1,b}]^*=E_{ba}$ for $b\leq n-1$ and
  $[E_{a,n-1},E_{n-1,n}]^*=[-E_{n,n-1},E_{n-1,a}]=-E_{na}$.
\end{proof}

\begin{theorem}\label{ghl:thm:relations}
  The operators of section \ref{ghl:sec:operators} satisfy (CS1)--(CS5) on $V_{\Rq}$. Hence $V_{\Rq}$ is a
  $\gl_n$-module, and the form \eqref{ghl:eq:NP} satisfies \eqref{ghl:eq:adjoint}.
\end{theorem}

\begin{proof}
  Every constraint of $\cP(\Rq)$ is an interval bound on one entry $m_{k,i}$ given by the rows $k\pm1$
  (by $\Rq$ for $k=n-1$). By Theorem \ref{ghl:thm:form} each $F$-relation follows from the corresponding
  $E$-relation by taking adjoints.

  For (CS1) and (CS2), the operator $H_a$ is diagonal, and $E_k$ $(F_k)$ changes $s_k$ by $+1$ $(-1)$
  and no other $s_r$.

  For (CS4) and (CS3) with $|k-l|\geq2$, the coefficients of $E_k$ depend on rows $k,k+1$ only, and
  validity of $P+\delta_{k,i}+\delta_{l,j}$ is equivalent to that of $P+\delta_{k,i}$ and $P+\delta_{l,j}$.

  For (CS3) with $k=l+1$, raising $m_{l+1,i}$ can only violate $m_{l,i-1}\geq m_{l+1,i}$, which no
  lowering in row $l$ repairs, and lowering $m_{l,j}$ can only violate $m_{l,j}\geq m_{l+1,j+1}$, which no
  raising in row $l+1$ repairs, so both intermediate patterns are valid and the two orders give the same
  coefficient. The case
  $l=k+1$ is the adjoint.

  For (CS3) with $k=l$, we show that $[E_k,F_k]=H_k-H_{k+1}$.
  Here $x_r:=l_{k,r}$ denotes the entries of row $k$, not those of the top
  row. For off-diagonal terms $Q=P+\delta_{k,i}-\delta_{k,j}$ $(i\neq j)$ both
  intermediate patterns are valid, and the two products agree except for one pair of denominators,
  where
  $$
    (x_j-x_i)(x_i-x_j+1)=(x_i-x_j)(x_j-x_i-1).
  $$
  For the diagonal term we need
  \begin{equation}\label{ghl:eq:CS3diag}
    \sum_{i:\,P-\delta_{k,i}\ \rm valid}\tilde A^-_{k,i}(P)\tilde A^+_{k,i}(P-\delta_{k,i})
    -\sum_{i:\,P+\delta_{k,i}\ \rm valid}\tilde A^+_{k,i}(P)\tilde A^-_{k,i}(P+\delta_{k,i})
    =(H_k-H_{k+1})(P),
  \end{equation}
  where $(H_k-H_{k+1})(P)$ is the eigenvalue of $H_k-H_{k+1}$ on $\xit_P$. For $k=n-1$ the factor $\gamma_i$ combines with $\prod_{j\neq g_i}(t-L_j)$ to $\prod_{j=1}^n(t-L_j)$. For
  small $\varepsilon\neq0$ put $x_s(\varepsilon):=x_s+s\varepsilon$ and
  $$
    F_\varepsilon(t):=\frac{\prod_{r=1}^{k+1}(t-l_{k+1,r})\prod_{r=1}^{k-1}(t+1-l_{k-1,r})}
    {\prod_{s=1}^k(t-x_s(\varepsilon))(t+1-x_s(\varepsilon))}.
  $$
  Its $2k$ poles are simple, and the residues at $x_i(\varepsilon)-1$ and $x_i(\varepsilon)$ are the
  deformed products of the first and second sum of \eqref{ghl:eq:CS3diag}, with signs $+$ and $-$: the
  other denominator factor with $s=i$ takes the value $t-x_i(\varepsilon)=-1$, resp.\ $t+1-x_i(\varepsilon)=1$,
  at these poles, and $\tilde A^+_{k,i}$ carries a minus sign. As a
  quotient of monic polynomials of degree $2k$, the sum of its residues is the coefficient of $t^{-1}$
  at $\infty$,
  $$
    2\sum_sx_s(\varepsilon)-\sum_rl_{k+1,r}-\sum_rl_{k-1,r}-1,
  $$
  whose limit is
  $(H_k-H_{k+1})(P)$ by $\sum_jl_{r,j}=s_r-r(r-1)/2$. As $\varepsilon\to0$: if $P-\delta_{k,i}$ is valid,
  $x_i-1>x_{i+1}$ and the residue at $x_i(\varepsilon)-1$ tends to the term of \eqref{ghl:eq:CS3diag};
  if not, $x_i$ lies on a lower wall ($x_i=l_{k-1,i}$, or $x_i-1=l_{k+1,i+1}$ for $k\leq n-2$, or
  $x_i-1=L_{g_i}$ for $k=n-1$). If moreover $x_i-1\neq x_{i+1}$, a numerator factor vanishes at
  $\varepsilon=0$. If $x_i-1=x_{i+1}$, the denominator has the simple zero $-\varepsilon$, but then
  interlacing forces $x_i=l_{k-1,i}$ and $x_i-1=l_{k+1,i+1}$ for $k\leq n-2$, resp.\ $x_i=l_{n-2,i}$ and
  $x_i-1=L_i$ $(i+1\leq q)$ or $x_i-1=L_{i+2}$ $(i>q)$ for $k=n-1$ ($i=q$ is impossible:
  $x_q\geq L_q+1\geq L_{q+2}+3\geq x_{q+1}+3$), so the numerator vanishes to order two, one factor from each
  of its two products, and the residue tends to $0$. The upper walls are symmetric. For $k=n-1$ the upper
  walls of $\Rq$ are $x_i=L_{i-1}$ $(2\leq i\leq q)$ and $x_i=L_{i+1}$ $(i>q)$, zeros of
  $\prod_{j=1}^n(t-L_j)$ at $t=x_i$.

  For (CS5), fix $k\leq n-2$ and put $x_i:=l_{k,i}$, $y_j:=l_{k+1,j}$.
  Constraints not coupling rows $k$, $k+1$ hold
  on every intermediate pattern.

  (B1) We show $(E_k^2E_{k+1}-2E_kE_{k+1}E_k+E_{k+1}E_k^2)\xit_P=0$ at the target
  $P+\delta_{k,i_1}+\delta_{k,i_2}+\delta_{k+1,j}$. After cancelling common factors, for $i_1=i_2$ we obtain
  $$
    (x-y-1)(x-y)-2(x-y)^2+(x-y)(x-y+1)=0\qquad(x=x_{i_1},\ y=y_j),
  $$
  and for $i_1<i_2$, with $u:=x_{i_1}-y_j$, $v:=x_{i_2}-y_j$ and $u-v\geq2$ (from the interlacing with row
  $k-1$), the six orderings give $-(u-v)^{-1}$ times
  $$
    \frac{(u-1)(v-1)-2u(v-1)+uv}{u-v+1}+\frac{(u-1)(v-1)-2v(u-1)+uv}{u-v-1}=1-1=0.
  $$

  (B2) We show $(E_{k+1}^2E_k-2E_{k+1}E_kE_{k+1}+E_kE_{k+1}^2)\xit_P=0$. For $j_1=j_2$ we obtain
  $$
    (x-y)-2(x-y-1)+(x-y-2)=0\qquad(x=x_i,\ y=y_{j_1}),
  $$
  and for $j_1<j_2$, with $u:=x_i-y_{j_1}$ and
  $v:=x_i-y_{j_2}$, one has $v-u\geq2$ (via row $k+2$, or via \eqref{ghl:eq:Rq} for $k+1=n-1$), and the six
  orderings give $-(v-u)^{-1}$ times
  $$
    \frac{uv-2(u-1)v+(u-1)(v-1)}{v-u+1}+\frac{uv-2u(v-1)+(u-1)(v-1)}{v-u-1}=1-1=0.
  $$
  These identities are for the formal sums over all orderings. By the separation bounds no denominator
  vanishes on any ordering. An ordering through an invalid intermediate pattern violates the interlacing
  of rows $k$, $k+1$ in one of two ways:

  \noindent ($\alpha$) row $k+1$ is raised at $j$ while $m_{k,j-1}=m_{k+1,j}$, i.\,e.\ $x_{j-1}=y_j+1$,
  and the later raising of $m_{k,j-1}$, which must occur, carries the factor $x_{j-1}-(y_j+1)=0$; or

  \noindent ($\beta$) row $k$ is raised at $i$ while $m_{k,i}=m_{k+1,i}$, i.\,e.\ $x_i=y_i$, and this step carries the factor
  $x_i-y_i=0$.

  Therefore every dropped ordering has formal value $0$, and (CS5) follows. Lemma
  \ref{ghl:lem:serre} gives the rest.
\end{proof}

\subsection{$K$-types and irreducibility}\label{ghl:sec:irred}

\begin{proposition}\label{ghl:prop:Ktypes}
  $V_{\Rq}=\bigoplus_{m\in\Rq}V_m$, $V_m$ the span of the $\xit_P$ with top row $m$, and $V_m\cong F_m\boxtimes
  \chi_{\mathrm{tot}-\sum m}$ as a module for $\liek=\gl_{n-1}\oplus\gl_1$. The action integrates to
  $K$, compatibly with $\gl_n$. Hence $V_{\Rq}$ is an admissible multiplicity free $(\gl_n,K)$-module
  with $K$-types $(m;\mathrm{tot}-\sum m)$, $m\in\Rq$.
\end{proposition}

\begin{proof}
  The operators $H_1,\dots,H_{n-1}$ and $E_k,F_k$ $(k\leq n-2)$ preserve $V_m$, and $H_n$ acts by
  $\mathrm{tot}-\sum m$.
  The basis of $V_m$ has $\dim F_m$ elements. The vector with $m_{k,i}=m_i$ has weight $m$ and is killed
  by $E_1,\dots,E_{n-2}$, and $V_m$ is $\U(n-1)$-unitary by Theorem \ref{ghl:thm:form}, hence completely
  reducible, so $V_m\cong F_m$. The weights are integral, so the action integrates to the connected
  group $K$, and the relations give compatibility with $\gl_n$.
\end{proof}

\begin{proposition}\label{ghl:prop:irred}
  $V_{\Rq}$ is irreducible.
\end{proposition}

\begin{proof}
  A non-zero submodule is $\bigoplus_{m\in T}V_m$. If $m,m+e_i\in\Rq$, then $m_{i-1}\geq m_i+1$ for
  $i\geq2$ (for $2\leq i\leq q$, $m_{i-1}\geq\tau_{i-1}\geq m_i+1$; for $i=q+1$,
  $m_q\geq L_q+q>L_{q+2}+q\geq m_{q+1}+1$; for $i\geq q+2$, $m_{i-1}\geq\tau_i\geq m_i+1$), the lower
  bounds by $m\in\Rq$ and the upper ones by $m+e_i\in\Rq$, so $m$ with the $i$-th entry removed interlaces both $m$
  and $m+e_i$ (for $i=1$ there is nothing to check). Let $P$ be the pattern with top row $m$, row $n-2$
  equal to $m$ without its $i$-th entry, and the highest pattern of that weight below it, and
  $Q:=P+\delta_{n-1,i}$. Theorem \ref{ghl:thm:form} (2) gives $\tilde A^+_{n-1,i}(P)\neq0\neq\tilde
  A^-_{n-1,i}(Q)$. Hence $m\in T$ if and only if $m+e_i\in T$, and $\Rq$ is connected by unit steps.
\end{proof}

\subsection{Identification with $\pi_K$}\label{ghl:sec:identification}

We use Schmid's characterization \cite[Theorem (1.3), pp.~536--537]{schmid75}: for $G$ connected,
semisimple, linear, with $\mathrm{rk}\,K=\mathrm{rk}\,G$, and $\pi_\lambda$ the discrete series with
non-singular Harish-Chandra parameter $\lambda\in\Lambda+\rho$, $\Psi=\{\alpha:\langle\lambda,\alpha\rangle>0\}$,
$\pi_\lambda$ is, up to infinitesimal equivalence, the only irreducible quasi-simple representation whose
restriction to $K$

\noindent (a) contains the $K$-type of highest weight $\lambda-\rho_c+\rho_n$, and

\noindent (b) contains no
$K$-type of highest weight $\lambda-\rho_c+\rho_n-A$, $A$ a non-empty sum of roots in $\Psi$, compact or
non-compact.

Highest weights are taken with respect to $\Phi^c\cap\Psi$ \cite[p.~537, footnote 1]{schmid75},
and every irreducible unitary representation is quasi-simple (p.~536 of loc.\ cit.).

\begin{theorem}\label{ghl:thm:identification}
  $V_{\Rq}\cong\pi_K$ as $(\gl_n,K)$-modules, and \eqref{ghl:eq:NP} is, up to a positive scalar, the unique
  invariant Hermitian form on $V_{\Rq}$.
\end{theorem}

\begin{proof}
  {\em Reduction to $G_1=\SU(n-1,1)$.} The group $G_1$ is connected, semisimple, linear and of equal rank,
  and $K_1=K\cap G_1\cong\U(n-1)$ is connected. We have $H=Z\cdot G_1$ and $K=Z\cdot K_1$ with $Z=\U(1)\cdot\one_n$.
  The group $Z$ acts on $V_{\Rq}$ by the character with differential $\sum_aH_a=\mathrm{tot}$, and on $\pi$
  by the same character, since $\sum_iA(q)_i=\sum_jA_j=\mathrm{tot}$. The module $V_{\Rq}$ is an irreducible
  $(\slie_n,K_1)$-module (Proposition \ref{ghl:prop:irred}), and $\pi|_{G_1}$ is the discrete series with
  parameter $\lambda=A(q)$ (restricted to the trace-zero part, which leaves $\Psi$ and the pairings
  unchanged, since roots have trace zero). Furthermore, $K$-types with a fixed central character restrict to
  distinct $K_1$-types. Hence it suffices to verify (a) and (b) for the $K$-types of $V_{\Rq}$, and an
  $(\slie_n,K_1)$-isomorphism $V_{\Rq}\cong\pi_K$ is a $(\gl_n,K)$-isomorphism. Finally, $\Phi^c\cap\Psi$ is
  the standard positive system of $\U(n-1)$, so the indexing of $K$-types agrees with ours.

  {\em Quasi-simplicity.} By Theorems \ref{ghl:thm:form}, \ref{ghl:thm:relations} and Propositions
  \ref{ghl:prop:Ktypes}, \ref{ghl:prop:irred}, $V_{\Rq}$ is irreducible, admissible and unitarizable, hence
  the Harish-Chandra module of an irreducible unitary, hence quasi-simple, representation of $G_1$.

  {\em (a).} We have $\lambda-\rho_c+\rho_n=\tau\in\Rq$.

  {\em (b), the cone argument.} Every $K$-type of $V_{\Rq}$ is $(m;\mathrm{tot}-\sum m)$, $m\in\Rq$, and with
  $\delta_i:=m_i-\tau_i$ (so $\delta_i\geq0$ for $i\leq q$, $\delta_i\leq0$ for $i>q$, by \eqref{ghl:eq:Rq}),
  we have
  $$
    (m;\mathrm{tot}-\textstyle\sum m)=\tau+\sum_ic_i\beta_i,\qquad c_i=\delta_i\ (i\leq q),\quad c_i=-\delta_i\ (i>q),
    \quad c_i\in\ZZ_{\geq0},
  $$
  the $n$-th coordinate being $-\sum_i\delta_i$. If $\tau-A=\tau+\sum_ic_i\beta_i$ with $A=\sum_{\alpha\in\Psi}a_\alpha\alpha$,
  $a_\alpha\in\ZZ_{\geq0}$, then $A+\sum_ic_i\beta_i=0$ is a non-negative combination of roots in $\Psi$,
  which lies in the open half space $\langle\lambda,\cdot\rangle>0$. Pairing with $\lambda$ gives $A=0$.

  By Schmid's theorem $V_{\Rq}\cong\pi_K$. An invariant Hermitian form is a $(\gl_n,K)$-map to the
  Hermitian dual. Since $V_{\Rq}$ is irreducible with multiplicity free $K$-types, these maps form a line,
  and \eqref{ghl:eq:NP} is a positive element.
\end{proof}

\begin{corollary}\label{ghl:cor:ktypes}
  The $K$-types of $\pi=\pi_{A(q)}$ are the $(m;\mathrm{tot}-\sum m)$ with $m\in\Rq$, each with
  multiplicity one.
\end{corollary}

\begin{remark}\label{ghl:rmk:confirmations}
  An irreducible
  $(\lieg,K)$-module is in general not determined by its infinitesimal character and its lowest $K$-type
  (for $\SL(3,\CC)$ and trivial infinitesimal character there are inequivalent modules with
  the same lowest $K$-type), which is why the cone argument for all sums $A$ in (b) is needed. Schmid's
  weaker condition with $A$ a single non-compact root does not suffice in general \cite[p.~537]{schmid75}.
\end{remark}
 \section{The Schmid vector and the term formula}\label{ghl:sec:schmid}

Throughout, $A$ satisfies \eqref{ghl:eq:gap}. By Theorem \ref{ghl:thm:identification} we fix an
isomorphism $\Phi\colon V_{\Rq}\to\pi_K$ of $(\gl_n,K)$-modules. For $m\in\Rq$ interlacing $\mu$ let
$\xi(m,\mu)\in V_m$ be the basis vector $\xit_P$ whose rows $1,\dots,n-2$ form the highest pattern of weight
$\mu$ ($m_{k,i}=\mu_i$ for $k\leq n-2$). It is a $K'$-highest weight vector of weight
$(\mu;\mathrm{tot}-\sum m)$. Similarly, for $m'$ interlacing $\mu-e_i$, $\Qhw(m';\mu-e_i)$ denotes the basis
vector $\xit_P$ with rows $1,\dots,n-2$ the highest pattern of weight $\mu-e_i$. We write $x_i=m_i-i+1$ as in
section \ref{ghl:sec:region} and $y_i:=\mu_i-i+1$ $(i\leq n-2)$, so that $y_i=L_i+1$ $(i\leq q)$ and
$y_i=L_{i+1}$ $(q<i\leq n-2)$.

\subsection{The support}\label{ghl:sec:support}

Put
$$
  a_r:=A_r-A_{q+1}\ (1\leq r\leq q),\qquad B:=A_{q+1}-A_n\geq1,
$$
and for a top row $m$: $d_r:=m_r-\tau_r$ $(1\leq r\leq q)$, $d_0:=\tau_{n-1}-m_{n-1}$, $u_r:=a_r+d_r$. Then
$a_1>\dots>a_q\geq2$ by \eqref{ghl:eq:gap}, and $u_r=x_r-L_{q+1}-1$.

\begin{lemma}\label{ghl:lem:support}
  Let $m\in\Rq$ interlace $\mu$ with $\sum_im_i=\sum_i\tau_i$. Then $m_i=\tau_i$ for $q<i\leq n-2$, and
  the set of such $m$ is
  \begin{equation}\label{ghl:eq:S}
    S:\qquad d_1\geq0,\qquad 0\leq d_r\leq a_{r-1}-a_r-1\ (2\leq r\leq q),\qquad d_0=d_1+\dots+d_q.
  \end{equation}
\end{lemma}

\begin{proof}
  For $q<i\leq n-2$, interlacing gives $m_i\geq\mu_i=\tau_i$ and $\Rq$ gives $m_i\leq\tau_i$. The sum
  condition gives $d_0=\sum d_r$. For $2\leq r\leq q$, interlacing gives $m_r\leq\mu_{r-1}=\tau_{r-1}$,
  i.\,e.\ $d_r\leq\tau_{r-1}-\tau_r=a_{r-1}-a_r-1$, which is also the bound of $\Rq$, and $m_r\geq\tau_r$ is
  $\Rq$. Finally $m_{n-1}\leq\mu_{n-2}$ follows from $m_{n-1}\leq\tau_{n-1}\leq\tau_{n-2}=\mu_{n-2}$.
\end{proof}

\begin{figure}[ht]
  \centering
  \begin{tikzpicture}[>=stealth,
      pt/.style={circle,fill=black,inner sep=0pt,minimum size=3.2pt},
      gpt/.style={circle,draw=black,fill=white,inner sep=0pt,minimum size=3.6pt},
      bl/.style={circle,draw=black,inner sep=0pt,minimum size=7pt}]
    \begin{scope}
      \foreach \e/\i in {{\tau_1+d_1}/0,{\tau_2+d_2}/1,{\tau_3}/2,{\tau_4-d_0}/3}
        \node at (1.6*\i,2.4) {$\e$};
      \foreach \e/\i in {{\tau_1}/0,{\tau_2}/1,{\tau_3}/2}
        \node at (0.8+1.6*\i,1.6) {$\e$};
      \foreach \e/\i in {{\tau_1}/0,{\tau_2}/1}
        \node at (1.6+1.6*\i,0.8) {$\e$};
      \node at (2.4,0) {$\tau_1$};
      \node[anchor=west] at (5.6,2.4) {$m$};
      \node[anchor=west] at (5.6,1.6) {$\mu$};
      \node[anchor=west] at (5.6,0.8) {$(\tau_1,\tau_2)$};
      \node[anchor=west] at (5.6,0) {$(\tau_1)$};
    \end{scope}
    \begin{scope}[xshift=8.9cm,x=0.62cm,y=0.62cm]
      \draw[->] (0,0) -- (7.6,0) node[right] {$d_1$};
      \draw[->] (0,0) -- (0,3.1) node[above] {$d_2$};
      \foreach \k in {0,1,2} \node[left] at (-0.15,\k) {$\scriptstyle\k$};
      \foreach \k in {0,...,6} \node[below] at (\k,-0.1) {$\scriptstyle\k$};
      \foreach \a in {0,...,6} {
        \node[pt] at (\a,0) {}; \node[pt] at (\a,1) {}; \node[gpt] at (\a,2) {};}
      \foreach \k in {0,1,2} \node at (7.05,\k+0.3) {$\cdots$};
      \node[bl] at (0,0) {};
      \node[anchor=south west] at (0.05,0.05) {$\scriptstyle\tau$};
      \node[anchor=west] at (7.5,2) {$\scriptstyle x_2=L_1$};
    \end{scope}
  \end{tikzpicture}
  \caption{The support for $n=5$, $q=2$ and $A=(8,5,2,0,-2)$, so that $\tau=(7,5,0,-1;2)$, $\mu=(7,5,0)$,
    $a=(6,3)$ and $B=4$. Left: the basis vector $\xi(m,\mu)$, rows $4,3,2,1$ from top to bottom, for $m\in S$
    with $d_0=d_1+d_2$. Interlacing with $\mu$ and $\Rq$ freeze $m_3=\tau_3$, and $m_2\leq\mu_1$ gives
    $d_2\leq a_1-a_2-1=2$. Right: the set $S$ \eqref{ghl:eq:S}, with $\tau$ at $d=0$. On the row $d_2=2$
    (open points) $x_2=L_1$, so that $x+e_2\notin\Rq$ (proof of Lemma \ref{ghl:lem:normS}). The rows
    $d_2=\mathrm{const}$ are the inner sums of \eqref{ghl:eq:split} for $p=2$.}\label{ghl:fig:S}
\end{figure}
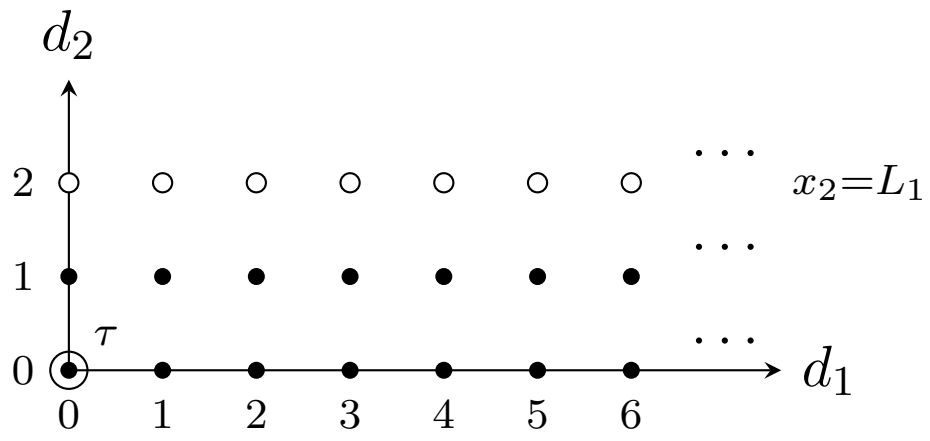

Figure \ref{ghl:fig:S} shows $\xi(m,\mu)$ and $S$ for $n=5$, $q=2$.

For $q=0$ all of $m_1,\dots,m_{n-2}$ are frozen and the sum condition gives $S=\{\tau\}$. {\em From now
on $q\geq1$.} The active indices are $I:=\{1,\dots,q,n-1\}$. For a function $x$ on $I$ put
$$
  \Delta_I(x):=\prod_{a<b\in I}(x_a-x_b),
$$
and
$$
  \Delta_q(u):=\prod_{r<s\leq q}(u_r-u_s).
$$
On $S$,
$u_r\leq a_{r-1}-1<u_{r-1}$, so $\Delta_q(u)>0$.

\subsection{Two lemmas on $\pi$ and $\sigma$}\label{ghl:sec:twolemmas}

We record two lemmas for section \ref{ghl:sec:conditions}: Lemma
\ref{ghl:lem:termwise} projects the $\lieh'_\CC$-action on an $H'$-smooth vector termwise onto a $K$-type,
and Lemma \ref{ghl:lem:Kprime} excludes $\gamma^{(i)}$ from $\sigma$ (also used in Propositions \ref{ghl:prop:Hprime}
and \ref{ghl:prop:Jq0}).

For a $K$-type $\delta$ (resp.\ $K'$-type $\gamma$) let $P_\delta$ (resp.\ $P'_\gamma$) be the orthogonal
projection of $\hatpi$ onto the $\delta$- (resp.\ $\gamma$-)isotypic subspace. Since $K'\subset K$, $P_\delta$
and $P'_\gamma$ commute.

\begin{lemma}\label{ghl:lem:termwise}
  Let $w\in\hatpi$ be $H'$-smooth, $X\in\lieh'_\CC$, $\delta'$ a $K$-type, and $F$ the finite set of
  $K$-types of $\pi$ occurring in $\lieh_\CC\otimes\delta'$. Then
  $$
    P_{\delta'}\pi(X)w=\sum_{\delta\in F}P_{\delta'}\pi(X)P_\delta w,
  $$
  where on the right $\pi(X)$ is the action
  on $\pi_K$.
\end{lemma}

\begin{proof}
  By linearity let $X\in\lieh'$ be real. All $K$-finite vectors are smooth. For $e\in\pi_K$ of type $\delta'$,
  $$
    \langle\pi(X)w,e\rangle=\frac{d}{dt}|_{t=0}\langle w,\pi(e^{-tX})e\rangle=-\langle w,\pi(X)e\rangle,
  $$
  and $\pi(X)e\in\bigoplus_{\delta\in F}\pi_K(\delta)$, so
  $$
    -\langle w,\pi(X)e\rangle=-\sum_{\delta\in F}\langle
    P_\delta w,\pi(X)e\rangle=\sum_{\delta\in F}\langle\pi(X)P_\delta w,e\rangle,
  $$
  by the same identity for the $K$-finite vector $P_\delta w$.
\end{proof}

\begin{lemma}\label{ghl:lem:Kprime}
  For $1\leq i\leq n-2$, $\gamma^{(i)}=\tau_\sigma-\beta_i$ is not a $K'$-type of $\sigma$. More precisely,
  every $K'$-type $\delta$ of $\sigma$ satisfies $\delta-\tau_\sigma\in\NN\Delta^+(\sigma)$.
\end{lemma}

\begin{proof}
  Let $\Lambda=A'(q)$. It is regular and integral and determines $\Delta^+=\Delta^+(\sigma)\supset
  \Delta_n^+(\sigma)$, with compact part $\Delta_c^+$. The group $H'$ is generated by its centre and
  $\SU(n-2,1)$, the centre acts on $\sigma$ by a scalar, $K'$-types with a given central character restrict
  to distinct $K'\cap\SU(n-2,1)$-types, and $\tau_\sigma-\beta_i$ has the central character of $\tau_\sigma$. The
  Weyl group $W_{K'}$, the roots and their pairings with $\Lambda$ are unchanged. Therefore Blattner's formula
  \cite[Theorem (1.3), p.~130]{hechtschmid} applies: if $\delta$ has non-zero multiplicity, then
  $$
    w(\delta+\rho_c)-\tau_\sigma-\rho_c\in\NN\Delta_n^+(\sigma)
  $$
  for some $w\in W_{K'}$. Since $\delta+\rho_c$
  is dominant, $\delta+\rho_c-w(\delta+\rho_c)\in\NN\Delta_c^+$, hence $\delta-\tau_\sigma\in\NN\Delta^+$. For
  $\delta=\tau_\sigma-\beta_i$ this gives $-\beta_i\in\NN\Delta^+$, and pairing with $\Lambda$ gives
  $\langle-\beta_i,\Lambda\rangle<0\leq\langle\delta-\tau_\sigma,\Lambda\rangle$, a contradiction.
\end{proof}

\subsection{The reduced equations}\label{ghl:sec:reduced}

We define linear equations in $V_{\Rq}$ for functions $c\colon S\to\CC$, the reduced equations, relate their left
sides to the $K'$-isotypic projections of $\hatpi$ (Lemma \ref{ghl:lem:coeff}), and compute the coefficients of
$E_{n,i}$ from \eqref{ghl:eq:GTcoeff} and \eqref{ghl:eq:GTtop} (Lemma \ref{ghl:lem:matrixelement}). The reduced
equations come from the Schmid conditions for $ju$ (Proposition \ref{ghl:prop:ju}).

For $1\leq i\leq q$ and $m'\in\Rq$ with $\sum m'=\sum\tau-1$ interlacing $\mu-e_i$, the {\em reduced equation}
$(m',i)$ for a function $c\colon S\to\CC$ is: the coefficient of $\Qhw(m';\mu-e_i)$ in
$$
  E_{n,i}\sum_{\delta\in S}c(\delta)\xi(\delta,\mu)
$$
vanishes. For such $m'$, $(\mu-e_i)_k=\tau_k$ for
$q<k\leq n-2$, so $m'_k=\tau_k$ there, and the $\delta\in S$ with a matrix element to $m'$ are the
$m'+e_r$, $r\in I$. If $\mu_i=\mu_{i+1}$, which happens for $i<q$ exactly when $A_i-A_{i+1}=1$, then $\mu-e_i$ is not
dominant and there is no such $m'$. The weight $\mu-e_q$ is always dominant, since $\mu_q-\mu_{q+1}=L_q-L_{q+2}\geq2$
for $q\leq n-3$, and $\mu_q$ is the last entry of $\mu$ for $q=n-2$.

\begin{lemma}\label{ghl:lem:coeff}
  Let $1\leq i\leq q$, let $m'$ be as in a reduced equation $(m',i)$, and let $x\in V_{m'}$ have $\U(n-2)$-weight
  $\mu-e_i$. Then
  $$
    P'_{\gamma^{(i)}}\Phi(x)=\beta\,\Phi(\Qhw(m';\mu-e_i)),
  $$
  where $\beta$ is the coefficient of $\Qhw(m';\mu-e_i)$ in $x$. In particular $e:=\Phi(\Qhw(m';\mu-e_i))$ is
  non-zero, and $\langle\Phi(x),e\rangle=\beta\|e\|^2$.
\end{lemma}

\begin{proof}
  On $V_{m'}$ the element $E_{nn}$ acts by $\mathrm{tot}-\sum m'=\tau_0+1$. Since $\Phi$ is $K$-equivariant,
  $P'_{\gamma^{(i)}}\Phi(x)=\Phi(P'_{\gamma^{(i)}}x)$, with the isotypic projection of the $K'$-module $V_{m'}$ on the
  right. Basis vectors of $V_{m'}$ with row $n-2$ different from $\mu-e_i$ lie in other $K'$-isotypic components.
  Those with row $n-2$ equal to $\mu-e_i$ span the copy of the $\U(n-2)$-type $\mu-e_i$, which is the
  $\gamma^{(i)}$-isotypic component of $V_{m'}$, and in it the weight space of weight $\mu-e_i$ is the line of
  $\Qhw(m';\mu-e_i)$. As $x$ has weight $\mu-e_i$, $P'_{\gamma^{(i)}}x=\beta\,\Qhw(m';\mu-e_i)$. Finally $e\neq0$, as
  $\Phi$ is injective, and $\langle\Phi(x),e\rangle=\langle P'_{\gamma^{(i)}}\Phi(x),e\rangle=\beta\|e\|^2$, as
  $P'_{\gamma^{(i)}}$ is an orthogonal projection with $P'_{\gamma^{(i)}}e=e$.
\end{proof}

Only the $K'$-invariance of the transported Hilbert form was used. We now compute the coefficients. For
$r\in I$ the lowering coefficient of $E_{n,n-1}$ from $m$ to $m-e_r$ on $\xi(m,\mu)$ is, by \eqref{ghl:eq:GTtop},
$$
  d_r^-(m;y)=\frac{\gamma_r(x_r)\prod_{j=1}^{n-2}(x_r-y_j)}{\prod_{s\neq r,\,s\leq n-1}(x_r-x_s)},\qquad
  \gamma_r(x_r)=x_r-y_r\ (r\leq q),\quad\gamma_{n-1}=1.
$$

\begin{lemma}\label{ghl:lem:matrixelement}
  For $m=m'+e_r\in S$, $r\in I$, and $z_r:=x'_r+1$, the coefficient of $\Qhw(m';\mu-e_i)$ in
  $E_{n,i}\xi(m,\mu)$ is
  \begin{equation}\label{ghl:eq:matrixelement}
    \alpha_i\,\frac{\gamma_r(z_r)\prod_{j\in\{1,\dots,q\}\smallsetminus\{i\}}(z_r-y_j)}
    {\prod_{s\in I\smallsetminus\{r\}}(z_r-x'_s)},\qquad
    \alpha_i=\prod_{t=i+1}^{n-2}(y_i-y_t)^{-1}\neq0.
  \end{equation}
\end{lemma}

\begin{proof}
  The lowering coefficients of $E_{k+1,k}$ $(k\leq n-2)$ depend only on rows $k$ and $k-1$. On
  $\xi(m,\mu)$, $l_{k,i}=l_{k-1,i}=y_i$ for $k>i$, so $E_{k+1,k}$ has zero coefficient at position $i$ of
  row $k$ until row $k-1$ has been lowered at $i$. Expanding
  $$
    E_{n-1,i}=[E_{n-1,n-2},[\dots,[E_{i+2,i+1},E_{i+1,i}]\dots]]
  $$
  (for $i=n-2$ just $E_{n-1,n-2}$), only the ordering
  $E_{n-1,n-2}\cdots E_{i+1,i}$, with sign $+$, reaches the highest pattern of weight $\mu-e_i$. Its
  coefficient is $1$ for the step in row $i$ and $1/(y_i-y_k)$ for the step in row $k$, $i<k\leq n-2$,
  independently of row $n-1$. Indeed, row $k-1$ having been lowered at $i$, the coefficient is
  $\tilde A^-_{k,i}=(y_i-(y_i-1))\prod_{j\neq i,\,j<k}(y_i-y_j)/\prod_{s\neq i,\,s\leq k}(y_i-y_s)=1/(y_i-y_k)$. The intermediate patterns are valid, $\mu-e_i$ being dominant. In
  $$
    E_{n,i}=E_{n,n-1}E_{n-1,i}-E_{n-1,i}E_{n,n-1}
  $$
  the two terms carry the factor $(x_r-y_i+1)$, resp.\
  $(x_r-y_i)$, in $d_r^-$, and their difference is $\alpha_i\,d_r^-(m;y)/(x_r-y_i)$. (If the intermediate
  pattern of the first term is not valid, i.\,e.\ $r=i+1$ and $m_{i+1}=\mu_i$, then $x_r=y_i-1$ and that
  term vanishes by itself. The intermediate pattern of the second term, with top row $m'$ and rows
  $1,\dots,n-2$ of weight $\mu$, is always valid, since $m'\in\Rq$ gives $m'_r\geq\tau_r=\mu_r$ for
  $r\leq q$, and lowering $m_{n-1}$ violates no interlacing with $\mu$.) Finally $x_r=z_r$, $x_s=x'_s$ $(s\neq r)$, and for each frozen
  $j\in(q,n-2]$, $y_j=x'_j$, so the factors $(z_r-y_j)/(z_r-x'_j)$ cancel.
\end{proof}

\subsection{Uniqueness}\label{ghl:sec:uniqueness}

Using Lemma \ref{ghl:lem:matrixelement}, we rewrite the reduced equations at each $m'$ as a Cauchy system
and show by induction on $d_0$ that their solutions on $S$ form a line (Proposition \ref{ghl:prop:unique}).

For $m'$ as above let $J(m'):=\{r\in I:m'+e_r\in S\}$ and $E(m'):=\{i\leq q:m'$ interlaces $\mu-e_i\}$. We
record four facts.
\begin{itemize}
  \item[(a)] We have $1,n-1\in J(m')$, since $d_1$ is unbounded and $m'\in\Rq$ forces $m'_r\geq\tau_r$
    $(r\leq q)$, hence $m'_{n-1}\leq\tau_{n-1}-1$.
  \item[(b)] We have $q\in E(m')$. For this one needs $m'_{q+1}\leq\mu_q-1$. For $q\leq n-3$,
    $$
      m'_{q+1}=\tau_{q+1}=L_{q+2}+q\leq L_q+q-2=\mu_q-2.
    $$
    For $q=n-2$,
    $$
      m'_{n-1}\leq\tau_{n-1}-1=L_n+n-3\leq L_{n-2}+n-3=\mu_q-1.
    $$
  \item[(c)] For $1\leq i\leq q-1$: $i\in E(m')\iff m'_{i+1}\leq\mu_i-1\iff i+1\in J(m')$, by
    \eqref{ghl:eq:S} and $\tau_i=\mu_i$. Hence $|E(m')|=|J(m')|-1$, and $z_{i+1}=y_i$ for $i\leq q-1$,
    $i\notin E(m')$.
  \item[(d)] For $r\in J(m')$, $y_r<z_r<y_{r-1}$ $(y_0:=+\infty)$, and $z_{n-1}\leq L_n<y_q$, so the $z_r$,
    $r\in J(m')$, are distinct and disjoint from $\{y_1,\dots,y_q\}$.
\end{itemize}
With
$$
  w_r:=c(m'+e_r)\,\frac{\gamma_r(z_r)\prod_{j=1}^q(z_r-y_j)}{\prod_{s\in I\smallsetminus\{r\}}(z_r-x'_s)},
$$
Lemma \ref{ghl:lem:matrixelement} turns the reduced equations at $m'$ into
$$
  \sum_{r\in J(m')}\frac{w_r}{z_r-y_i}=0\qquad(i\in E(m')).
$$
This is a Cauchy system of size $|E|\times(|E|+1)$ with distinct and disjoint nodes. Every maximal minor is
non-zero, so the kernel is a line and its non-zero vectors have no zero component. The residues of
$$
  R(t)=\frac{\prod_{j\in E(m')}(t-y_j)}{\prod_{s\in J(m')}(t-z_s)}
  =\frac{\prod_{j=1}^q(t-y_j)}{\prod_{s\in I}(t-z_s)}\qquad(z_s:=x'_s+1)
$$
span it (the two products agree since $t-y_k=t-z_{k+1}$ for $k\notin E(m')$, by (c)): $R(t)/(t-y_i)$ is
regular at $y_i$ and $O(t^{-2})$, so its residues at the $z_r$ sum to $0$. Comparing $w_r$ with
$\Res_{t=z_r}R(t)$ and using
$$
  \prod_{s\in I\smallsetminus\{r\}}\frac{x'_r-x'_s+1}{x'_r-x'_s}=\frac{\Delta_I(x'+e_r)}{\Delta_I(x')},
$$
one finds that the reduced equations at $m'$ are equivalent to
\begin{equation}\label{ghl:eq:cauchysol}
  c(m'+e_r)\,\gamma_r(x'_r+1)=C(m')\,\Delta_I(x'+e_r)\qquad(r\in J(m'))
\end{equation}
for some $C(m')\in\CC$.

\begin{proposition}\label{ghl:prop:unique}
  Every solution $c\colon S\to\CC$ of the reduced equations is a scalar multiple of
  \begin{equation}\label{ghl:eq:cdelta}
    c_\delta:=\frac{\Delta_I(x)}{\prod_{r=1}^qd_r!}.
  \end{equation}
\end{proposition}

\begin{proof}
  We proceed by induction on $d_0$. For $d_0=0$ we have $m=\tau$.
  For $m\in S$ with $d_0>0$ take $r\leq q$ with $d_r>0$ and
  $m'=m-e_r$. Then $m'\in\Rq$ (as $d_r>0$) with $\sum m'=\sum\tau-1$ interlaces $\mu-e_q$ (by (b)), $m'+e_{n-1}\in S$ has $d_0-1$, and by
  \eqref{ghl:eq:cauchysol} the equations at $m'$ (the kernel has no zero component) fix
  $$
    \frac{c(m)}{c(m'+e_{n-1})}=\frac{\Delta_I(x)}{d_r\Delta_I(x-e_r+e_{n-1})},
  $$
  since $\gamma_r(x_r)=d_r$ and
  $\gamma_{n-1}=1$. This is the ratio of the values of \eqref{ghl:eq:cdelta}.
\end{proof}

\subsection{The norm on $S$}\label{ghl:sec:norm}

We determine the squared norms $\|\Phi\xi(m,\mu)\|^2$ on $S$ up to a single positive constant (Lemma
\ref{ghl:lem:normS}), using the adjointness \eqref{ghl:eq:adjoint} and the top-row coefficients
\eqref{ghl:eq:GTtop}.

For $x\in\Rq$ interlacing $\mu$ (with any $E_{nn}$-eigenvalue) put $N_\pi(x):=\|\Phi\xi(x)\|^2$,
$\xi(x)=\xi(m,\mu)$. Let
$$
  d^+_r(x;L)=-\frac{\prod_{j\neq g_r}(x_r-L_j)}{\prod_{s\neq r}(x_r-x_s)}
$$
be the raising
coefficient of $E_{n-1,n}$ \eqref{ghl:eq:GTtop}. Regard $d_0:=L_n-x_{n-1}$ as an independent coordinate,
and define on the box $d_1,\dots,d_q,d_0$ as in \eqref{ghl:eq:S} but with $d_0\geq0$ free
\begin{equation}\label{ghl:eq:NS}
  \mathcal N(x):=\frac{1}{\Delta_I(x)}\prod_{r=1}^q\frac{(d_r!)^2}{(a_r+1)_{d_r}(a_r+B+1)_{d_r}}
  \cdot d_0!\,(B+1)_{d_0},
\end{equation}
where $(t)_k$ is the rising factorial.

\begin{lemma}\label{ghl:lem:normS}
  There is $\kappa>0$ with $N_\pi=\kappa\,\mathcal N$ on $S$.
\end{lemma}

\begin{proof}
  By \eqref{ghl:eq:adjoint}, $E_{n-1,n}^*=-E_{n,n-1}$ on $\pi_K$, and $\Phi$ intertwines the actions. Basis
  vectors with different top rows lie in different $K$-types and are orthogonal, and the operators
  $E_{n-1,n}$, $E_{n,n-1}$ change the top row only. Hence, whenever $x$ and $x+e_r$ $(r\in I)$ both belong to the model
  and interlace $\mu$,
  \begin{equation}\label{ghl:eq:steprel}
    d^+_r(x;L)\,N_\pi(x+e_r)=-d^-_r(x+e_r;y)\,N_\pi(x).
  \end{equation}
  The function $\mathcal N$ satisfies the same relations: for $r\leq q$ the frozen factors of $d_r^\pm$
  cancel, and $x_r+1-y_j=x_r-L_j$ $(j\leq q)$ cancels the factors $j\in\{1,\dots,q\}\smallsetminus\{r\}$,
  so that \eqref{ghl:eq:steprel} reads
  \begin{align*}
    \frac{\mathcal N(x+e_r)}{\mathcal N(x)}&=\frac{(x_r-L_r)^2}{(x_r-L_{q+1})(x_r-L_n)}\cdot
    \frac{\Delta_I(x)}{\Delta_I(x+e_r)}\qquad(r\leq q),\\
    \frac{\mathcal N(x-e_{n-1})}{\mathcal N(x)}&=(L_{q+1}-x_{n-1}+1)(L_n-x_{n-1}+1)\cdot
    \frac{\Delta_I(x)}{\Delta_I(x-e_{n-1})},
  \end{align*}
  and with $x_r-L_r=d_r+1$, $x_r-L_{q+1}=u_r+1$, $x_r-L_n=u_r+B+1$, $L_{q+1}-x_{n-1}=B+d_0$ these are the
  ratios of \eqref{ghl:eq:NS}.

  It remains to connect every point of $S$ to $x(\tau)$ by steps along which both coefficients in
  \eqref{ghl:eq:steprel} are non-zero. Since $S$ is a box in $(d_1,\dots,d_q)$ containing $0$, every point
  is reached from $x(\tau)$ by raising one $d_r$ by one at a time through points of $S$. Each unit step is
  $x\to x+e_r\to x+e_r-e_{n-1}$ $(r\leq q)$, taken only with endpoint $x+e_r-e_{n-1}\in S$. Then
  $m_r+1\leq\tau_{r-1}=\mu_{r-1}$ for $r\geq2$, by \eqref{ghl:eq:S} at the endpoint, i.\,e.\ $x_r+1\leq L_{r-1}$,
  so the intermediate top row $x+e_r$ lies in $\Rq$ and interlaces $\mu$. Along these steps
  $x_r\in[L_r+1,L_{r-1}-1]$, so $x_r\neq L_j$ for all $j$, $x_r-L_r\neq0$, and $x_r+1-y_j\neq0$ for all $j$.
  For the $e_{n-1}$-step, $x_{n-1}\leq L_n<y_j$ and $x_{n-1}-1<L_j$ for all $j$. Hence both coefficients are
  non-zero, and $N_\pi/\mathcal N$ is constant on $S$, equal to $N_\pi(x(\tau))/\mathcal N(x(\tau))>0$.
  (For a general point of $S$ the bound $x_r\leq L_{r-1}-1$ fails: $x_r=L_{r-1}$ occurs for $q\geq2$, and
  then $x+e_r\notin\Rq$. This is why the steps are taken towards the endpoint.)
\end{proof}

By Theorem \ref{ghl:thm:identification} the same holds for the form \eqref{ghl:eq:NP} of the model. On $S$
it is a positive multiple of $\mathcal N$.

\subsection{The term formula}\label{ghl:sec:term}

\begin{theorem}\label{ghl:thm:term}
  Let $1\leq q\leq n-2$. For $\delta=(d_1,\dots,d_q)\in S$ put
  $t_\delta:=|c_\delta|^2\mathcal N(\xi_\delta)/(|c_\tau|^2\mathcal N(\xi_\tau))$. Then
  \begin{equation}\label{ghl:eq:tdelta}
    t_\delta=\frac{\Delta_q(u)}{\Delta_q(a)}\prod_{r=1}^q\frac{u_r+B+1+d_0}{a_r+B+1}\cdot
    \frac{d_0!\,(B+1)_{d_0}}{\prod_{r=1}^q(a_r+1)_{d_r}(a_r+B+1)_{d_r}}\;\in\;\QQ_{>0}.
  \end{equation}
  It depends only on $a$ and $B$, not on $n$ or $A_{q+2},\dots,A_{n-1}$, and it is unchanged under a
  rescaling of the basis vectors $\xi_\delta$.
\end{theorem}

\begin{proof}
  By \eqref{ghl:eq:cdelta} and \eqref{ghl:eq:NS} the factors $(d_r!)^2$ cancel. Moreover
  $$
    \frac{\Delta_I(x)^2}{\Delta_I(x)}=\Delta_I(x)=\Delta_q(u)\prod_{r\leq q}(x_r-x_{n-1}),\qquad
    x_r-x_{n-1}=u_r+B+1+d_0.
  $$ Normalize at $\delta=0$. Under $\xi_\delta\mapsto\lambda_\delta\xi_\delta$,
  $c_\delta\mapsto c_\delta/\lambda_\delta$ and $\mathcal N\mapsto|\lambda_\delta|^2\mathcal N$.
\end{proof}
 \section{The summation}\label{ghl:sec:summation}

In this section $q\geq1$, $a_1>\dots>a_q\geq2$ and $B\geq1$ are integers, and $S$, $u_r$, $d_0$ are as in
\eqref{ghl:eq:S}. All terms below are positive, so the order of summation is free.

We evaluate $\sum_{\delta\in S}t_\delta$ (Theorem \ref{ghl:thm:term}) as the product \eqref{ghl:eq:P}
(Corollary \ref{ghl:cor:summation}). The evaluation is an identity of series in the integers $a_r$ and $B$, and
it will be used in the proof of Lemma \ref{ghl:lem:X}.

\subsection{Normalization}\label{ghl:sec:normalization}

Multiplying \eqref{ghl:eq:tdelta} by
$$
  \frac{B!\,\Delta_q(a)\prod_{r=1}^q(a_r+B+1)}{\prod_{r=1}^qa_r!\,(a_r+B)!}
$$
gives
\begin{equation}\label{ghl:eq:Sq}
  S^{(q)}(a;B):=\sum_{\delta\in S}\frac{d_0!\,(B+d_0)!\,\Delta_q(u)\prod_{r=1}^q(u_r+B+1+d_0)}
  {\prod_{r=1}^qu_r!\,(u_r+B)!},
\end{equation}
and since
$$
  \prod_{r=1}^q\frac{a_r(a_r+B)}{(a_r-1)(a_r+B+1)}
$$
is \eqref{ghl:eq:P} ($a_r+B=A_r-A_n$), the identity
$\sum_\delta t_\delta=s(q)$ of \eqref{ghl:eq:P} is equivalent to
\begin{equation}\label{ghl:eq:Sqclosed}
  S^{(q)}(a;B)=\frac{B!\,\Delta_q(a)}{\prod_{r=1}^q(a_r-1)\,(a_r-1)!\,(a_r+B-1)!}.
\end{equation}

\subsection{An auxiliary family}\label{ghl:sec:family}

To prove \eqref{ghl:eq:Sqclosed} by induction on the number of summation variables, we embed $S^{(q)}(a;B)$
in a family of $p$-fold sums with auxiliary variables and a shift of $d_0$. Separating the last summation
variable yields an inner sum of the same shape, which gives the induction step in section \ref{ghl:sec:inductionsec}.

For integers $p\geq1$, $m\geq0$, $B\geq1$, $D\geq0$, $a_1>\dots>a_p$ with $a_p-D-1-m\geq1$, and
indeterminates $z_1,\dots,z_m$, put (sum over $d_1\geq0$, $0\leq d_r\leq a_{r-1}-a_r-1$ $(2\leq r\leq p)$,
$u_r=a_r+d_r$, $d_0=\sum_{r\leq p}d_r$)
\begin{multline*}
  \cI_{p,m}(a;B,D;z):=\sum\frac{(d_0+D)!\,(B+d_0+D)!\,\Delta_p(u)\prod_{r=1}^p(u_r+B+1+d_0+D)}
  {\prod_{r=1}^pu_r!\,(u_r+B)!}\\
  \times\prod_{i=1}^m\Big[(z_i+B+1+d_0+D)\prod_{r=1}^p(u_r-z_i)\Big].
\end{multline*}
Then $\cI_{q,0}(a;B,0)=S^{(q)}(a;B)$, the argument $z$ being omitted for $m=0$. Splitting off the sum over $d_p$
($\Delta_p(u)=\Delta_{p-1}(u')\prod_{r<p}(u_r-u_p)$ with $u'=(u_1,\dots,u_{p-1})$, and the factor of $r=p$
becomes that of a new variable $z_{m+1}=u_p$, with $D'=D+d_p$) gives
\begin{equation}\label{ghl:eq:split}
  \cI_{p,m}(a;B,D;z)=\sum_{d_p=0}^{a_{p-1}-a_p-1}\frac{\prod_{i=1}^m(u_p-z_i)}{u_p!\,(u_p+B)!}\,
  \cI_{p-1,m+1}(a_1,\dots,a_{p-1};B,D+d_p;z,u_p)\qquad(p\geq2).
\end{equation}
For $T\subseteq\{1,\dots,m\}$, $r=|T|$, $\bar T$ its complement and $\one_T$ its indicator function, put
\begin{equation}\label{ghl:eq:CS}
  C_T(B,D;z):=(-1)^r\frac{(D+1)_r(B+D+1)_r\prod_{i\in\bar T}z_i(z_i+B)}
  {\prod_{i\in\bar T}(z_i-D-1-r)\prod_{j\in T}(z_j-D-r)}\prod_{i\in\bar T,\,j\in T}\frac{z_i-z_j-1}{z_i-z_j},
\end{equation}
\begin{equation}\label{ghl:eq:Np}
  \cM_p(a;B):=\frac{\Delta_p(a)}{\prod_{t=1}^p(a_t-1)!\,(a_t+B-1)!},
\end{equation}
\begin{equation}\label{ghl:eq:Phi}
  \Omega_{p,m}(a;B,D;z):=\sum_{T\subseteq\{1,\dots,m\}}C_T(B,D;z)\prod_{t=1}^p
  \frac{\prod_{i=1}^m(a_t-z_i-\one_T(i))}{a_t-D-1-|T|}.
\end{equation}

The coefficients \eqref{ghl:eq:CS} are chosen so that the sum over the last summation variable telescopes.
For $m=1$, $x=z_1$ and $y=D+1$ they are $C_\emptyset=x(x+B)/(x-y)$ and $C_{\{1\}}=-y(y+B)/(x-y)$, and Lemma
\ref{ghl:lem:sumC} is the identity $x(x+B)-y(y+B)=(x-y)(x+y+B)$. For $x=u_1$ and $D+d_1$ in place of $D$ this
writes the factor $u_1+B+1+d_1+D$ of the summand of $\cI_{1,0}$ as a sum of two terms with the denominator
$a_1-D-1$, which does not depend on $d_1$. After multiplication by $(d_1+D)!\,(B+d_1+D)!/(u_1!\,(u_1+B)!)$
the two terms are, with opposite signs, the values at $d_1$ and $d_1+1$ of one function. Hence the sum over $d_1$
telescopes, with value $\cM_1(a_1;B)\,D!\,(B+D)!\,(a_1-D-1)^{-1}$.

For general $m$ and $p\geq2$ the terms come in pairs $T$, $T\cup\{m+1\}$ with $T\subseteq\{1,\dots,m\}$, where $z_{m+1}=u$
is the variable adjoined in \eqref{ghl:eq:split}. Let $a'=(a_1,\dots,a_{p-1})$. By \eqref{ghl:eq:CS}, the
terms of $\prod_{i\leq m}(u-z_i)\,\Omega_{p-1,m+1}(a';B,D;z,u)$ for $T$ and for $T\cup\{m+1\}$ are $u(u+B)$
times a function of $u$ and $D$, and $-(D+1)(B+D+1)$ times the same function at $(u+1,D+1)$. This function is
$\prod_{t<p}(a_t-u)$ times the $T$-term of $\Omega_{p,m}(a',u;B,D;z)$. In \eqref{ghl:eq:split}, $u=u_p$ and
$D+d_p$ increase together with $d_p$, and the factors $u(u+B)$ and $(D+1)(B+D+1)$ are absorbed by
$1/(u_p!\,(u_p+B)!)$ and by the factor $(D+d_p)!\,(B+D+d_p)!$ of \eqref{ghl:eq:main} for $p-1$. Hence the sum
over $d_p$ telescopes for each $T$, the upper boundary term vanishes by the factor $a_{p-1}-u_p$, and the term at $d_p=0$ is $\cM_p(a;B)\,D!\,(B+D)!$ times the $T$-term of
$\Omega_{p,m}(a;B,D;z)$. For $p=1$ the same pairing turns Lemma \ref{ghl:lem:sumC} for the variables
$z_1,\dots,z_m,u$, multiplied by $\prod_{i\leq m}(u-z_i)$, into the recurrence \eqref{ghl:eq:recurrence} for
$a_1=u$.

\subsection{Three lemmas}\label{ghl:sec:threelemmas}

\begin{lemma}\label{ghl:lem:sumC}
  $$
    \sum_{T\subseteq\{1,\dots,m\}}C_T(B,D;z)=\prod_{i=1}^m(z_i+B+D+1).
  $$
\end{lemma}

\begin{proof}
  We proceed by induction on $m$. For $m=0$, $C_\emptyset=1$. Consider the sum as a rational function of $z_1$.
  
  (1) Poles at $z_1=z_k$ $(k\geq2)$ occur only for $T=T_0\cup\{k\}$ and $T'=T_0\cup\{1\}$ $(1,k\notin T_0)$.
  These have the same size and the same factors except $(z_1-z_k-1)/(z_1-z_k)$ against
  $(z_k-z_1-1)/(z_k-z_1)$, with residues $-1$ and $+1$, which cancel.

  (2) Poles at $z_1=D+1+r$
  $(0\leq r\leq m-1)$ occur only for $T\subseteq\{2,\dots,m\}$ with $|T|=r$ and for $T\cup\{1\}$. At $z_1=D+1+r$
  the $j$-product of $T$ becomes $\prod_{j\in T}(z_j-D-1-r)^{-1}$, and for $T\cup\{1\}$ (sign $(-1)^{r+1}$)
  the $i$-product becomes $\prod_{i\in\bar T\smallsetminus\{1\}}(z_i-D-1-r)^{-1}$. For $T$ this comes from the cross
  factors $(z_1-z_j-1)/(z_1-z_j)$, $j\in T$, and the factor $z_1(z_1+B)$ becomes $(D+1+r)(B+D+1+r)$. For
  $T\cup\{1\}$ it comes from the cross factors $(z_i-z_1-1)/(z_i-z_1)$, $i\in\bar T\smallsetminus\{1\}$. Both carry
  $(D+1)_{r+1}(B+D+1)_{r+1}$, and the residues cancel.

  (3) As $z_1\to\infty$, terms with $1\in T$ are
  $O(1/z_1)$, and a term with $1\notin T$, $|T|=r$, is $C_T(B,D;z_2,\dots,z_m)$ times
  $$
    \frac{z_1(z_1+B)}{z_1-D-1-r}\cdot\prod_{j\in T}\frac{z_1-z_j-1}{z_1-z_j}
    =z_1+B+D+1+O(1/z_1).
  $$
  Hence the sum is a polynomial in $z_1$, equal to $(z_1+B+D+1)$ times the sum for $m-1$.
\end{proof}

\begin{lemma}\label{ghl:lem:poly}
  $\Omega_{p,m}(a;B,D;z)$ is a polynomial in $z$ over $\QQ(a,B,D)$, and its denominators are products of
  $a_t-D-1-r$ $(0\leq r\leq m)$, which are invariant under $(a,D)\mapsto(a+k,D+k)$.
\end{lemma}

\begin{proof}
  We argue as in Lemma \ref{ghl:lem:sumC}, (1)--(2). At $z_1=z_k$, the extra factor $\prod_t(\cdots)$ of
  \eqref{ghl:eq:Phi} is the same for $T_0\cup\{k\}$ and $T_0\cup\{1\}$. At $z_1=D+1+r$,
  $$
    \frac{a_t-z_1}{a_t-D-1-r}=1=\frac{a_t-z_1-1}{a_t-D-2-r},
  $$
  so it is the same for $T$ and $T\cup\{1\}$. Since $\Omega$ is symmetric under simultaneous permutation of the
  $z_i$ and the poles are simple and along linear forms, there are none.
\end{proof}

\begin{lemma}\label{ghl:lem:base}
  As rational functions in $a_1,B,D,z$,
  \begin{multline}\label{ghl:eq:recurrence}
    a_1(a_1+B)\Omega_{1,m}(a_1;B,D;z)-(D+1)(B+D+1)\Omega_{1,m}(a_1+1;B,D+1;z)\\
    =(a_1+B+D+1)\prod_{i=1}^m(a_1-z_i)(z_i+B+D+1).
  \end{multline}
  If $a_1-D-1-m\geq1$, the series $\cI_{1,m}$ converges coefficientwise in $z$, and
  $$
    \cI_{1,m}=\cM_1(a_1;B)\,D!\,(B+D)!\,\Omega_{1,m}(a_1;B,D;z)
  $$
  in $\QQ[z]$.
\end{lemma}

\begin{proof}
  (1) We apply Lemma \ref{ghl:lem:sumC} to the $m+1$ variables $z_1,\dots,z_m,z_{m+1}=a_1$ and multiply by
  $\prod_{i\leq m}(a_1-z_i)$. The right side becomes the right side of \eqref{ghl:eq:recurrence}. For
  $T\subseteq\{1,\dots,m\}$, $r=|T|$, the terms for $T$ and $T\cup\{m+1\}$ become, by \eqref{ghl:eq:CS},
  $a_1(a_1+B)$ times the $T$-term of $\Omega_{1,m}(a_1;B,D;z)$ and $-(D+1)(B+D+1)$ times the $T$-term of
  $\Omega_{1,m}(a_1+1;B,D+1;z)$. For the second term we use $(x)_{r+1}=x\,(x+1)_r$ and
  $(a_1-z_i)(z_i-a_1-1)/(z_i-a_1)=a_1+1-z_i$ $(i\notin T)$. This proves \eqref{ghl:eq:recurrence}.

  (2) To telescope the series, put
  $$
    G(k):=\frac{(k+D)!\,(B+k+D)!}{(a_1+k-1)!\,(a_1+B+k-1)!}\Omega_{1,m}(a_1+k;B,D+k;z).
  $$
  By
  \eqref{ghl:eq:recurrence} at $(a_1+k,D+k)$, $G(k)-G(k+1)$ is the $k$-th summand of $\cI_{1,m}$. For distinct
  non-integral $z$ all denominators are non-zero (by Lemma \ref{ghl:lem:poly} they do not involve $k$), and
  a $T$-term of $G(k)$ is $O(k^{2(D+1+r-a_1)})\to0$, since $a_1-D-1-r\geq a_1-D-1-m\geq1$. Hence
  $\cI_{1,m}=G(0)$ on a Zariski dense set. Each $z$-monomial of the summand is $O(k^{-3})$, so both sides are
  polynomials in $z$ (Lemma \ref{ghl:lem:poly}) and agree. Since $G(0)=\cM_1D!(B+D)!\Omega_{1,m}$, the claim
  follows.
\end{proof}

\subsection{The induction}\label{ghl:sec:inductionsec}

\begin{theorem}\label{ghl:thm:induction}
  For all $p\geq1$, $m\geq0$, $B\geq1$, $D\geq0$ and integers $a_1>\dots>a_p$ with $a_p-D-1-m\geq1$,
  \begin{equation}\label{ghl:eq:main}
    \cI_{p,m}(a;B,D;z)=\cM_p(a;B)\,D!\,(B+D)!\,\Omega_{p,m}(a;B,D;z)\qquad\text{in }\QQ[z_1,\dots,z_m].
  \end{equation}
\end{theorem}

\begin{proof}
  We proceed by induction on $p$, for all $m$ simultaneously. The case $p=1$ is Lemma \ref{ghl:lem:base}.
  Let $p\geq2$, $a'=(a_1,\dots,a_{p-1})$, $D'=D+d_p$, $u=u_p=a_p+d_p$, $0\leq d_p\leq a_{p-1}-a_p-1$.

  The hypotheses propagate, since $a_{p-1}-D'-1-(m+1)\geq a_p-D-1-m\geq1$. Therefore \eqref{ghl:eq:main} holds
  for $(p-1,m+1,a',D')$ as a polynomial identity in $z_1,\dots,z_{m+1}$ and may be specialized at $z_{m+1}=u$.
  On the left this gives the summand of \eqref{ghl:eq:split}. On the right $\Omega_{p-1,m+1}$ may be evaluated
  termwise: with $r=|\tilde T\smallsetminus\{m+1\}|$ for the term $\tilde T\subseteq\{1,\dots,m+1\}$, the
  denominators containing $z_{m+1}$ become $u-D'-1-r=a_p-D-1-r\neq0$ $(r\leq m)$, whether $m+1\in\tilde T$ or
  not, and the others are $a_t-D'-1-\tilde r\geq a_p-D-m-1\geq1$ $(t<p,\ \tilde r\leq m+1)$, or $z_i-u$ with
  $z_i$ formal.

  {\em Pairing.} Split $\tilde T\subseteq\{1,\dots,m+1\}$ into $T$ and $T\cup\{m+1\}$, $T\subseteq\{1,\dots,m\}$,
  $|T|=r$. After multiplication by
  $$
    \frac{\prod_i(u-z_i)}{u!\,(u+B)!}\cdot\cM_{p-1}(a')D'!(B+D')!,
  $$
  the term $T$
  gives $H_T(d_p)$ and the term $T\cup\{m+1\}$ gives $-H_T(d_p+1)$, where, with $u=a_p+d$,
  $$
    H_T(d):=\frac{(-1)^r\cM_{p-1}(a';B)}{a_p-D-1-r}\,\frac{(D+d+r)!\,(B+D+d+r)!}{(u-1)!\,(u+B-1)!}
  $$
  $$
    \times\frac{\prod_{i\in\bar T}z_i(z_i+B)\prod_{i\in\bar T,j\in T}\frac{z_i-z_j-1}{z_i-z_j}}
    {\prod_{i\in\bar T}(z_i-D-d-1-r)\prod_{j\in T}(z_j-D-d-r)}
  $$
  $$
    \times\prod_{i=1}^m(u-z_i-\one_T(i))\prod_{t=1}^{p-1}\frac{(a_t-u)\prod_{i=1}^m(a_t-z_i-\one_T(i))}
    {a_t-D-d-1-r}.
  $$
  Indeed, for $T$, we have
  $$
    \frac{u(u+B)}{u!(u+B)!}=\frac{1}{(u-1)!(u+B-1)!}
  $$
  and
  $$
    \prod_i(u-z_i)\prod_{j\in T}\frac{u-z_j-1}{u-z_j}=\prod_i(u-z_i-\one_T(i)).
  $$
  For $T\cup\{m+1\}$ the sign flips,
  \begin{gather*}
    (D'+1)_{r+1}(B+D'+1)_{r+1}D'!(B+D')!=(D+d_p+1+r)!(B+D+d_p+1+r)!,\\
    u!(u+B)!=((u+1)-1)!((u+1)+B-1)!,\\
    \prod_i(u-z_i)\prod_{i\in\bar T}\frac{z_i-u-1}{z_i-u}
    =\prod_i((u+1)-z_i-\one_T(i)),
  \end{gather*}
  and $u-D'-(r+1)=a_p-D-1-r$ does not depend on $d_p$.

  At $d=a_{p-1}-a_p$, we have $u=a_{p-1}$ and the factor $a_{p-1}-u$ vanishes, while all denominators stay
  $\geq a_p-D-1-m\geq1$, so $H_T(a_{p-1}-a_p)=0$. At $d=0$,
  $$
    \cM_{p-1}(a')\prod_{t<p}\frac{a_t-a_p}{(a_p-1)!(a_p+B-1)!}=\cM_p(a),
  $$
  and $H_T(0)$ is $\cM_p(a)D!(B+D)!$ times the $T$-term of $\Omega_{p,m}(a;B,D;z)$.

  Summing up, we obtain
  $$
  \cI_{p,m}=\sum_T\sum_{d_p}(H_T(d_p)-H_T(d_p+1))=\sum_TH_T(0),
  $$
  which is \eqref{ghl:eq:main}. Only finite sums are interchanged: the $d_p$-sum is finite,
  and the inner series converge coefficientwise.
\end{proof}

\begin{corollary}\label{ghl:cor:summation}
  For $q\geq1$, \eqref{ghl:eq:Sqclosed} holds, i.\,e.\ $\sum_{\delta\in S}t_\delta$ equals \eqref{ghl:eq:P}.
\end{corollary}

\begin{proof}
  Take $p=q$, $m=0$, $D=0$ in \eqref{ghl:eq:main}. This is allowed since $a_q-1\geq1$, and
  $\Omega_{q,0}=\prod_r(a_r-1)^{-1}$.
\end{proof}
 \section{The Schmid vector generates $\sigma$}\label{ghl:sec:converse}

In this section we prove \hypref{ghl:hyp:J}{J} for $q\geq1$. Throughout, $1\leq q\leq n-2$ and $A$ satisfies
\eqref{ghl:eq:gap}. The objects $\Phi$, $\xi(m,\mu)$, $y_i$, $S$, $I$, $J(m')$, $E(m')$ and $c_\delta$ are as in section
\ref{ghl:sec:schmid}, and $\rho=\hatpi|_{H'}$ as in section \ref{ghl:sec:lemmaP}. We write $\xi_\delta:=\xi(\delta,\mu)$
$(\delta\in S)$, $u_0$ for a highest weight vector of $\tau_\sigma$, and $p:=n-2$, so that $H'\cong\U(p,1)$.

Section \ref{ghl:sec:Hprime} describes the $K'$-types and transitions of all discrete series of $H'$.
Sections \ref{ghl:sec:weak} to \ref{ghl:sec:Id} construct the Schmid vector, prove that its closed
$H'$-span is an irreducible discrete series and identify it with $\hatsigma$. Theorem
\ref{ghl:thm:J} gives \hypref{ghl:hyp:J}{J}.

\subsection{The discrete series of $H'$}\label{ghl:sec:Hprime}

In section \ref{ghl:sec:span} the closed span of the Schmid vector is first identified only as some discrete
series of $H'$ (Lemma \ref{ghl:lem:classif}). Proposition \ref{ghl:prop:Hprime} therefore records the
$K'$-types and $\liep'_\CC$-transitions of every discrete series of $H'$, using section \ref{ghl:sec:gt},
conjugation (Lemma \ref{ghl:lem:conj}) and Blattner's formula for $n=3$, as needed in Lemma \ref{ghl:lem:Id}.

The coordinates of the diagonal torus $T'$ of $K'$ are indexed by $1,\dots,p$ and $0$ (for $\U(1)_{e_n}$), as for
$K'$-types. Let $A''=(A''_1>\dots>A''_{p+1})\in(\ZZ+\frac p2)^{p+1}$ and $0\leq q'\leq p$, define $A''(q')$ by
\eqref{ghl:eq:Aq} with $p+1$ in place of $n$, and let $\zeta=\pi_{A''(q')}$ be the discrete series of $H'$ with
Harish-Chandra parameter $A''(q')$. The centre $\U(1)\cdot\one$ of $H'$ acts on it by $z\mapsto z^{\mathrm{tot}'}$,
$\mathrm{tot}':=\sum_jA''_j$. Put $L''_j:=A''_j-\frac p2\in\ZZ$, $L''_0:=+\infty$, $L''_{p+2}:=-\infty$, write
$x_i(m):=m_i-i+1$ for a top row $m=(m_1\geq\dots\geq m_p)$, and let
$$
  J_i:=[L''_i+1,L''_{i-1}]\ (i\leq q'),\qquad J_i:=[L''_{i+2}+1,L''_{i+1}]\ (q'<i\leq p),\qquad
  R':=J_1\times\dots\times J_p.
$$
The $J_i$ are non-empty integer intervals and every element of $J_i$ exceeds every element of $J_{i+1}$, so
$x\in R'$ is strictly decreasing, i.\,e.\ $m$ is dominant. For a $K'$-type with top row $m$ we write $P'_{(m)}$ for
its isotypic projection when the $\U(1)$-coordinate is fixed by the central character.

\begin{proposition}\label{ghl:prop:Hprime}
  \noindent {\rm(i)} The $K'$-types of $\zeta$ are the $(m;\mathrm{tot}'-\sum m)$ with $x(m)\in R'$, each with multiplicity
  one.

  \noindent{\rm(ii)} If $(m;\mathrm{tot}'-\sum m)$ and $(m+e_i;\mathrm{tot}'-\sum m-1)$ are both $K'$-types of $\zeta$, then
  $$
    P'_{(m+e_i)}\,\zeta(\liep'_\CC)\,P'_{(m)}\neq0\quad\text{and}\quad P'_{(m)}\,\zeta(\liep'_\CC)\,P'_{(m+e_i)}\neq0
  $$
  on $\zeta_{K'}$.
\end{proposition}

For $n\geq4$ and $q'\leq p-1$ this is section \ref{ghl:sec:gt} for $H'$. The chamber $q'=p$, which is excluded there,
is reduced to the case $q=0$ by complex conjugation.

\begin{lemma}\label{ghl:lem:conj}
  Let $p\geq2$ and $q'=p$. Put $B:=(-A''_{p+1}>\dots>-A''_1)$, let $V_B$ be the module $V_{R_0}$ of section
  \ref{ghl:sec:gt} for $\U(p,1)$, i.\,e.\ with $p+1$ in place of $n$, $q=0$, and $L^B_j:=B_j-\frac p2$ in place of $L_j$,
  and let $M$ be its complex conjugate, defined below. Then $M\cong\zeta_{K'}$, and {\rm(i)}, {\rm(ii)} of Proposition
  \ref{ghl:prop:Hprime} hold for $\zeta$.
\end{lemma}

\begin{proof}
  (1) {\em The model.} Section \ref{ghl:sec:gt} assumes no gap condition, so it applies to $\U(p,1)$, $p+1\geq3$,
  $q=0$, and the parameter $L^B$: $V_B$ is an irreducible admissible multiplicity free $(\gl_{p+1},K')$-module
  with a positive definite form satisfying \eqref{ghl:eq:adjoint} (Theorems \ref{ghl:thm:form},
  \ref{ghl:thm:relations}, Propositions \ref{ghl:prop:Ktypes}, \ref{ghl:prop:irred}), i.\,e.\ $\langle Yw,w'\rangle=
  -\langle w,Yw'\rangle$ for $Y\in\lieh'$. Since $L^B_j=-L''_{p+2-j}-p$, its $K'$-types are the
  $(m;\mathrm{tot}_B-\sum m)$ with $x_k(m)\in[L^B_{k+2}+1,L^B_{k+1}]$ $(1\leq k\leq p$, $L^B_{p+2}=-\infty)$, where
  $\mathrm{tot}_B=\sum_jB_j=-\mathrm{tot}'$.

  (2) {\em The conjugate module.} Let $M:=\overline{V_B}$ be the conjugate vector space, with $c\cdot\bar w:=
  \overline{\bar cw}$, and put $k\cdot\bar w:=\overline{kw}$ $(k\in K')$ and $Y\cdot\bar w:=\overline{Yw}$
  $(Y\in\lieh')$, extended $\CC$-linearly, so that $Z\cdot\bar w=\overline{\bar Zw}$ for $Z\in\gl_{p+1}$, where
  $Z\mapsto\bar Z$ is the conjugation with respect to $\lieh'$. Since $Y^*=-Y$ for $Y\in\lieh'$ by
  \eqref{ghl:eq:adjoint}, $\bar Z=-Z^*$ for the adjoint $Z^*$ on $V_B$. As $Z\mapsto\bar Z$ is an automorphism of the real Lie algebra $\gl_{p+1}$
  commuting with $\mathrm{Ad}(K')$ and equal to the identity on $\liek'$, $M$ is a $(\gl_{p+1},K')$-module. The form
  $\langle\bar w,\bar w'\rangle_M:=\langle w',w\rangle$ is positive definite, and for $Y\in\lieh'$
  $$
    \langle Y\cdot\bar w,\bar w'\rangle_M=\langle w',Yw\rangle=-\langle Yw',w\rangle=-\langle\bar w,Y\cdot\bar w'\rangle_M,
  $$
  so $M$ satisfies \eqref{ghl:eq:adjoint}. The complex subspaces of $M$ are the conjugates of those of $V_B$, with the
  same stability properties, so $M$ is irreducible. For a $K'$-type $\delta$ the isotypic projection
  $$
  P_\delta=\dim\delta\int_{K'}\overline{\chi_\delta(k)}\,k\,dk
  $$
  satisfies $P^M_\delta\bar w=\overline{P^{V_B}_{\delta^\vee}w}$,
  since $\chi_{\delta^\vee}=\overline{\chi_\delta}$. Hence $M$ is admissible and multiplicity free, and its $K'$-types are
  the duals $(m;m_0)^\vee=(-m_p,\dots,-m_1;-m_0)$ of those of $V_B$. As an irreducible admissible module with a
  positive definite invariant form, $M$ is the Harish-Chandra module of an irreducible unitary, hence quasi-simple,
  representation.

  (3) {\em The dual box.} For $m^\vee:=(-m_p,\dots,-m_1)$ and $k:=p+1-i$ we have $x_i(m^\vee)=-x_k(m)-p+1$, so
  $$
    x_k(m)\in[L^B_{k+2}+1,L^B_{k+1}]\iff x_i(m^\vee)\in[-L^B_{k+1}-p+1,-L^B_{k+2}-p]=[L''_i+1,L''_{i-1}]=J_i,
  $$
  with $L^B_{p+2}=-\infty$ corresponding to $L''_0=+\infty$. The $\U(1)$-coordinate of the dual is
  $$
    -\left(\mathrm{tot}_B-\sum m\right)=\mathrm{tot}'-\sum m^\vee.
  $$
  Therefore the $K'$-types of $M$ are the $(m;\mathrm{tot}'-\sum m)$ with
  $x(m)\in R'$, each once.

  (4) {\em Schmid's criterion.} We reduce to $G_1=\SU(p,1)$ as in the proof of Theorem \ref{ghl:thm:identification}. The
  centre acts on $M$ by $z\mapsto\overline{z^{\mathrm{tot}_B}}=z^{\mathrm{tot}'}$, as on $\zeta$. Let
  $\lambda:=A''(p)=(A''_1,\dots,A''_p;A''_{p+1})$. As $A''_i>A''_{p+1}$ for $i\leq p$, the positive system
  $\Psi=\{\alpha:\langle\lambda,\alpha\rangle>0\}$ has the standard compact part and the non-compact part
  $\beta'_i=e_i-e_n$ $(i\leq p)$. With $\rho_{c,i}=\frac{p+1}2-i$ and $\rho_n=\frac12\sum_i\beta'_i$, the weight
  $\tau':=\lambda-\rho_c+\rho_n$ has $\tau'_i=L''_i+i$ $(i\leq p)$ and $\tau'_0=L''_{p+1}$.

  \noindent (a) We have $x_i(\tau')=L''_i+1\in J_i$
  and $\mathrm{tot}'-\sum_i\tau'_i=L''_{p+1}=\tau'_0$, so $\tau'$ is a $K'$-type of $M$.

  \noindent (b) Every $K'$-type of $M$ is
  $\tau'+\sum_ic_i\beta'_i$ with $c_i=x_i(m)-L''_i-1\in\ZZ_{\geq0}$, the $\U(1)$-coordinate being $\tau'_0-\sum_ic_i$. If
  $\tau'-\Sigma=\tau'+\sum_ic_i\beta'_i$ with $\Sigma$ a non-empty sum of roots in $\Psi$, then $\Sigma+\sum_ic_i\beta'_i=0$ is a non-negative
  combination of roots in $\Psi$, and pairing with $\lambda$ gives $\Sigma=0$, a contradiction. By \cite[Theorem (1.3),
  pp.~536--537]{schmid75}, $M\cong\zeta_{K'}$, which gives (i).

  (5) {\em Transitions.} Let $(m;\cdot)$, $(m+e_i;\cdot)$ be $K'$-types of $\zeta$. For $w\in V_B$ and $Z\in\liep'_\CC$,
  $$
    P^M_{\delta'}(Z\cdot\bar w)=\overline{P^{V_B}_{\delta'^\vee}\bar Zw},
  $$
  and $\bar Z$ runs through $\liep'_\CC$ with $Z$, as $\liep'$ is real. Hence
  $$
    P^M_{(m+e_i)}\liep'_\CC P^M_{(m)}\neq0
  $$
  if and only if
  $$
    P^{V_B}_{(m^\vee-e_k)}\liep'_\CC P^{V_B}_{(m^\vee)}\neq0,
  $$
  where $k=p+1-i$, and likewise in the other direction. By (3), $m^\vee-e_k$ and $m^\vee$ both
  lie in $R_0$ for $V_B$, and the argument for the case $q'\leq p-1$ in the proof of Proposition
  \ref{ghl:prop:Hprime} (ii) below, applied to the model $V_B$ (parameter $L^B$, $q'=0$; it is the proof of
  Proposition \ref{ghl:prop:irred} and uses only the model), gives the non-zero transitions in $V_B$.
\end{proof}

\begin{proof}[Proof of Proposition \ref{ghl:prop:Hprime}]
  {\em The case $n\geq4$, $q'\leq p-1$.} Section \ref{ghl:sec:gt} applies to $H'=\U(n-2,1)$, with $n-1\geq3$ in place of $n$,
  $L''$ in place of $L$ and $q'\leq(n-1)-2$ in place of $q$. The indices $1,\dots,n-1$ there are the indices
  $1,\dots,n-2,n$ here. Condition \eqref{ghl:eq:Rq} for $n-1$ is the definition of the $J_i$, so (i) is Corollary
  \ref{ghl:cor:ktypes} for $H'$, and $\zeta_{K'}\cong V'':=V_{R_{q'}}$ by Theorem \ref{ghl:thm:identification}. For
  (ii) we work in $V''$. If $m$, $m+e_i\in R'$ and $i\geq2$, then $m_{i-1}\geq m_i+1$ (proof of Proposition
  \ref{ghl:prop:irred}; for $i=q'+1$ it uses only $L''_{q'}>L''_{q'+2}$), so $m$ with its $i$-th entry removed
  interlaces both $m$ and $m+e_i$. For $i=1$ this holds without condition. Let $P$ be the pattern with top row $m$, row $n-3$ equal to $m$ without its $i$-th
  entry, and the highest pattern of that weight below it, and $Q:=P+\delta_{n-2,i}$. Both lie in $\cP(R')$. By
  Theorem \ref{ghl:thm:form} (2), $\tilde A^+_{n-2,i}(P)<0$ and $\tilde A^-_{n-2,i}(Q)>0$. Hence, in the indices
  of section \ref{ghl:sec:gt} for $n-1$, $E_{n-2,n-1}\xit_P$ has the non-zero $\xit_Q$-component
  $\tilde A^+_{n-2,i}(P)$, and $E_{n-1,n-2}\xit_Q$ the non-zero $\xit_P$-component $\tilde A^-_{n-2,i}(Q)$.
  Here $\xit_Q$ lies in the
  $K'$-type $(m+e_i;\cdot)$, which has multiplicity one. In the indices of $H'$ the two operators are $E_{n-2,n}$ and
  $E_{n,n-2}$, both in $\liep'_\CC$.

  {\em The case $n\geq4$, $q'=p$.} This follows from Lemma \ref{ghl:lem:conj}.

  {\em The case $n=3$.} Here $H'=\U(1,1)$, $K'=\U(1)\times\U(1)_{e_3}$ is a torus, $W_{K'}$ is trivial, $\rho_c=0$, and there is
  one positive non-compact root, $\beta'=e_1-e_3$ if $q'=1$ and $\beta'=e_3-e_1$ if $q'=0$. Blattner's formula
  \cite[Theorem (1.3), p.~130]{hechtschmid}, applied to $G_1=\SU(1,1)$ as in Lemma \ref{ghl:lem:Kprime}, gives the
  multiplicity $Q(\delta-\tau')\in\{0,1\}$ of a $K'$-type $\delta$ with the central character of $\zeta$, where
  $\tau'=\lambda+\rho_n$, so the $K'$-types are the $\tau'+j\beta'$, $j\geq0$, each once. For $q'=1$,
  $\tau'=(L''_1+1;L''_2)$, and these are the $(m_1;\mathrm{tot}'-m_1)$ with $m_1\in[L''_1+1,\infty)=J_1$. For $q'=0$,
  $\tau'=(L''_2;L''_1+1)$, and $m_1\in(-\infty,L''_2]=J_1$. This is (i). For (ii) let $W:=\zeta_{K'}$. By (i) $\zeta$ is admissible, so
  $W$ is an irreducible $(\lieh'_\CC,K')$-module, and
  $W=\bigoplus_{j\geq0}W_j$, $W_j:=W(\tau'+j\beta')$, $\dim W_j=1$. Let $Z_\pm\in\liep'_\CC$ be root vectors of
  $\pm\beta'$ ($E_{13}$, $E_{31}$ in some order). Then $\liep'_\CC=\langle Z_+,Z_-\rangle$, $Z_\pm W_j\subset W_{j\pm1}$
  ($W_{-1}:=0$), $\liek'_\CC$ preserves each $W_j$, and $Z_+^*=-Z_-$ on $W$ by \eqref{ghl:eq:adjoint}. Fix $k\geq0$. If
  $Z_+W_k=0$, then $\langle Z_-w',w\rangle=-\langle w',Z_+w\rangle=0$ for $w\in W_k$, $w'\in W_{k+1}$, and as
  $Z_-W_{k+1}\subset W_k$, $Z_-W_{k+1}=0$. Then $\bigoplus_{j\geq k+1}W_j$ is a non-zero proper submodule, a
  contradiction. Vice versa, if $Z_-W_{k+1}=0$, the same adjointness gives $Z_+W_k=0$. Therefore both transitions between
  $\tau'+k\beta'$ and $\tau'+(k+1)\beta'$ are non-zero. These are the pairs $m$, $m+e_1$ in $J_1$.
\end{proof}

\subsection{The vector $v$ and the weak Schmid conditions}\label{ghl:sec:weak}

We pass from the model $V_{\Rq}$ (section \ref{ghl:sec:gt}) and the formulas of section \ref{ghl:sec:schmid}
to the Hilbert space $\hatpi$. Using Lemma \ref{ghl:lem:normS}, Theorem \ref{ghl:thm:term} and Corollary
\ref{ghl:cor:summation}, Lemma \ref{ghl:lem:X} constructs the non-zero vector $v\in\hatpi$, which is then
characterized by weak Schmid conditions in Lemmas \ref{ghl:lem:WS} and \ref{ghl:lem:sol}.

Let $\rho^\infty$ be the space of $H'$-smooth vectors of $\hatpi$ with its Fr\'echet topology. For $e\in\rho^\infty$
and $Y\in\lieh'$ put
$$
  \rho(Y)e:=\frac d{dt}|_{t=0}\pi(\exp tY)e,
$$
extended $\CC$-linearly to $\lieh'_\CC$. On $\pi_K$ this
is the action of $\lieh'_\CC$. The projections $P'_\gamma$ are integrals over $K'\subset H'$. They preserve
$\rho^\infty$ and are continuous on it.

\begin{lemma}\label{ghl:lem:density}
  $\pi_K$ is dense in $\rho^\infty$.
\end{lemma}

\begin{proof}
  For $\varphi\in C_c^\infty(H')$ the operator
  $$
    \rho(\varphi)=\int_{H'}\varphi(h)\rho(h)\,dh'
  $$
  maps $\hatpi$ continuously to $\rho^\infty$, and it maps $\pi^\infty$ to $\pi^\infty$, the integral converging in the Fr\'echet
  space $\pi^\infty$. As $\pi^\infty$ is dense in $\hatpi$, $\rho(\varphi)\hatpi$ lies in the closure of $\pi^\infty$
  in $\rho^\infty$. The G{\aa}rding space, spanned by the $\rho(\varphi)\hatpi$, is dense in $\rho^\infty$.
  Furthermore, $\pi_K$ is dense in $\pi^\infty$, and on $\pi^\infty$ the topology of $\pi^\infty$ is finer than
  that of $\rho^\infty$, since $U(\lieh'_\CC)\subset U(\lieh_\CC)$.
\end{proof}

Let $\mathcal X\subset\hatpi$ be the closed subspace of vectors of $K'$-type $\tau_\sigma$ and $T'$-weight
$(\mu;\tau_0)$, the range of the product of $P'_{\tau_\sigma}$ and the $T'$-weight projection. As $(\mu;\tau_0)$ is
the highest weight of $\tau_\sigma$, $\mathcal X$ consists of $K'$-highest weight vectors.

\begin{lemma}\label{ghl:lem:X}
  \noindent {\rm(i)} $\mathcal X$ is the Hilbert sum of the lines $\CC\Phi(\xi_\delta)$, $\delta\in S$, so every $x\in\mathcal X$
  is
  $$
    x=\sum_{\delta\in S}c_x(\delta)\Phi(\xi_\delta),
  $$
  norm convergent, with $c_x\colon S\to\CC$.
    
  \noindent{\rm(ii)} The series
  $$
    v:=\sum_{\delta\in S}c_\delta\,\Phi(\xi_\delta)
  $$
  converges in $\hatpi$, $0\neq v\in\mathcal X$, and $\|v\|^2=|c_\tau|^2N_\pi(\xi_\tau)\,s$ with $s$ the
  product~\eqref{ghl:eq:P}.
\end{lemma}

\begin{proof}
  (i) The two projections commute with every $P_\delta$, and $\hatpi$ is the Hilbert sum of its $K$-isotypic
  subspaces. By Corollary \ref{ghl:cor:ktypes} and the multiplicity free branching, $\pi_K(\delta)$, $\delta=(m;m_0)$,
  contains the $K'$-type $\tau_\sigma$ if and only if $m\in\Rq$ interlaces $\mu$ and $m_0=\tau_0$, and then once, with
  $T'$-weight line $\CC\Phi(\xi(m,\mu))$. Here $m_0=\tau_0$ means $\sum m=\sum\tau$, so $m\in S$ by Lemma
  \ref{ghl:lem:support}.

  (ii) The $\Phi(\xi_\delta)$ are orthogonal, and by Lemma \ref{ghl:lem:normS} and the
  definition of $t_\delta$ in Theorem \ref{ghl:thm:term},
  $$
    \sum_\delta|c_\delta|^2N_\pi(\xi_\delta)=\kappa\sum_\delta|c_\delta|^2\mathcal N(\xi_\delta)=
    |c_\tau|^2N_\pi(\xi_\tau)\sum_\delta t_\delta,
  $$
  which is finite and equal to $|c_\tau|^2N_\pi(\xi_\tau)s$ by Corollary
  \ref{ghl:cor:summation}. We also have $c_\tau=\Delta_I(x(\tau))\neq0$.
\end{proof}

Since $v$ is not known a priori to be $H'$-smooth, we formulate the Schmid conditions weakly against $H'$-smooth
vectors, a condition preserved by the commutant of $\rho(H')$ (Proposition \ref{ghl:prop:Virred}). By Lemma
\ref{ghl:lem:density} they reduce on $\mathcal X$ to the reduced equations of section \ref{ghl:sec:reduced},
so Proposition \ref{ghl:prop:unique} fixes the line $\CC v$ in Lemma \ref{ghl:lem:sol}.

We say that $x\in\hatpi$ satisfies the {\em weak Schmid conditions} if
\begin{equation}\label{ghl:eq:WS}
  \langle x,\rho(Y)e\rangle=0\qquad(Y\in\liep'_\CC,\ e\in P'_{\gamma^{(i)}}\rho^\infty,\ 1\leq i\leq q),
\end{equation}
and the {\em $K$-finite weak Schmid conditions} if this holds for $e\in P'_{\gamma^{(i)}}\pi_K$. For $q<i\leq n-2$ no
condition is imposed: $\gamma^{(i)}=(\mu+e_i;\tau_0-1)$ is not a $K'$-type of $\rho$, since a $K$-type $(m;\cdot)$ containing
it would have $m_i\geq\mu_i+1=\tau_i+1$, against \eqref{ghl:eq:Rq}. Let $\mathrm{Sol}$ be the set of $x\in\mathcal X$
satisfying \eqref{ghl:eq:WS}.

\begin{lemma}\label{ghl:lem:WS}
  \noindent {\rm(a)} For $x\in\hatpi$, \eqref{ghl:eq:WS} is equivalent to the $K$-finite weak Schmid conditions.

  \noindent {\rm(b)} For $x\in\mathcal X$, the $K$-finite weak Schmid conditions hold if and only if $c_x$ satisfies every
  reduced equation $(m',i)$ of section \ref{ghl:sec:reduced}.

  \noindent {\rm(c)} $c_\bullet=(c_\delta)_{\delta\in S}$ satisfies every reduced equation.
\end{lemma}

\begin{proof}
  (a) One direction is trivial. Let $e\in P'_\gamma\rho^\infty$, $\gamma=\gamma^{(i)}$. By Lemma
  \ref{ghl:lem:density} there are $e_k\in\pi_K$ with $e_k\to e$ in $\rho^\infty$. Then $P'_\gamma e_k\in P'_\gamma\pi_K$
  and $P'_\gamma e_k\to P'_\gamma e=e$ in $\rho^\infty$, $\rho(Y)$ is continuous on $\rho^\infty$, and
  $\rho^\infty\to\hatpi$ is continuous.

  (b) Fix $i\leq q$. The $K$-types of $\pi$ containing $\gamma^{(i)}=(\mu-e_i;\tau_0+1)$ are the
  $(m';\mathrm{tot}-\sum m')$ with $m'\in\Rq$ interlacing $\mu-e_i$ and $\sum m'=\sum\tau-1$, each containing it once.
  These are the $m'$ of the reduced equations $(m',i)$. (If $\mu-e_i$ is not dominant there is no such $m'$, and
  nothing to prove.) By linearity in $e$, the $K$-finite conditions for $i$ say that $\langle x,\pi(Y)e\rangle=0$ for
  all $Y\in\liep'_\CC$, all such $m'$ and all $e\in W':=P_{m'}P'_{\gamma^{(i)}}\pi_K\cong\gamma^{(i)}$. Fix $m'$. By
  \eqref{ghl:eq:adjoint} every $Y\in\liep'_\CC$ has an adjoint $Y^\dagger\in\liep'_\CC$ on $\pi_K$
  ($E_{an}^\dagger=-E_{na}$, $E_{na}^\dagger=-E_{an}$, $a\leq n-2$), and $Y\mapsto Y^\dagger$ is a bijection of
  $\liep'_\CC$. A $\delta=(m;\tau_0)\in S$ with $P_{m'}\pi(\liep'_\CC)\pi_K(\delta)\neq0$ has $m=m'+e_r$ for some $r$, as
  $\liep'_\CC$ changes the $\U(1)_{e_n}$-weight by $\pm1$ and the top row by a unit vector, so it is one of the
  $m'+e_r$, $r\in J(m')$ (section \ref{ghl:sec:reduced}). Let $x_F$ be the corresponding finite partial sum of $x$. For $e\in W'$,
  $\langle x,\pi(Y)e\rangle=\langle\pi(Y^\dagger)x_F,e\rangle$.

  For $\delta\in S$ let $\iota_\delta\colon\tau_\sigma\to\pi_K(\delta)$ be the $K'$-map with $\iota_\delta(u_0)=
  \Phi(\xi_\delta)$, and define
  $$
    T\colon\liep'_\CC\otimes\tau_\sigma\to W',\qquad
    T(Z\otimes w):=P_{m'}P'_{\gamma^{(i)}}\,\pi(Z)\sum_{r\in J(m')}c_x(m'+e_r)\,\iota_{m'+e_r}(w).
  $$
  The map $T$ is $K'$-equivariant, as $\pi(\mathrm{Ad}(k)Z)\pi(k)=\pi(k)\pi(Z)$ and $\iota_\delta$, $P_{m'}$,
  $P'_{\gamma^{(i)}}$ commute with $K'$. Furthermore,
  $\langle x,\pi(Y)e\rangle=\langle T(Y^\dagger\otimes u_0),e\rangle$ for $e\in W'$. As $W'$ is
  irreducible, $T$ vanishes on every $K'$-isotypic component of $\liep'_\CC\otimes\tau_\sigma$ other than the
  $\gamma^{(i)}$-component. The latter is irreducible: with $\liep'^+:=\langle E_{n,k}:k\leq n-2\rangle$ ($E_{n,k}$ has
  weight $e_n-e_k$) and $\liep'^-:=\langle E_{k,n}:k\leq n-2\rangle$,
  $$
    \liep'^+\otimes\tau_\sigma\cong((\CC^p)^*\otimes F_\mu)\boxtimes\chi_{\tau_0+1}
  $$
  contains $F_{\mu-e_i}\boxtimes\chi_{\tau_0+1}$ once by Pieri's rule,
  and $\liep'^-\otimes\tau_\sigma$ has $\U(1)_{e_n}$-weight $\tau_0-1$, so $T$ is determined by $T(h)$,
  where $h$ is a highest weight vector of the $\gamma^{(i)}$-component.

  {\em The $\gamma^{(i)}$-component of $E_{n,i}\otimes u_0$ is a non-zero multiple of $h$.} Fix the tensor product of
  $K'$-invariant inner products on $\liep'^+\otimes\tau_\sigma$. Vectors of distinct weights are orthogonal, and so
  are distinct isotypic components. Write
  $$
    h=c\,E_{n,i}\otimes u_0+\sum_kZ_k\otimes w_k
  $$
  with $Z_k\in\liep'^+$ and
  $w_k\in\tau_\sigma$ of weights $<\mu$. This is possible since the weights $e_n-e_k$ of $\liep'^+$ have multiplicity
  one: the second factor of the component of $h$ in $E_{n,k}\otimes\tau_\sigma$ has $\U(p)$-weight $\mu-e_i+e_k$, which
  is $\mu$ for $k=i$, with weight space $\CC u_0$, and $\neq\mu$ for $k\neq i$,
  hence $<\mu$ if the component is non-zero.
  We have $c\neq0$, since in $V_\nu\otimes V_\lambda$ the component in $V_\nu\otimes v_\lambda$ of a non-zero
  $\lieu$-invariant vector is non-zero, where $\lieu$ is the nilradical of the standard Borel subalgebra of
  $\gl_p$. Indeed, regard the vector
  as a $\lieu$-map $\varphi\colon V_\lambda^*\to V_\nu$, and let $v^*_{-\lambda}\in V_\lambda^*$ be the lowest weight
  vector with $v^*_{-\lambda}(v_\lambda)=1$, which vanishes on the weight spaces of weights $<\lambda$. The component
  of the vector in $V_\nu\otimes v_\lambda$ is $\varphi(v^*_{-\lambda})\otimes v_\lambda$. As $V_\lambda^*=U(\lieu)v^*_{-\lambda}$, $\varphi(v^*_{-\lambda})=0$
  forces $\varphi=0$. Since $E_{n,i}\otimes u_0$ has the weight of $h$, and the weight space of that weight in the
  $\gamma^{(i)}$-component is $\CC h$, its $\gamma^{(i)}$-component is $\langle E_{n,i}\otimes u_0,h\rangle h/\|h\|^2$, and
  $$
    \langle E_{n,i}\otimes u_0,h\rangle=\bar c\,\|E_{n,i}\|^2\|u_0\|^2\neq0.
  $$

  Hence $T=0$ if and only if $T(E_{n,i}\otimes u_0)=0$. Now
  $$
    T(E_{n,i}\otimes u_0)=P_{m'}P'_{\gamma^{(i)}}\pi(E_{n,i})x_F,
  $$
  and by Lemma \ref{ghl:lem:coeff}, applied to the component in $V_{m'}$ of $E_{n,i}\sum_{r\in J(m')}c_x(m'+e_r)\xi(m'+e_r,\mu)$, this is the left side of the reduced equation $(m',i)$ for
  $c_x$ times $\Phi(\Qhw(m';\mu-e_i))$. If $c_x$ satisfies $(m',i)$, then $T=0$ and $\langle x,\pi(Y)e\rangle=0$ for all
  $Y\in\liep'_\CC$, $e\in W'$. Conversely, the condition with $Y=E_{i,n}$, $Y^\dagger=-E_{n,i}$ and
  $e=\Phi(\Qhw(m';\mu-e_i))$ gives $\langle T(E_{n,i}\otimes u_0),e\rangle=0$, i.\,e.\ $(m',i)$.

  (c) Let $m'$ be as in section \ref{ghl:sec:reduced} and $x'=x(m')$, $z_r=x'_r+1$. By Lemma
  \ref{ghl:lem:matrixelement} and fact (d) of section \ref{ghl:sec:uniqueness}, the equations at $m'$ are
  $$
    \sum_{r\in J(m')}\frac{w_r}{z_r-y_i}=0,\quad i\in E(m'),
  $$
  with $w_r$ as there. Let $c=c_\bullet$ and $r\in J(m')$. For
  $r\leq q$, $\gamma_r(z_r)=x'_r-L_r=d_r(m')+1$, $d_r(m'+e_r)=d_r(m')+1$ and $d_s(m'+e_r)=d_s(m')$ $(s\neq r)$. For
  $r=n-1$, $\gamma_{n-1}=1$ and no $d_s$ changes. Therefore
  $$
    c(m'+e_r)\gamma_r(z_r)=\frac{\Delta_I(x'+e_r)}{\prod_{s\leq q}d_s(m')!}.
  $$
  In $\Delta_I$ the factors containing $r$ change by $(x'_r+1-x'_s)/(x'_r-x'_s)$, and $x'_r+1\neq x'_s$ as $x'$ and
  $x'+e_r$ are strictly decreasing. Hence
  $$
    w_r=K(m')\frac{\prod_{j\leq q}(z_r-y_j)}{\prod_{s\in I\smallsetminus\{r\}}(z_r-z_s)}=K(m')\Res_{t=z_r}R(t),
    \qquad K(m'):=\frac{\Delta_I(x')}{\prod_{s\leq q}d_s(m')!},
  $$
  with
  $$
    R(t)=\frac{\prod_{j\leq q}(t-y_j)}{\prod_{s\in I}(t-z_s)},
  $$
  whose poles are simple since the $z_s$ are distinct. Fix
  $i\in E(m')$. By fact (a), $I\smallsetminus J(m')\subset\{2,\dots,q\}$, and for $r=i'+1\in I\smallsetminus J(m')$ fact
  (c) gives $i'\notin E(m')$, so $i'\neq i$, and $z_r=y_{i'}$. The numerator $\prod_{j\neq i}(t-y_j)$ of
  $R(t)/(t-y_i)$ contains $t-y_{i'}$, and distinct such $r$ use distinct factors. At $t=y_i$ the function
  $R(t)/(t-y_i)$ is regular: $y_i\neq z_r$ for $r\in J(m')$ by fact (d), and $y_i\neq y_{i'}=z_{i'+1}$ otherwise.
  Its only poles are therefore the $z_r$, $r\in J(m')$, all simple, with residues $\Res_{z_r}R/(z_r-y_i)$, and as
  it is $O(t^{-2})$ (degrees $q-1$ and $q+1$), these residues sum to $0$. This is the equation $(m',i)$ for $c_\bullet$.
\end{proof}

\begin{lemma}\label{ghl:lem:sol}
  $\mathrm{Sol}=\CC v$, and $v\neq0$.
\end{lemma}

\begin{proof}
  We have $0\neq v\in\mathcal X$ by Lemma \ref{ghl:lem:X}, and $v\in\mathrm{Sol}$ by Lemma \ref{ghl:lem:WS} (c), (b),
  (a). If $x\in\mathrm{Sol}$, then $c_x$ satisfies the reduced equations by Lemma \ref{ghl:lem:WS} (a), (b), so
  $c_x\in\CC c_\bullet$ by Proposition \ref{ghl:prop:unique}, and $x\in\CC v$.
\end{proof}

\subsection{The closed span of $v$}\label{ghl:sec:span}

Let $\mathcal A$ be the commutant of $\rho(H')$, the algebra of bounded operators on $\hatpi$ commuting with
$\rho(H')$, and let $V$ be the closed $H'$-span of $v$.

\begin{proposition}\label{ghl:prop:Virred}
  $a v\in\CC v$ for every $a\in\mathcal A$, and $V$ is irreducible.
\end{proposition}

\begin{proof}
  The operator $a$ commutes with $\rho(K')$, hence with the two projections defining $\mathcal X$, so
  $a\mathcal X\subset\mathcal X$. Both $a$ and $a^*\in\mathcal A$ preserve $\rho^\infty$, as
  $h\mapsto\rho(h)ax=a\rho(h)x$, and commute there with $\rho(Y)$ and with $P'_\gamma$.
  For $x\in\mathrm{Sol}$ and $e\in P'_{\gamma^{(i)}}\rho^\infty$ we have $a^*e\in
  P'_{\gamma^{(i)}}\rho^\infty$, so $\langle ax,\rho(Y)e\rangle=\langle x,\rho(Y)a^*e\rangle=0$.
  Hence $\mathcal A$ preserves $\mathrm{Sol}=\CC v$ (Lemma \ref{ghl:lem:sol}).
  Let $Q\in\mathcal A$ be the projection onto $V$. The commutant of
  $\rho|_V$ is $Q\mathcal AQ$. If $Q'\in Q\mathcal AQ$ is a projection, then $Q'v=cv$ with $c^2=c$. If
  $c=1$, $Q'\rho(h)v=\rho(h)v$ for all $h$ and $Q'=Q$, and if $c=0$, $Q'=0$. A von Neumann algebra without
  non-trivial projections is $\CC$, so $V$ is irreducible by Schur's lemma.
\end{proof}

\begin{lemma}\label{ghl:lem:Vds}
  $V$ is a discrete series representation of $H'$.
\end{lemma}

\begin{proof}
  This is the argument of Remark \ref{ghl:rmk:discrete}. As $\pi_K$ is dense and $Q\neq0$, there is $y\in\pi_K$ with
  $Qy\neq0$. By Lemma \ref{ghl:lem:L2}, $h\mapsto\langle\rho(h)y,y\rangle$ lies in $L^2(H')$, and by Lemma
  \ref{ghl:lem:closed} the closed $H'$-span $W$ of $y$ has an isometric $H'$-embedding $U\colon W\to L^2(H')$. The
  orthogonal projection $P_W$ restricts to an $H'$-map $V\to W$, non-zero since $\langle P_WQy,y\rangle=\|Qy\|^2$.
  Therefore $U\circ P_W|_V$ is a non-zero $H'$-map from the irreducible unitary $V$ to $L^2(H')$, hence a multiple
  of an isometry.
\end{proof}

\begin{lemma}\label{ghl:lem:classif}
  $V\cong\pi_{A''(q')}$ for some $A''\in(\ZZ+\frac p2)^{p+1}$ strictly decreasing and $0\leq q'\leq p$, with
  $\sum_jA''_j=\sum_iA'(q)_i$.
\end{lemma}

\begin{proof}
  Let $G_1=\SU(n-2,1)$. It is connected, semisimple, linear and of equal rank, $K_1=K'\cap G_1$ is connected, and
  $H'=Z\cdot G_1$ with $Z=\U(1)\cdot\one$ the centre and $Z\cap G_1$ finite. The centre $Z$ acts on the
  irreducible $V$ by a character $z\mapsto z^c$, so closed $G_1$-invariant subspaces are $H'$-invariant and
  $V|_{G_1}$ is irreducible. Furthermore, $|\langle\rho(zg)x,y\rangle|=|\langle\rho(g)x,y\rangle|$ for $x,y\in V$, with Haar measure
  on $H'$ the image of that on $Z\times G_1$, so $V|_{G_1}$ is square integrable. By \cite[p.~536]{schmid75} the
  character of $V|_{G_1}$ is $\Theta_\lambda$ for a non-singular $\lambda\in\Lambda+\rho$, unique up to $W_{K_1}$ by
  \cite[p.~130]{hechtschmid}; and as an irreducible unitary representation is determined by its distribution
  character, $V|_{G_1}\cong\pi^{G_1}_\lambda$.

  By \cite[Theorem (1.3)]{schmid75} (a), $V|_{G_1}$ contains the $K_1$-type $\lambda-\rho_c+\rho_n$. It lies in a
  $K'$-type $\delta=(\delta_1,\dots,\delta_p;\delta_0)$ of $V$, and $\sum\delta=c$. Put $\Lambda'':=\delta+\rho_c-\rho_n$,
  with $\rho_c$, $\rho_n$ those of the positive system of $\lambda$, in the coordinates of $T'$. Every entry of
  $\rho_c-\rho_n$ lies in $\ZZ+\frac p2$ ($\rho_{c,i}=\frac{p+1}2-i$, $\rho_{c,0}=0$, $\rho_{n,i}=\pm\frac12$,
  $\rho_{n,0}=-\frac12\sum_i(\pm1)$), so $\Lambda''\in(\ZZ+\frac p2)^{p+1}$. Its image modulo $\CC(1,\dots,1)$ is
  $\lambda$, and $\sum\Lambda''=\sum\delta=c$, as roots have coordinate sum $0$. These two data determine $\Lambda''$,
  which is therefore unique up to $W_{K'}=S_p$, and regular since $\lambda$ is non-singular. Sorting its first $p$
  entries, $\Lambda''=A''(q')$ in the form \eqref{ghl:eq:Aq} with $A''$ strictly decreasing and $q'$ the number of
  compact entries larger than the non-compact one. A unitary $G_1$-isomorphism $V|_{G_1}\to\pi_{A''(q')}|_{G_1}$
  commutes with $Z$, which acts on both by $z^c$, hence is an $H'$-isomorphism.

  Finally $Z\subset T'$ acts on $v\in\mathcal X$ by $z^{\sum\mu+\tau_0}$, and $\sum\mu+\tau_0=\sum\tau_\sigma=
  \sum_iA'(q)_i$, since $\rho_n-\rho_c$ has coordinate sum $0$. Hence $c=\sum_jA''_j=\sum_iA'(q)_i$.
\end{proof}

\begin{lemma}\label{ghl:lem:Vschmid}
  $v$ is a smooth vector of $V$, and

  \noindent {\rm(i)} $P'_{\gamma^{(i)}}\rho(Y)v=0$ for $Y\in\liep'_\CC$ and $1\leq i\leq q$;

  \noindent {\rm(ii)} for $q<i\leq n-2$, $\gamma^{(i)}$ is not a $K'$-type of $V$.
\end{lemma}

\begin{proof}
  By Lemma \ref{ghl:lem:classif} and Proposition \ref{ghl:prop:Hprime} (i), $V(\tau_\sigma):=P'_{\tau_\sigma}V$ is
  finite-dimensional. The subspace $P'_{\tau_\sigma}V^\infty$ is dense in it, hence equal to it, so
  $V(\tau_\sigma)\subset V^\infty$, and $v\in V(\tau_\sigma)$. Smooth vectors of $V$ are smooth in $\rho$.

  \noindent (i) For $e\in P'_{\gamma^{(i)}}\rho^\infty$ and
  $Y\in\liep'_\CC$, \eqref{ghl:eq:WS} for $v$ (Lemma \ref{ghl:lem:sol}) and
  $\langle v,\rho(Y)e\rangle=-\langle\rho(\bar Y)v,e\rangle$, where $\bar Y\in\liep'_\CC$ is the conjugate of $Y$
  with respect to $\lieh'$, give $\langle P'_{\gamma^{(i)}}\rho(\bar Y)v,e\rangle=0$. Such $e$ are dense in
  $P'_{\gamma^{(i)}}\hatpi$, so $P'_{\gamma^{(i)}}\rho(\bar Y)v=0$, and $\bar Y$ runs through $\liep'_\CC$.

  \noindent (ii) Here $\gamma^{(i)}$ is not a $K'$-type of $\rho$ (section \ref{ghl:sec:weak}), and
  $V\subset\rho$.
\end{proof}

\subsection{Identification}\label{ghl:sec:Id}

We now combine the $K'$-type and transition description of Proposition \ref{ghl:prop:Hprime} with the
conditions of Lemmas \ref{ghl:lem:classif} and \ref{ghl:lem:Vschmid} to identify $V$ with $\hatsigma$
(Lemma \ref{ghl:lem:Id}).

\begin{lemma}\label{ghl:lem:Id}
  Let $V\subset\rho$ be an irreducible closed subrepresentation, $V\cong\pi_{A''(q')}$ as in Lemma
  \ref{ghl:lem:classif}, with $\sum_jA''_j=\sum_iA'(q)_i$, and assume
  \begin{itemize}
    \item[(H1)] $0\neq v\in V(\tau_\sigma)$ (then $v$ is a smooth vector of $V$, $V(\tau_\sigma)$ being
      finite-dimensional);
    \item[(H2)] $P'_{\gamma^{(i)}}\rho(Y)v=0$ for $Y\in\liep'_\CC$ and $1\leq i\leq q$;
    \item[(H3)] for $q<i\leq n-2$, $\gamma^{(i)}$ is not a $K'$-type of $V$.
  \end{itemize}
  Then $q'=q$ and $A''(q')=A'(q)$, i.\,e.\ $V\cong\hatsigma$.
\end{lemma}

\begin{proof}
  Recall $y_i=x_i(\mu)$, i.\,e.\ $y_i=L_i+1$ $(i\leq q)$, $y_i=L_{i+1}$ $(q<i\leq p)$. By (H1) and Proposition
  \ref{ghl:prop:Hprime} (i), $x(\mu)\in R'$, i.\,e.\ $y_i\in J_i$ for all $i$, and $\tau_\sigma$ has multiplicity one in $V$,
  so $V(\tau_\sigma)\cong\tau_\sigma$ is the $K'$-span of $v$.

  {\em From $v$ to $V(\tau_\sigma)$.} For $i\leq q$ the map
  $$
  \liep'_\CC\otimes V(\tau_\sigma)\to V(\gamma^{(i)}),\qquad Y\otimes w\mapsto P'_{\gamma^{(i)}}\rho(Y)w,
  $$
  is $K'$-equivariant. By (H2) it vanishes on $\liep'_\CC\otimes v$, hence
  on $\liep'_\CC\otimes\rho(k)v=\rho(k)(\mathrm{Ad}(k)^{-1}\liep'_\CC\otimes v)$ for all $k\in K'$, hence on
  $\liep'_\CC\otimes V(\tau_\sigma)$. Therefore $P'_{\gamma^{(i)}}\rho(\liep'_\CC)P'_{\tau_\sigma}=0$
  on $V_{K'}$.

  {\em The case $i\leq q$.} The central character fixes the $\U(1)$-coordinates: $\mathrm{tot}'-\sum\mu=\tau_0$. If
  $x(\mu)-e_i\in R'$, then $\gamma^{(i)}=(\mu-e_i;\tau_0+1)$ is a $K'$-type of $V$, and Proposition
  \ref{ghl:prop:Hprime} (ii) with $m=\mu-e_i$ gives $P'_{\gamma^{(i)}}\rho(\liep'_\CC)P'_{\tau_\sigma}\neq0$, against the
  above. Hence $x(\mu)-e_i\notin R'$. As $R'$ is a product of intervals in the coordinates $x_i$ and $x(\mu)\in R'$,
  this
  means $y_i=\min J_i$.

  {\em The case $q<i\leq p$.} If $x(\mu)+e_i\in R'$, then $\gamma^{(i)}=(\mu+e_i;\tau_0-1)$ is a $K'$-type of $V$, against
  (H3). Hence $y_i=\max J_i$.

  We have $y_q-y_{q+1}=L_q+1-L_{q+2}\geq2$ if $q\leq p-1$. The four cases are:
  \begin{itemize}
    \item[(a)] If $q'<q=p$, then $J_p=(-\infty,L''_{p+1}]$, as $p>q'$ and $L''_{p+2}=-\infty$, so $y_p=\min J_p$ is
      impossible.
    \item[(b)] If $q'<q\leq p-1$, then as $q$, $q+1>q'$ we have $y_q=\min J_q=L''_{q+2}+1$ and
      $y_{q+1}=\max J_{q+1}=L''_{q+2}$, so $y_q-y_{q+1}=1$.
    \item[(c)] If $q'>q$, so that $q\leq p-1$, then as $q$, $q+1\leq q'$ we have $y_q=\min J_q=L''_q+1$ and
      $y_{q+1}=\max J_{q+1}=L''_q$, so $y_q-y_{q+1}=1$.
    \item[(d)] If $q'=q$, then $y_i=L''_i+1$ $(i\leq q)$ and $y_i=L''_{i+1}$ $(q<i\leq p)$, i.\,e.\ $L''_j=L_j$ for
      $j\leq q$ and $q+2\leq j\leq p+1$.
  \end{itemize}
  The cases (a) to (c) are contradictory. In the case (d), the entries of $A'(q)$ minus $\frac p2$ are
  $L^\sigma_j=L_j$ $(j\neq q+1)$ and $L^\sigma_{q+1}=L_{q+1}+1$ ($j\leq n-1$), strictly decreasing by
  \eqref{ghl:eq:gap}, with the non-compact entry in position $q+1$, and $\sum_jL''_j=\sum_jL^\sigma_j$ gives
  $L''_{q+1}=L_{q+1}+1$. We conclude that $A''(q')=A'(q)$.
\end{proof}

\begin{theorem}\label{ghl:thm:J}
  Let $1\leq q\leq n-2$ and let $A$ satisfy \eqref{ghl:eq:gap}. Then the closed $H'$-span $V$ of
  $$
    v=\sum_{\delta\in S}c_\delta\Phi(\xi_\delta)
  $$
  is an irreducible subrepresentation of $\hatpi|_{H'}$ isomorphic to
  $\hatsigma$. In particular \hypref{ghl:hyp:J}{J} holds, and $v=ju$ for an isometric $H'$-embedding $j\colon\hatsigma\to\hatpi$ and a
  $K'$-highest weight vector $u\in\tau_\sigma$.
\end{theorem}

\begin{proof}
  By Proposition \ref{ghl:prop:Virred} and Lemmas \ref{ghl:lem:Vds}, \ref{ghl:lem:classif}, $V$ is irreducible and
  $V\cong\pi_{A''(q')}$ with $\sum_jA''_j=\sum_iA'(q)_i$. (H1) holds as $v\in V\cap\mathcal X$, $v\neq0$, and (H2) and (H3)
  are Lemma \ref{ghl:lem:Vschmid}. By Lemma \ref{ghl:lem:Id}, $V\cong\hatsigma$, and an isomorphism of irreducible
  unitary representations is a multiple of a unitary one. This gives \hypref{ghl:hyp:J}{J}. By Lemma \ref{ghl:lem:P} (i) every
  isometric $H'$-embedding $j$ has image $V$. By Proposition \ref{ghl:prop:Hprime} (i) for $\sigma$, $\tau_\sigma$ has
  multiplicity one in $\sigma$, so $V(\tau_\sigma)=j(\tau_\sigma)$, and its $T'$-weight line of weight $(\mu;\tau_0)$,
  which contains $v$, is $j(\CC u_0)$. Hence $v=ju$ for some $u\in\CC u_0$, and $\|u\|=\|v\|$, $j$ being isometric.
\end{proof}

\begin{remark}\label{ghl:rmk:J}\leavevmode\par
  \noindent {\rm(1)} The closed span $V$ is irreducible before it is known to be $\sigma$. The uniqueness of the solutions of the
  Schmid conditions in $\mathcal X$, a statement about $\pi_K$, is what makes this so. The infinitesimal character would not suffice in place of Lemma \ref{ghl:lem:Id}: it does not separate
  $\sigma$ from the other members of its packet.

  \smallskip
  \noindent {\rm(2)} With $u$ as in Theorem \ref{ghl:thm:J} and $f:=P_\tau v\in\tau$, Lemma \ref{ghl:lem:P} (iii) gives
  $I(f,f')=d(\sigma)^{-1}\|P_\tau v\|^4>0$ for $f'=\bar u$, so the integral \eqref{ghl:eq:I} on the minimal types is
  positive for this pair.
\end{remark}
 \section{Proof of Theorem \ref{ghl:thm:A}}\label{ghl:sec:proofA}

We prove Theorem \ref{ghl:thm:A}. Section \ref{ghl:sec:q0} treats $q=0$, and in section
\ref{ghl:sec:factorization} we factor $I^*$ (Corollary \ref{ghl:cor:Cinf}). For $q\geq1$, Theorem \ref{ghl:thm:J}
and Lemma \ref{ghl:lem:X} give the components of $ju$ and the value of $s(q)$ (Theorem \ref{ghl:thm:schmid}, section
\ref{ghl:sec:ju}). In section \ref{ghl:sec:conditions} we derive the reduced equations from the Schmid conditions
for $ju$, which gives a second proof of Theorem \ref{ghl:thm:schmid} for a given $j$ (Remark
\ref{ghl:rmk:Jholds} (3)).

\subsection{The case $q=0$}\label{ghl:sec:q0}

For $q=0$ we prove \hypref{ghl:hyp:J}{J} directly. By \eqref{ghl:eq:tau} and \eqref{ghl:eq:tausigma}, $\tau_i=L_{i+1}+i-1$
$(i\leq n-1)$, $\tau_0=L_1+n-1$, and $\tau_\sigma=(\tau_1,\dots,\tau_{n-2};\tau_0)$. Since
$\mu=(\tau_1,\dots,\tau_{n-2})$ interlaces $(\tau_1,\dots,\tau_{n-1})$, $\tau_\sigma$ occurs in $\tau|_{K'}$,
with multiplicity one. Put $\liep^-:=\langle E_{in}:i\leq n-1\rangle$, $\liep^+:=\langle E_{ni}:i\leq n-1\rangle$
and $\liep'^\pm:=\liep^\pm\cap\lieh'_\CC$.

\begin{proposition}\label{ghl:prop:Jq0}
  Let $q=0$. Then $M:=U(\lieh'_\CC)\tau_\sigma\subset\pi_K$ is irreducible and isomorphic to
  $\sigma_K$. Its closure $W\subset\hatpi$ is $H'$-stable, and with the restricted inner product
  $W\cong\hatsigma$. In particular \hypref{ghl:hyp:J}{J} holds, $j(\tau_\sigma)=\tau_\sigma\subset\tau$, and $s(0)=1$.
\end{proposition}

\begin{proof}
  The $K$-types of $\pi$ are $(m;\mathrm{tot}-\sum m)$ with $m\in R_0$, i.\,e.\ $m_i\leq\tau_i$
  (Corollary \ref{ghl:cor:ktypes}), so their $\U(1)_{e_n}$-weight is at least $\tau_0$. Since $E_{in}$
  lowers this weight by one, we have $\liep^-\tau=0$. Likewise every $K'$-type of $\sigma$ has
  $\U(1)_{e_n}$-weight at least $\tau_0$ (it lies in $\tau_\sigma+\NN\Delta^+(\sigma)$ by Lemma
  \ref{ghl:lem:Kprime}), so $\liep'^-$ kills $\tau_\sigma\subset\sigma_K$. By Poincar\'e--Birkhoff--Witt,
  $M=U(\liep'^+)\tau_\sigma$. The module $M$ is graded by the eigenvalue $\tau_0+d$ of $E_{nn}$, which is
  central in $\liek'_\CC$. Its graded pieces $M_d$ are finite-dimensional, $\liep'^+$ raises the degree, and
  $M_0=\tau_\sigma$ is $K'$-irreducible.

  We show that $M$ is irreducible. A submodule $N$ is $K'$-stable, hence graded. The inner product of
  $\hatpi$ is positive and invariant on $M$, and distinct degrees are orthogonal, so
  $M=N\oplus(N^\perp\cap M)$ with
  both summands submodules. Then $\tau_\sigma=N_0\oplus(N^\perp\cap M)_0$ lies in $N$ or in
  $N^\perp\cap M$, and as $\tau_\sigma$ generates $M$, $N=M$ or $N=0$.

  Next we show that $M\cong\sigma_K$. The generalized Verma module
  $U(\lieh'_\CC)\otimes_{U(\liek'_\CC+\liep'^-)}\tau_\sigma$ has a unique irreducible quotient, since the
  sum of its submodules with zero degree-$0$ part is its unique maximal proper submodule. Both $M$ and
  $\sigma_K$ ($\sigma_K$ is irreducible and generated by $\tau_\sigma$, which is killed by $\liep'^-$) are
  irreducible quotients, hence isomorphic.

  The $K$-finite vectors of $\pi$ are analytic and $H'$ is connected, so $W$ is $H'$-stable. Its
  $K'$-finite vectors are $M$ (a $K'$-type with $\U(1)_{e_n}$-character $\tau_0+d$ occurs in $M$ only in $M_d$,
  so the $K'$-isotypic parts of $M$ are finite-dimensional). Therefore $W$ is an
  irreducible unitary representation of $H'$ infinitesimally equivalent to $\sigma$, hence unitarily
  equivalent to it. Finally $j(\tau_\sigma)=\tau_\sigma\subset\tau$ gives $P_\tau ju=ju$ and $s(0)=1$.
\end{proof}

\subsection{The factorization of $I^*$}\label{ghl:sec:factorization}

\begin{corollary}\label{ghl:cor:Cinf}
  Let $\pi=\pi_{A(q)}$, $\sigma=\pi_{A'(q)}$, with \eqref{ghl:eq:gap}. Then there is an isometric
  $H'$-embedding $j\colon\hatsigma\to\hatpi$, unique up to a phase. If $P_\tau ju\neq0$ for
  $0\neq u\in\tau_\sigma$, then for $f\in\tau$ and $f'$ in the minimal $K'$-type of $\pi'$
  $$
    I^*(f,f')=\frac{1}{d(\sigma)\,s(q)}\,R^\circ(f,f').
  $$
\end{corollary}

\begin{proof}
  By Proposition \ref{ghl:prop:Jq0} $(q=0)$ and Theorem \ref{ghl:thm:J} $(q\geq1)$, \hypref{ghl:hyp:J}{J} holds.
  A non-zero bounded $H'$-map $\hatsigma\to\hatpi$ is a multiple of an isometry by Schur's lemma for
  the irreducible $\sigma$, and Lemma \ref{ghl:lem:P} (i) gives uniqueness. For $f\in\tau$,
  $\langle f,ju\rangle=\langle f,P_\tau ju\rangle=\langle f,w\rangle$, so by Lemma \ref{ghl:lem:P} (iii)
  $$
    I^*(f,f')=d(\sigma)^{-1}|\langle f,w\rangle|^2/(\|f\|^2\|u\|^2)=d(\sigma)^{-1}s(q)^{-1}R^\circ(f,f').
  $$
\end{proof}

We have $P_\tau ju\neq0$ for $q=0$ by Proposition \ref{ghl:prop:Jq0}, and for $q\geq1$ by Theorem \ref{ghl:thm:schmid} below.

\begin{remark}\label{ghl:rmk:Jholds}\leavevmode\par
  \noindent {\rm(1)} The statement \hypref{ghl:hyp:J}{J} is known. By the theorem of He \cite[Theorem 1.1]{he},
  the discrete spectrum of the restriction of a discrete series of $\U(p,q)$ to $\U(p-1,q)$ is the Hilbert sum of
  the discrete series whose Harish-Chandra parameters satisfy the interlacing relation of \cite[Definition
  1]{he} with the parameter of the given representation. For $\pi=\pi_{A(q)}$ and $\sigma=\pi_{A'(q)}$ we have
  $A_i>A_i-\frac12>A_{i+1}$ for $i\leq n-1$, $i\neq q+1$, and for the non-compact entries
  $$
    A_{q+1}+\tfrac12>A_{q+1}>A_{q+2},\qquad A_q-\tfrac12>A_{q+1}+\tfrac12\quad\text{if }q\geq1,
  $$
  the last inequality being \eqref{ghl:eq:gap}. In the notation of \cite{he} the sequence of signs is
  $(+\oplus)^q\,\ominus\,-\,(+\oplus)^{n-q-2}\,+$, which contains only the admissible pairs $(+\oplus)$,
  $(\oplus+)$, $(\oplus\ominus)$, $(\ominus-)$ and $(-+)$. Hence $\sigma$ occurs in $\hatpi|_{H'}$ for all
  $0\leq q\leq n-2$ under \eqref{ghl:eq:gap}. The theorem of Beuzart-Plessis \cite[Theorem
  12.4.1]{beuzartplessis} shows that among the pairs of representations of the pure inner forms with a given
  tempered Langlands parameter exactly one admits a non-zero invariant form on the smooth vectors. It does not determine
  this pair. The proof of \cite{he} proceeds by the theta correspondence and an induction on the rank, and does
  not give the components of $ju$ in the $K$-types of $\pi$. Theorem \ref{ghl:thm:J} proves \hypref{ghl:hyp:J}{J}
  for $q\geq1$ within the model of section \ref{ghl:sec:gt}, and identifies $ju$ with the formal Schmid vector,
  whose components are explicit. The value of $s(q)$ rests on the latter.

  \smallskip
  \noindent {\rm(2)} In the setting of \cite[Conjecture 4.15]{ghl}, \hypref{ghl:hyp:J}{J} also follows from
  \cite[Theorem 4.14]{ghl} by Lemma \ref{ghl:lem:P} (ii) and \cite[Proposition 6.12]{ghl}. That route depends on \cite[Corollary 1.3 (b)]{harris14}, whose hypothesis
  \cite[Hypothesis 1.1]{harris14} is, at $\vz$, equivalent to \hypref{ghl:hyp:J}{J}. By (1), \cite[Theorem 1.1]{he} gives \cite[Hypothesis
  1.1]{harris14} at $\vz$ for the representations of \cite[Theorem 4.14]{ghl}.

  \smallskip
  \noindent {\rm(3)} If \hypref{ghl:hyp:J}{J} is taken from \cite[Theorem 1.1]{he}, then Theorem \ref{ghl:thm:A}
  follows without section \ref{ghl:sec:converse}. Indeed, the embedding $j$ of Corollary \ref{ghl:cor:Cinf} then
  exists, Propositions \ref{ghl:prop:ju} and \ref{ghl:prop:unique} give $P_\tau ju\neq0$ and
  $s(q)=\sum_{\delta\in S}t_\delta$ (Remark \ref{ghl:rmk:schmidroute}), and Corollary \ref{ghl:cor:summation}
  gives \eqref{ghl:eq:P}. This uses sections \ref{ghl:sec:lemmaP} to \ref{ghl:sec:summation}, Corollary \ref{ghl:cor:Cinf} and section
  \ref{ghl:sec:conditions}. With section \ref{ghl:sec:converse} the proof of Theorem \ref{ghl:thm:A} does not depend
  on \cite{he}.
\end{remark}

\subsection{The components of $ju$}\label{ghl:sec:ju}

Let $1\leq q\leq n-2$. The notation is that of sections \ref{ghl:sec:schmid} and \ref{ghl:sec:converse}. Let $j$
be as in Corollary \ref{ghl:cor:Cinf}, let $u\in\tau_\sigma$ be a $K'$-highest weight vector, $u\neq0$, and let $v$
be as in Lemma \ref{ghl:lem:X}.

\begin{theorem}\label{ghl:thm:schmid}
  Let $1\leq q\leq n-2$ and $A$ with \eqref{ghl:eq:gap}. Then there is $\lambda\in\CC^\times$
  with
  $$
    P_\delta ju=\lambda c_\delta\Phi(\xi(\delta,\mu))
  $$
  for $\delta\in S$ and $P_\delta ju=0$ otherwise. In
  particular $P_\tau ju\neq0$ and
  $$
    s(q)=\frac{\|ju\|^2}{\|P_\tau ju\|^2}=\sum_{\delta\in S}t_\delta.
  $$
\end{theorem}

\begin{proof}
  By Theorem \ref{ghl:thm:J} there are an isometric $H'$-embedding $j_0\colon\hatsigma\to\hatpi$ and a
  $K'$-highest weight vector $u'\in\tau_\sigma$ with $v=j_0u'$. By Lemma \ref{ghl:lem:P} (i), $j=\varepsilon j_0$
  with $|\varepsilon|=1$, both maps being isometric. The $K'$-highest weight vectors of $\tau_\sigma$ form a line,
  and $u'\neq0$ as $v\neq0$, so $u=\alpha u'$ with $\alpha\in\CC^\times$. Hence $ju=\lambda v$ with
  $\lambda:=\varepsilon\alpha\neq0$. By Lemma \ref{ghl:lem:X} (ii) the series $v=\sum_{\delta\in S}c_\delta\Phi(\xi_\delta)$
  converges in norm, and its terms lie in the distinct $K$-types $\delta\in S$, so $P_\delta ju=\lambda
  c_\delta\Phi(\xi_\delta)$ for $\delta\in S$ and $P_\delta ju=0$ otherwise. As $c_\tau=\Delta_I(x(\tau))\neq0$,
  $P_\tau ju\neq0$. Finally $\|P_\tau ju\|^2=|\lambda|^2|c_\tau|^2N_\pi(\xi_\tau)$, and by Lemma \ref{ghl:lem:X} (ii)
  and Corollary \ref{ghl:cor:summation}
  $$
    \|ju\|^2=|\lambda|^2\|v\|^2=|\lambda|^2|c_\tau|^2N_\pi(\xi_\tau)\sum_{\delta\in S}t_\delta,
  $$
  which gives the claim.
\end{proof}

\begin{proof}[Proof of Theorem \ref{ghl:thm:A}]
  For $q=0$, \hypref{ghl:hyp:J}{J} holds and $s(0)=1$ by Proposition \ref{ghl:prop:Jq0}. For $q\geq1$, \hypref{ghl:hyp:J}{J} holds by Theorem \ref{ghl:thm:J}, and
  Theorem \ref{ghl:thm:schmid} gives $P_\tau ju\neq0$ and $s(q)=\sum_\delta t_\delta$, which is
  \eqref{ghl:eq:P} by Corollary \ref{ghl:cor:summation}. By Corollary \ref{ghl:cor:Cinf} we obtain (ii), and
  \eqref{ghl:eq:fdeg} gives the second expression for $C_\infty(q)$. Finally $\dim F_{A'(q)}\in\ZZ_{>0}$.
\end{proof}

\begin{remark}\label{ghl:rmk:sharp}
  The condition $A_q-A_{q+1}\geq2$ is used three times in the computation of $s(q)$: for $\tau_\sigma$ in \eqref{ghl:eq:tausigma}, in the
  tail estimate of Lemma \ref{ghl:lem:base} (the bound $a_p-D-1-m\geq1$ propagates to $a_q\geq2$), and in
  $\Omega_{q,0}=\prod(a_r-1)^{-1}$. It is also used in the proof of \hypref{ghl:hyp:J}{J}, for $\sigma$ to be a discrete series and in Lemma \ref{ghl:lem:Id}. At $a_q=1$ the parameter $A'(q)$ is singular, $\sigma$ is not a discrete
  series, \eqref{ghl:eq:P} has a pole, and the formal series \eqref{ghl:eq:tdelta} diverges: for $q=1$,
  $a_1=1$, its terms are
  $$
    \frac{(2k+B+2)(B+1)}{(B+2)(k+1)(k+B+1)}\sim\frac{2(B+1)}{(B+2)k}.
  $$
  No condition on the other
  gaps is needed: a unit gap $A_r-A_{r+1}=1$ $(r\neq q)$ empties the range of $d_r$, matched by the missing
  equation for the non-dominant $\mu-e_r$.
\end{remark}

\subsection{The Schmid conditions for $ju$}\label{ghl:sec:conditions}

We derive the reduced equations of section \ref{ghl:sec:reduced} from the Schmid conditions for $ju$
(Proposition \ref{ghl:prop:ju}). Together with Proposition \ref{ghl:prop:unique} this gives a second proof of
Theorem \ref{ghl:thm:schmid} (Remark \ref{ghl:rmk:schmidroute}).

Let $1\leq q\leq n-2$, and let $j$ be as in Corollary \ref{ghl:cor:Cinf}. It exists by Theorem \ref{ghl:thm:J}, and also by \cite[Theorem 1.1]{he} (Remark \ref{ghl:rmk:Jholds} (1)).
Fix a $K'$-highest weight vector
$u\in\tau_\sigma$, $u\neq0$, and put $w_\delta:=P_\delta ju\in\pi_K(\delta)$.

Since $j$ is $K'$-equivariant, $ju$ is a $K'$-highest weight vector of weight $(\mu;\tau_0)$.
As $P_\delta$ commutes with $\pi(K')$, $w_\delta$ is $0$ or a $K'$-highest weight vector of that weight in
$\pi_K(\delta)$. By Corollary \ref{ghl:cor:ktypes} each $K$-type of $\pi$ has multiplicity one, so by the
multiplicity free branching such a vector exists only if $\delta=(m;m_0)$ with $m$ interlacing $\mu$ and
$m_0=\tau_0$, and is then unique up to a scalar. Moreover $m\in\Rq$ (Corollary \ref{ghl:cor:ktypes}) and
$\sum m=\sum\tau$, because $\sum m+m_0=\mathrm{tot}=\sum\tau+\tau_0$. By Lemma \ref{ghl:lem:support},
$\delta\in S$. Hence
\begin{equation}\label{ghl:eq:wdelta}
  w_\delta=c(\delta)\,\Phi(\xi(\delta,\mu))\quad(\delta\in S),\qquad w_\delta=0\quad(\delta\notin S),
\end{equation}
for a function $c\colon S\to\CC$.

Now $u\in\sigma_K$ is $H'$-smooth and $j(\sigma(h)u)=\pi(h)ju$, so $ju$ is $H'$-smooth
and $\pi(X)ju=j\sigma(X)u$ for $X\in\lieh'_\CC$. By Lemma \ref{ghl:lem:Kprime}, $P'_{\gamma^{(i)}}\sigma(X)u=0$,
and $P'_\gamma j=jP'_\gamma$, hence
\begin{equation}\label{ghl:eq:schmidcond}
  P'_{\gamma^{(i)}}\,\pi(X)\,ju=0\qquad(X\in\liep'_\CC,\ 1\leq i\leq n-2).
\end{equation}

\begin{proposition}\label{ghl:prop:ju}
  The function $c$ of \eqref{ghl:eq:wdelta} satisfies all reduced equations.
\end{proposition}

\begin{proof}
  Let $(m',i)$ be a reduced equation and $e:=\Phi(\Qhw(m';\mu-e_i))\in\pi_K(m')$. It lies in the $K'$-type $\gamma^{(i)}=(\mu-e_i;\tau_0+1)$.
  By \eqref{ghl:eq:schmidcond} with $X=E_{n,i}$, $\langle\pi(E_{n,i})ju,e\rangle=0$. By Lemma
  \ref{ghl:lem:termwise} and \eqref{ghl:eq:wdelta}, $0=\langle\Phi(x),e\rangle$, where $x\in V_{m'}$ is the component in $V_{m'}$ of
  $$
    \sum_{\delta}c(\delta)E_{n,i}\xi(\delta,\mu),
  $$
  the sum running over $\delta=m'+e_r\in S$. As $x$ has weight $\mu-e_i$, Lemma \ref{ghl:lem:coeff} shows
  that $\langle\Phi(x),e\rangle$ is the coefficient of $\Qhw(m';\mu-e_i)$ in $x$ times $\|e\|^2\neq0$.
\end{proof}

\begin{remark}\label{ghl:rmk:schmidroute}
  Propositions \ref{ghl:prop:ju} and \ref{ghl:prop:unique} give a proof of Theorem \ref{ghl:thm:schmid} which uses
  the existence of $j$, but not that $ju$ is a multiple of $v$, and hence no result of section
  \ref{ghl:sec:converse} if $j$ is given (Remark \ref{ghl:rmk:Jholds} (3)). By these propositions, $c=\lambda c_\bullet$ for some
  $\lambda\in\CC$. Since $ju\neq0$ and $\hatpi$ is the Hilbert sum of its
  $K$-isotypic subspaces, some $w_\delta\neq0$ and $\lambda\neq0$. As
  $c_\tau=\Delta_I(x(\tau))\neq0$, $P_\tau ju\neq0$. By Parseval and Lemma \ref{ghl:lem:normS},
  $$
    \|ju\|^2=\sum_{\delta\in S}|\lambda|^2|c_\delta|^2N_\pi(\xi_\delta)=
    |\lambda|^2\kappa\sum_\delta|c_\delta|^2\mathcal N(\xi_\delta),
  $$
  and dividing by $\|P_\tau ju\|^2=|\lambda|^2\kappa|c_\tau|^2\mathcal N(\xi_\tau)$ gives $s(q)=\sum_{\delta\in S}t_\delta$.
\end{remark}
 \section{Consequences for the conjecture of Grobner, Harris and Lin}\label{ghl:sec:consequences}

We return to the global setting of \cite[\S\S3.3, 4.2]{ghl} recalled in the introduction, and write
$H_v$, $H'_v$, $\pi_v=\pi(q)_v$, $\pi'_v=\pi'(q)_v$ for the local data at $v\in\Sinf$. For $v\neq\vz$,
$H_v=\U(n)$ and $H'_v=\U(n-1)$ are compact, $\pi_v$ has infinitesimal character $A_v=A_{\lambda,v}$, and
$\pi'_v$ is the dual of the representation $\sigma_v$ of $\U(n-1)$ with infinitesimal character
$\Lambda'_v=(A_{v,i}-\frac12)_{i\leq n-1}$ \cite[\S4.2]{ghl}. At $\vz$, $\sigma_{\vz}=\sigma$ and
$\Lambda'_{\vz}=A'(q)$.

\subsection{Measures}\label{ghl:sec:globalmeasures}

We use the measures of \cite[p.~28]{ghl}. There the measure on $H(\bfA_{F^+})$ is normalized as in
\cite[\S5.2, p.~2043]{harris13a}, with
$$
\vol(\U(V_1)(F^+)\backslash\U(V_1)(\bfA_{F^+}))=1,
$$
and this defines the measures on $H'(\bfA_{F^+})$ and $H''(\bfA_{F^+})$. In \cite{harris13a} the measures on
$H$ and on $H'$ are both Tamagawa measures, and for any factorization
the local measure at a finite place takes rational values, and gives volume $1$ to a hyperspecial maximal compact
subgroup where the group is unramified. We take $dh'$ to be the Tamagawa measure of
$H'$, factored in this way. Its local factors at the finite places take rational values on all compact open
subgroups, since $\vol(U)/\vol(U\cap U_0)=[U:U\cap U_0]$ for compact open $U$, $U_0$. If $dh'$ is instead obtained
as the quotient of the Tamagawa measure on $H''=H'\times\U(V'_1)$ by the measure on $\U(V'_1)\cong\U(V_1)$ of volume $1$,
it is twice the Tamagawa measure, as the Tamagawa number of $\U(V_1)$ is $2$, and its local factors at the finite
places have the same properties if the measure on $\U(V'_1)$ is factored in the same way. In either case $dh'$ is a rational
multiple of the Tamagawa measure. The archimedean factors are not constrained.

Consequently $\prod_{v\in\Sinf}dh'_v$ is determined up to $\QQ^\times$, but the individual factors are not:
$dh'_{v_1}\mapsto t\,dh'_{v_1}$, $dh'_{v_2}\mapsto t^{-1}dh'_{v_2}$ respects every convention. For
$[F^+:\QQ]\geq2$ the per-place statement of \cite[Conjecture 4.15]{ghl} is therefore not invariant under the
admitted choices, while the product $\prod_{v\in\Sinf}I^*_v$, which is what the proof of
\cite[Theorem 6.1]{ghl} uses, is well posed up to $\QQ^\times$.

\subsection{The de~Rham structures at the CM point}\label{ghl:sec:DR}

We take de~Rham structures as in \cite[\S2.3]{harris13a}, defined by a model of the standard principal bundle as below, on $\pi(q)$ and $\pi'(q)$ at a common CM point, as in \cite[pp.~2045--2046]{harris13a} and \cite[\S3.1, p.~87]{harris14}. We fix the point as follows. Write $V'=V_1\oplus\dots\oplus V_{n-1}$ as an orthogonal sum of lines and put $V_n=V'_1$, the orthogonal complement of $V'$ in $V$. Let $x$ be the point attached to $V=V_1\oplus\dots\oplus V_n$ as in \cite[\S4.2, pp.~2035--2036]{harris13a}, $x'$ the point attached to $V'=V_1\oplus\dots\oplus V_{n-1}$, and
$$
T=\prod_i\U(V_i)\subset H'',\quad T'=T\cap H'=\prod_{i<n}\U(V_i).
$$
Exactly one line, $V_-$ say, is negative at $\vz$. It lies in $V'$, and every line is positive at $v\neq\vz$. Write $T_x\subset T$ and $T'_{x'}\subset T'$ for the Mumford--Tate groups of $x$ and $x'$.

Let $U\subset H(\bfA_f)$ and $U'\subset H'(\bfA_f)$ be neat open compact subgroups, $s_0$ the image of $x$ in $\mathrm{Sh}_U(H,Y_V)$, and $s_{T_x}$ that of $x$ in $\mathrm{Sh}(T_x,x)$ at the level $U\cap T_x(\bfA_f)$. Write likewise $s'_0$ and $s'_{T'_{x'}}$ for $x'$. We denote by $k(\cdot)$ residue fields. As $s_{T_x}$ maps to $s_0$, $k(s_0)\subset k(s_{T_x})$. Let $Z_H$ and $Z_{H'}$ be the centres of $H$ and $H'$, $\mu_x$ and $\mu_{x'}$ the cocharacters attached to $x$ and $x'$, $K_x$ and $K'_{x'}$ the centralizers of $\mu_x$ in $H$ and of $\mu_{x'}$ in $H'$, and $L_0\subset\bar\QQ$ a finite extension of $F^{\rm Gal}$ over which the irreducible representations of $H$,
$H'$, $K_x$ and $K'_{x'}$ are defined. The de~Rham structures of \cite[\S2.3]{harris13a} are defined by a model of the standard principal bundle over $\mathrm{Sh}_U(H,Y_V)$, with its $H(\bfA_f)$-action and its equivariant map to $\check Y_V$. By a {\em model} we mean such a model over a number field $E_1\supset E(H,Y_V)$, compatible with the analytic description of \cite[\S3.4]{harris85}. It defines models over $E_1L_0$ of the automorphic vector bundles $[\mathcal V_W]$, and the de~Rham structure on $W$ at $x$ is the fibre of the model of $[\mathcal V_W]$ at $s_0$, a vector space over $E_1L_0k(s_0)$, identified with $W(\CC)$ by the analytic trivialization. A model over $E(H,Y_V)$ is given by \cite[Theorem 4.8]{harris85}.
We fix a model for $H$ and one for $H'$, over $E_1$ and $E'_1$.

\begin{proposition}\label{ghl:prop:DR}
  Put $L=E_1L_0k(s_{T_x})$.
  \begin{itemize}
  \item[(i)] There are $t_0\in T_x(\CC)$ and $z\in Z_H(\CC)$ such that for every algebraic representation $W$ of $K_x$ over $L_0$, with Betti structure $W_B=W(L_0)$, the de~Rham structure on $W\otimes L$ at $x$ is $t_0z\cdot(W_B\otimes L)$. The element $t_0$ may be chosen independently of the model.
  \item[(ii)] There is a model over $E(H,Y_V)$ for which (i) holds with $z=1$.
  \end{itemize}
  The same holds for $H'$, with $x'$, $T'_{x'}$, $Z_{H'}$, $E'_1$, $L'=E'_1L_0k(s'_{T'_{x'}})$, $t'_0$ and $z'$ in place of $x$, $T_x$, $Z_H$, $E_1$, $L$, $t_0$ and $z$.
\end{proposition}

We prove Proposition \ref{ghl:prop:DR} using the canonical local systems of \cite{harris85} and the CM motive of
$(T_x,x)$.

\begin{proof}[Proof of Proposition \ref{ghl:prop:DR}]
We give the proof for $H$. The proof for $H'$ proceeds mutatis mutandis. We view $H$, $T$ and $T_x$
as groups over $\QQ$ by restriction of scalars. The group $K_x$ is defined over $F^{\rm Gal}$, and over $L_0$ all of them are
split, so the irreducible representations of $H$ and
$K_x$ over $L_0$ are absolutely irreducible, and $E(H,Y_V)\subset F^{\rm Gal}\subset L_0$. Here $T$ is a maximal
torus of $H$, the roots of $T$ in $H$ span the lattice $Q$ of section \ref{ghl:sec:Rrational},
and $X^*(T)/Q=X^*(Z_H)$ is
free, $X^*$ denoting the characters over $L_0$, over which $T$ and $Z_H$ split. The hypotheses (1.1.1)--(1.1.4) of \cite[pp.~155--156]{harris85} hold for $(H,Y_V)$ and for all unitary
groups of section \ref{ghl:sec:global}, which have signature $(m,1)$ at $\vz$: the weight is trivial, the centre is
anisotropic over $\RR$, and the derived group is not compact at $\vz$. For $n=2$, $H'$ is a torus, Theorem 4.2 of
\cite{harris85} holds for the CM pair $(H',x')$ \cite[Remark 4.5.5]{harris85}, and the argument below applies with
$K'_{x'}=T'_{x'}=H'$.

(a) {\em The canonical local systems.} We first take the model of \cite[Theorem 4.8]{harris85}. For an
irreducible representation $\rho$ of $H$ on $V_\rho$, over $L_0$, write $[\mathcal V_\rho]:=[\mathcal V_{V_\rho}]$ and let $V_\rho^\nabla\subset[\mathcal V_\rho]_\CC$
be the canonical local system in $L_0$-vector spaces \cite[Theorem 4.2, Remark 4.5.1, Theorem 4.6.5]{harris85}, and
$V_0\subset V_\rho(\CC)$ its fibre at $s_0$ in the analytic trivialization. By \cite[(4.2.3), Proposition
4.3]{harris85}, $V_0$ is an $H(\QQ)$-stable $L_0$-subspace. As $H(\QQ)$ is Zariski dense in $H$, Schur's lemma gives
$V_0=c_\rho V_\rho(L_0)$, with $c_\rho\in\CC^\times$ determined modulo $L_0^\times$. The canonical local systems are
compatible with tensor products and direct factors: by the uniqueness
\cite[Remark 4.5.1]{harris85} this follows from the corresponding compatibility of the CM motives, which holds by
their construction \cite[\S3.15]{harris85}. Let $\lambda$, $\mu$ be dominant weights. Then
$V_{\lambda+\mu}$ occurs in $V_\lambda\otimes V_\mu$ with multiplicity one, hence $c_{\lambda+\mu}=c_\lambda c_\mu$
modulo $L_0^\times$. If $\lambda$, $\mu$ are regular and $\alpha$ is a simple root, $\lambda+\mu-\alpha$ is dominant
and $V_{\lambda+\mu-\alpha}$ occurs with multiplicity one as well, hence $c_{\lambda+\mu-\alpha}=c_\lambda c_\mu$.
Therefore $\lambda\mapsto c_\lambda$ extends to a homomorphism $X^*(T)\to
\CC^\times/L_0^\times$ which is trivial on $Q$. As $X^*(Z_H)$ is free, the induced homomorphism on $X^*(Z_H)$ lifts to
a homomorphism $X^*(Z_H)\to\CC^\times$, i.\,e.\ to an element $z\in Z_H(\CC)$, and $V_0=z\,V_\rho(L_0)$ for every
irreducible $\rho$.

(b) {\em The torus $T_x$.} The image of $x$ in $\mathrm{Sh}(T_x,x)$ defines the CM pair $(T_x,x)\subset(H,Y_V)$, with
reflex field $E_x=E(T_x,x)\subset F^{\rm Gal}$. It satisfies (1.1.3), (1.1.4), $T_x(\RR)$ being compact. By
\cite[Theorem 4.2 (iii)]{harris85} there is an $L_0$-rational isomorphism of the restriction of $[\mathcal V_\rho]$
to $\mathrm{Sh}(T_x,x)$ with the de~Rham realization of the CM motive $\mathfrak M(\rho|_{T_x})$ which carries
$V_\rho^\nabla$ to its Betti realization (4.2.1), (4.2.2). The motive is defined through a surjection
$(T_0,h_0)\to(T_x,x)$ from a product of CM pairs admitting symplectic embeddings \cite[Proposition 2.4.2,
Lemmas 3.15.4, 3.15.5]{harris85}, and the Mumford--Tate group $T_{h_0}$ of $h_0$ acts on the Betti realization
through its image $T_x$. Let $\rho$ be faithful. The functor of \cite[\S3.15]{harris85} is compatible with the
decomposition of $\rho|_{T_x}$ into characters, so the $T_x$-invariant tensors in the spaces $V_0^{\otimes a}\otimes
V_0^{\vee\otimes b}$ are the elements of the Betti realization of the trivial part of the motive, i.\,e.\ Hodge
classes. They are invariant under $T_x(\bfA_f)$, since $T_x(\bfA_f)$ acts trivially on
the analytic trivialization. By the argument of the proof of \cite[Lemma 3.15.3, p.~177]{harris85}, which rests on
\cite[I, Theorem 2.11, Proposition 2.9]{dmos}, their de~Rham components are
$L_0$-rational sections over $\mathrm{Sh}(T_x,x)$, so their values at $s_{T_x}$ are rational over
$L=L_0k(s_{T_x})$. Since $T_x$ is reductive, it is the subgroup of $H$ which fixes these tensors
\cite[I, Proposition 3.1 (c)]{dmos}. Hence the $L$-linear isomorphisms from $V_0\otimes L$ to the de~Rham structure of
$[\mathcal V_\rho]$ at $s_0$ which preserve them form, after the comparison isomorphism, a torsor under $T_x$ over
$L$ (the push-out along $T_{h_0}\to T_x$ of the torsor of the CM motive). It is trivial by Hilbert's Theorem 90,
$T_x$ being split over $L_0$, and an $L$-point gives $t_0\in T_x(\CC)$ with de~Rham structure
$t_0(V_0\otimes L)=t_0z\,V_\rho(L)$. The same $t_0$ serves for every representation of $H$, all of which are direct
factors of tensor constructions on $\rho$, and the rational direct factors give rational subbundles
\cite[Lemma 3.5.1 (b)]{harris85}.

(c) {\em Representations of $K_x$.} The cocharacter $\mu_x$ is central in $K_x$, so
$V_\rho|_{K_x}=\bigoplus_p{\rm gr}^pV_\rho$ for the Hodge filtration $F^\bullet$, which is rational on
$\check Y_V$ \cite[3.1.3]{harris85} and on $[\mathcal V_\rho]$
\cite[Lemma 3.5.1 (a)]{harris85}. The functor $\mathcal V\mapsto[\mathcal V]$, defined by descent of
$\pi^*\mathcal V$ \cite[\S3.4]{harris85}, is exact, so the
de~Rham structure on $[\mathcal V_{{\rm gr}^pV_\rho}]$ at $s_0$ is ${\rm gr}^p$ of that on $[\mathcal V_\rho]$. On the Betti
side the analytic trivialization carries $F^\bullet$ to the filtration at the point defined by $x$. As $t_0z\in T(\CC)$
preserves $F^\bullet$, the de~Rham structure on ${\rm gr}^pV_\rho\otimes L$ is $t_0z\,{\rm gr}^pV_\rho(L)$. We may
assume $W$ irreducible. Every such representation $W$ of $K_x$ is a subquotient of the restriction of a representation
of $H$ \cite[Proposition 2.21 (b)]{delignemilne}, hence, $K_x$ being reductive, a direct summand of some ${\rm
gr}^pV_\rho$.
The unipotent radical of the stabilizer $P_x$ of the point acts trivially on ${\rm gr}^pV_\rho$, so an $L_0$-rational
$K_x$-equivariant projector onto $W$ is a morphism of homogeneous bundles on $\check Y_V$ over $L_0$, and defines an
$L$-rational projector on the fibre at $s_0$. Since $t_0z\in K_x(\CC)$ preserves $W$, the de~Rham structure on
$W\otimes L$ is $t_0z\,(W_B\otimes L)$.

(d) {\em Other models.} Let a second model be given, over $E_1$. For irreducible $\rho$ the two models of the
flat bundle $[\mathcal V_\rho]$ with its connection and $H(\bfA_f)$-action are isomorphic over $E_1L_0$
\cite[Lemma 3.5.1 (d)]{harris85}, and the automorphisms of the analytic triple are the elements of
$\mathrm{End}_H(V_\rho)(\CC)=\CC$, as in the proof of \cite[Lemma 3.5.1 (d)]{harris85}.
Hence the de~Rham structures of the second model on $V_\rho(\CC)$ at $s_0$ are $a_\rho$ times those of the first,
with $a_\rho\in\CC^\times$ determined modulo
$(E_1L_0)^\times$. Both models define tensor functors $\rho\mapsto[\mathcal V_\rho]$, and rational direct factors give
rational subbundles, so the argument of (a) gives $z_1\in Z_H(\CC)$ with $a_\rho\equiv\omega_\rho(z_1)$ modulo $(E_1L_0)^\times$ for all
irreducible $\rho$, where $\omega_\rho$ is the central character of $\rho$. Part (c) applies to the second model,
by \cite[Lemma 3.5.1 (a)]{harris85} and the exactness of the descent, and
gives (i) with $zz_1$ in place of $z$, over $E_1L_0k(s_{T_x})$. Finally, $z$ commutes with $H$ and with the action of
$H(\bfA_f)$ and acts trivially on $\check Y_V$, so transporting the model of \cite[Theorem 4.8]{harris85} by $z^{-1}$
gives a model over $E(H,Y_V)$. Its de~Rham structures are $z^{-1}$ times the former ones, i.\,e.\ $t_0\,(W_B\otimes L)$,
and its canonical local systems are the constant ones $V_\rho(L_0)$. This is (ii).
\end{proof}

For connected Shimura varieties and representations trivial on the centre the corresponding description, of the
conjugates of the automorphic vector bundles at special points through the comparison isomorphism of the CM motive,
is \cite[Remark 5.4 (b)]{milne88}.

We put $t=t_0z$ and $t'=t'_0z'$. Then $\beta(t)=1$ for every character $\beta$ of $T$ which is trivial on $T_xZ_H$, and likewise for $t'$. In Theorem \ref{ghl:thm:B} and in sections \ref{ghl:sec:Rrational} to \ref{ghl:sec:proofB}, $t$ and $t'$ enter only through such values and through values of characters
trivial on $Z_H$, resp.\ $Z_{H'}$, and the statements there are independent of the choice of the models. The classes of the Petersson norms of de~Rham rational forms in section \ref{ghl:sec:global} depend on the models, through $z$ and $z'$ (section \ref{ghl:sec:Fplus}), and there the models are fixed by \hypref{ghl:hyp:dR}{dR}. For the model of (ii), $t\in T_x(\CC)$. This is the description of the comparison isomorphisms in \cite[\S2.3, p.~2011]{harris13a}.

The fields $E(\pi(q))$, $E(\pi'(q))$ of \cite[\S3.2.2]{ghl} are not canonical, and they may be enlarged \cite[p.~26]{ghl}. We write $E$ for the field generated by $E(\pi(q))E(\pi'(q))$, the fields $L$ and $L'$ of Proposition \ref{ghl:prop:DR}, and their complex conjugates, and replace the de~Rham structures by their extensions of scalars to $E$, so that $\bar E=E$. The field $E$ contains the residue fields $k(s_0)$, $k(s'_0)$, and $E$ need not be CM (Remark \ref{ghl:rmk:nonCM}).

\begin{remark}\label{ghl:rmk:nonCM}
  For $F^+=\QQ$, $F=\QQ(\sqrt{-23})$ and $n=3$, the residue field of every point of $\mathrm{Sh}_U(H,Y_V)$ contains the Hilbert class field $H_F$ of $F$. Indeed, $\det$ defines a morphism from $(H,Y_V)$ to the Shimura datum of $C=\U(1)_{F/\QQ}$ with $h_C(u)=(u/\bar u)^k$, where $k=2$ or $k=-1$ according to the normalization of $h$, and by the reciprocity law \cite[II, Definition 5.1, Proposition 5.2]{milne90} ${\rm Gal}(F^{\rm ab}/F)$ acts on $C(\QQ)\backslash C(\bfA_f)/U_C$, $U_C\supset\det U$ the maximal compact subgroup of $C(\bfA_f)$, through $a\mapsto(a/\bar a)^{\pm k}$. The ideles in $F^\times F_\infty^\times\hat{\mathcal O}_F^\times$ act trivially, and if $(a/\bar a)^k\in C(\QQ)U_C$, then the class of the ideal of $a$ has order dividing $2k$, hence is trivial, the class number being $3$. Therefore ${\rm Gal}(H_F/F)$ acts freely on this finite set, which gives the assertion. The field $H_F$ is the splitting field of $X^3-X-1$, with Galois group $S_3$ over $\QQ$. It is not totally real, and it is not CM, complex conjugation not being central in $S_3$, hence it lies in no CM field.

  The enlargement by $k(s_0)$ is not a Brauer obstruction in the sense of \cite[\S\S5--6]{januszewskirtl}: these obstructions vanish here, since $K_x$ splits over $F$ and is a product of unitary groups \cite[Propositions 5.8 and 6.14]{januszewskirtl}. The field $k(s_0)$ is imposed by the point $s_0$, through the reciprocity law.

  In this case $E$, which contains $k(s_0)$, is neither CM nor totally real, \hypref{ghl:hyp:Gal}{Gal} of section \ref{ghl:sec:global} cannot hold, and Theorem \ref{ghl:thm:C} (ii) is vacuous. In particular the field $E_1L_0k(s_0)$, over which the fibres at $s_0$ are canonically vector spaces, need not lie in a CM field.
\end{remark}

\subsection{Rationality of the branching coefficients}\label{ghl:sec:Rrational}

We determine which characters in $Q$, resp.\ $Q'$, are trivial on $T_x$, resp.\ $T'_{x'}$ (Lemma \ref{ghl:lem:Tx}), combine this with
Proposition \ref{ghl:prop:DR} to compare de~Rham and Betti structures on $\U(V_-)(F^+_v)$-eigenspaces
(Lemma \ref{ghl:lem:DRv}) and deduce $R^\circ_v\in E$ on eigenvectors (Lemma \ref{ghl:lem:Rrational}).

Let $x$, $x'$, $T$, $T'$, $V_-$, $T_x$, $T'_{x'}$, $t$ and $t'$ be as in section \ref{ghl:sec:DR}. The torus $T$
splits over $F^{\rm Gal}\subset E$ \cite[p.~26]{ghl}. For $v\in\Sinf$ fix an embedding $\iota_v$ of $F$ above $v$
and let $e_{i,v}$ be the character $\iota_v$ of the factor $\U(V_i)(F^+_v)$ of $T$. The $e_{i,v}$ form a basis of
the characters of $T$, and the $e_{i,v}$ with $i<n$ one of the characters of $T'$. Let $Q$ be the lattice of
characters $\sum m_{i,v}e_{i,v}$ of $T$ with $\sum_im_{i,v}=0$ for every $v$, and $Q'$ the corresponding
lattice for $T'$. They contain the roots of $T$ in $H$, resp.\ of $T'$ in $H'$.

\begin{lemma}\label{ghl:lem:Tx}
  A character $\beta=\sum m_{i,v}e_{i,v}\in Q$ is trivial on $T_x$ if and only if $m_{-,v}=0$ for every
  $v\in\Sinf$. The same holds for $\beta\in Q'$ and $T'_{x'}$.
\end{lemma}

\begin{proof}
  The character $\beta$ is trivial on $T_x$ if and only if $\langle\sigma\beta,\mu_x\rangle=0$ for all
  $\sigma\in{\rm Gal}(\bar\QQ/\QQ)$. Now $\sigma e_{i,v}=\pm e_{i,v'}$, where $v'$ is the place of
  $\sigma\circ\iota_v$,
  and $\langle e_{i,w},\mu_x\rangle-\langle e_{j,w},\mu_x\rangle$ is $0$ if $V_i$, $V_j$ have the same sign at
  $w$ and non-zero otherwise. Fix $j\neq-$. Then $V_j$ is positive everywhere. Since
  $\sum_im_{i,v}=0$, $\langle\sigma\beta,\mu_x\rangle$ is a non-zero multiple of $m_{-,v}$ for the place $v$ with
  $v'=\vz$, and every $v$ occurs. For $T'$ the argument is the same, $V'$ being a sum of $n-1\geq2$ lines.
\end{proof}

In particular the compact roots of $T$ at $\vz$ are trivial on $T_x$, and at $v\neq\vz$ the roots
$e_{i,v}-e_{-,v}$ are not. Every $e_{i,v}$ restricts to $\iota_v$ on $Z_H$, so the elements of $Q$ are trivial on
$Z_H$, and those of $Q'$ on $Z_{H'}$. By Proposition \ref{ghl:prop:DR}, $\beta(t)=1$ for $\beta\in Q$ trivial on
$T_x$, and $\beta(t')=1$ for $\beta\in Q'$ trivial on $T'_{x'}$. For a representation $W$ of $H_v$ or $H'_v$ and $m\in\ZZ$ let $W[m]$ be the subspace on
which $\U(V_-)(F^+_v)$ acts by $\iota_v^m$.

\begin{lemma}\label{ghl:lem:DRv}
  Let $v\in\Sinf$, and let $W$ be the minimal $K$-type of $\pi_v$ or the minimal $K'$-type of
  $\pi'_v$, $W_E$ its de~Rham structure and $W_B$ its Betti structure. Then
  $$
    W_E=\bigoplus_m(W_E\cap W[m])
  $$
  and
  $$
    W_E\cap W[m]=c_m\cdot(W_B\cap W[m])
  $$
  with $c_m\in\CC^\times$. At $\vz$, $W=W[m]$ for a single $m$.
\end{lemma}

\begin{proof}
  Both structures are stable under $T_E$, resp.\ $T'_E$. By Proposition \ref{ghl:prop:DR},
  $$
    W_E=\bigoplus_\beta\beta(t)\,W_{B,\beta},
  $$
  where the sum is over the weights of $W$, with $t'$ in place of $t$
  for $\pi'_v$. Two weights of $W$ with the same
  eigenvalue of $\U(V_-)(F^+_v)$ differ by an element of $Q$, resp.\ $Q'$, with $m_{-,w}=0$ for all $w$, which is
  trivial on $T_x$, resp.\ $T'_{x'}$, by Lemma \ref{ghl:lem:Tx}. Its value at $t$,
  resp.\ $t'$, is $1$. At $\vz$, $\U(V_-)(F^+_{\vz})$ is a central factor of $K_{\vz}$ and of
  $K'_{\vz}$.
\end{proof}

\begin{lemma}\label{ghl:lem:Rrational}
  Let $v\in\Sinf$, $f_v$ and $f'_v$ rational elements of the minimal $K$-types of $\pi_v$ and
  $\pi'_v$ for the de~Rham structures which are eigenvectors of $\U(V_-)(F^+_v)$, and $R^\circ_v(f_v,f'_v)$ as in
  section \ref{ghl:sec:R} (at $v\neq\vz$ with $K_v=H_v$, $\tau=\pi_v$). Then $R^\circ_v(f_v,f'_v)\in E$, and
  $R^\circ_v\not\equiv0$ on such vectors.
\end{lemma}

\begin{proof}
  By Lemma \ref{ghl:lem:DRv}, $f_v$ and $f'_v$ are multiples of Betti rational vectors, and $R^\circ_v$ depends
  only on the lines $\CC f_v$ and $\CC\bar f'_v$. We may therefore work with the Betti structures. They are
  stable under $\liek_E={\rm Lie}(K_x)_E$, resp.\ $\liek'_E=\liek_E\cap\liek'_{v,\CC}$, $K_x$ and $K'_{x'}=K_x\cap H'$
  being the centralizers of $\mu_x$ in $H$ and of $\mu_{x'}$ in $H'$. The conjugate of $\mu_x$ is $\bar\mu_x$, and
  $\mu_x\bar\mu_x$ is central, so $K_x$ is stable under complex conjugation, and so is $\liek_E$, as $\bar E=E$.
  The same holds for $x'$. The action of $X\in\liek'_{v,\CC}$ on $\bar\pi'_v$ is
  $\bar w\mapsto\overline{c(X)w}$, where $c$ is the
  conjugation with respect to $\liek'_v$. Hence the conjugate of the $E$-structure on the minimal $K'$-type of
  $\pi'_v$ is a $\liek'_E$-stable $E$-structure on $\tau_\sigma$, and $u=\bar f'_v$ is rational for it.
  The space $\Hom_{K'}(\tau_\sigma,\tau)$ is a line, and its $\liek'_E$-invariant rational elements
  form an $E$-line,
  spanned by $\varphi_0$ say. Then $P_\tau\circ j|_{\tau_\sigma}=\kappa\varphi_0$, $\kappa\in\CC^\times$, and
  $w=\kappa\varphi_0(u)$ with $\varphi_0(u)\in\tau_E$. The same argument, applied to $\tau\otimes\bar\tau$, gives a
  non-zero invariant $E$-valued form on $\tau_E\times\tau_E$, hermitian since $\bar E=E$. By Schur's lemma the
  inner product of $\pi_v$ on $\tau$ is a positive multiple $H_\tau$ of it. Since $R^\circ_v$ is homogeneous of
  degree $0$ in $w$ and in the form,
  $$
    R^\circ_v(f_v,f'_v)=\frac{H_\tau(f_v,\varphi_0u)\,\overline{H_\tau(f_v,\varphi_0u)}}{H_\tau(f_v,f_v)\,
    H_\tau(\varphi_0u,\varphi_0u)}\in E,
  $$
  as $\bar E=E$. The vector $\varphi_0u$ is an eigenvector of $\U(V_-)(F^+_v)\subset K'_v$, and
  taking $f_v=\varphi_0u$ gives $R^\circ_v=1$.
\end{proof}

\begin{remark}\label{ghl:rmk:definiteperiods}\leavevmode\par
  \noindent {\rm(1)} At $\vz$ every de~Rham rational vector of the minimal types is admissible in Lemma \ref{ghl:lem:Rrational},
  and the de~Rham structures there are multiples of the Betti structures.

  \smallskip
  \noindent {\rm(2)} Let $v\neq\vz$, $\gamma=e_{-,v}-e_{j,v}$ with $j\neq-$, and let $f_v=f_1+f_2$ with
  $f_k\in\tau_E\cap\tau[m_k]$ non-zero, $m_1\neq m_2$, $f_1$ and $f'_v$ as in Lemma \ref{ghl:lem:Rrational} with
  $R^\circ_v(f_1,f'_v)\neq0$. Then $f_2$ is orthogonal to $f_1$ and to $w$, and
  $$
    R^\circ_v(f_v,f'_v)=R^\circ_v(f_1,f'_v)\,\frac{\|f_1\|^2}{\|f_1\|^2+\|f_2\|^2},\qquad
    \frac{\|f_2\|^2}{\|f_1\|^2}\in|\gamma(t)|^{2(m_2-m_1)}\,E^\times .
  $$
  Hence $R^\circ_v(f_v,f'_v)\in E$ if and only if $|\gamma(t)|^{2(m_2-m_1)}\in E$.
  Complex conjugation maps $\gamma(t)$
  into $\bar\QQ^\times\gamma(t)$, i.\,e.\ $\gamma(\bar t\,t)\in\bar\QQ^\times$, since
  $\overline{\gamma(t)}=\gamma(\bar t)^{-1}$, the conjugate of the character $\iota_v$ of $\U(V_i)(F^+_v)$ being
  $\iota_v^{-1}$. Indeed, the conjugate CM datum is $(T,h_x^{-1})$, $h_x$ denoting the homomorphism $x$ with values in
  $T_\RR$, which the inversion of $T$ identifies with $(T,h_x)$ \cite[Proposition 2.1, p.~18]{ghl}. The condition
  therefore implies
  $\gamma(t)^{2(m_2-m_1)}\in\bar\QQ^\times$.
  The character $\gamma$ is not trivial on $T_x$ (Lemma \ref{ghl:lem:Tx}), and one expects
  $\gamma(t)$ to be transcendental of
  infinite order modulo $\bar\QQ^\times$. We do not prove this. In this sense \cite[Conjecture 4.15]{ghl} for all
  rational $f_v$ is not expected at the definite places when $F^+\neq\QQ$.

  \smallskip
  \noindent {\rm(3)} The same computation, for a non-compact root $a$ at $\vz$, shows that, with the description of
  $U(\lieg)_{{\rm DR},x}$ in \cite[p.~2011]{harris13a}, the statement of \cite[Proposition 3.19]{harris13a}
  cannot hold for all $v,v_1,v_2\in\pi_{{\rm DR},x}$. There $U(\lieg)_{{\rm DR},x}$ is the $E_x$-algebra generated
  by ${\rm Lie}(T)$ and the elements $p(a)X_a$, $X_a$ the standard root vectors for $T$ and $p(a)$ the CM period of the
  character $a$ of $T$. Let $G=\mathrm{GU}(1,1)$ be the unitary similitude group of the hermitian form $\diag(1,-1)$
  on $\QQ(i)^2$, $x$ the CM point attached to the two coordinate lines, and $T$ the diagonal torus, so that the
  $K_x$-types are lines. We identify the derived algebra of $\lieg_\CC$ with $\mathfrak{sl}_2(\CC)$ through the first
  embedding of $\QQ(i)$, and take $X_a=E_{12}$, $X_{-a}=E_{21}$, $H_a=[X_a,X_{-a}]=\diag(1,-1)$; other
  choices of the root vectors change them by factors in $\bar\QQ^\times$, which do not affect the argument. The
  conjugation $c$ of $\lieg_\CC$ with respect to ${\rm Lie}\,G(\RR)$ is $X\mapsto-J\,{}^t\bar X J$, $J=\diag(1,-1)$,
  on $\mathfrak{sl}_2(\CC)$, so $c(X_a)=X_{-a}$, and on a unitary representation the adjoint of $X$ is $-c(X)$. Let $\pi$ be as
  in \cite[\S3.4.2]{harris13a}, with holomorphic discrete series component at infinity, and let $w\in\pi_{{\rm DR},x}$ be
  non-zero in the minimal $K_x$-type. Interchanging the two coordinates if necessary, we have $X_{-a}w=0$. Let
  $w_1=p(a)X_aw\in\pi_{{\rm DR},x}$.
  Then $H_aw=kw$ with $k\in\ZZ$, $X_{-a}X_aw=[X_{-a},X_a]w=-H_aw$, and
  $$
  \begin{aligned}
    \langle w_1,w_1\rangle&=p(a)\overline{p(a)}\,\langle X_aw,X_aw\rangle
      =-|p(a)|^2\langle w,c(X_a)X_aw\rangle\\
    &=-|p(a)|^2\langle w,X_{-a}X_aw\rangle=|p(a)|^2\langle w,H_aw\rangle=k\,|p(a)|^2\langle w,w\rangle .
  \end{aligned}
  $$
  Here $k>0$, since $k\langle w,w\rangle=\|X_aw\|^2$ and $X_aw\neq0$, as otherwise $\CC w$ would be stable under
  $\lieg_\CC$. Hence $Q(w)^{-1}\langle w_1,w_1\rangle=k|p(a)|^2$, and \cite[Proposition 3.19]{harris13a} would give
  $|p(a)|^2\in\bar\QQ$. Equivalently, since $p(a)p(-a)\in\bar\QQ^\times$,
  $$
    c(p(a)X_a)=\overline{p(a)}\,X_{-a}\in|p(a)|^2\,\bar\QQ^\times\cdot p(-a)X_{-a},
  $$
  so $U(\lieg)_{{\rm DR},x}\otimes\bar\QQ$ is stable under $c$ only if $|p(a)|^2\in\bar\QQ^\times$.
  The character $a\circ h_x$ is $z\mapsto(z/\bar z)^{\pm1}$. The CM periods are multiplicative modulo
  $\bar\QQ^\times$, and, up to $\bar\QQ^\times$ and to the normalization, the period of the character $z$ of
  $\Res_{\QQ(i)/\QQ}\mathbb G_m$ is $\Omega=\Gamma(\frac14)^2/(2\sqrt{2\pi})$, which is real, and that of $z\bar z$ is a
  power of $2\pi i$. Hence $p(a)\in\bar\QQ^\times\,\Omega^{2\epsilon}(2\pi i)^j$ with $\epsilon=\pm1$, $j\in\ZZ$, and
  $$
    |p(a)|^2\in\bar\QQ^\times\,\Gamma(\tfrac14)^{8\epsilon}\pi^{2j-2\epsilon},
  $$
  which is transcendental, $\pi$ and $\Gamma(\frac14)$ being algebraically independent
  \cite{chudnovsky}. The exponents $\epsilon$ and $j$ depend on the normalization of $p(a)$. By
  \cite[Remark 3.21]{harris13a} the proposition is meant for all vectors.
\end{remark}

\subsection{The definite places}\label{ghl:sec:definite}

\begin{lemma}\label{ghl:lem:definite}
  For $v\neq\vz$ and all $f_v\in\pi_v$, $f'_v\in\pi'_v$,
  $$
    I^*_v(f_v,f'_v)=\vol_{dh'_v}(H'_v)\,(\dim F_{\Lambda'_v})^{-1}R^\circ_v(f_v,f'_v).
  $$
\end{lemma}

\begin{proof}
  As in section \ref{ghl:sec:integral},
  $$
    I_v(f_v,f'_v)=\int_{H'_v}\langle\pi_v(h)f_v,f_v\rangle\overline{\langle\sigma_v(h)u,u\rangle}dh'_v
  $$
  with $u=\bar f'_v$. The restriction $\pi_v|_{H'_v}$ is multiplicity
  free and contains $\sigma_v$ by interlacing. Let $j$ be the isometric embedding, unique up to a phase.
  By Schur orthogonality for the compact group $H'_v$ only the $\sigma_v$-component $jj^*f_v$ contributes, and
  $$
    I_v=\vol(H'_v)(\dim\sigma_v)^{-1}|\langle f_v,ju\rangle|^2,
  $$
  with $\dim\sigma_v=\dim F_{\Lambda'_v}$.
\end{proof}

\subsection{The product of the archimedean volumes}\label{ghl:sec:lemmaM}

Using a totally definite inner form of $H'$, we prove in Lemma \ref{ghl:lem:M} that the product of the
volumes from Theorem \ref{ghl:thm:A} and Lemma \ref{ghl:lem:definite} is rational for the measures of
section \ref{ghl:sec:globalmeasures}.

For $v\in\Sinf$ let $\volc_{dh'_v}(\U(n-1))$ be the volume of the compact dual for the measure matched
with $dh'_v$ (section \ref{ghl:sec:measures}). At the definite places it is $\vol_{dh'_v}(H'_v)$. We use three
standard facts.
\begin{itemize}
  \item[(T1)] There is a Hermitian space $V''$ over $F/F^+$ of dimension $n-1$, positive definite at every
    $v\in\Sinf$ and isomorphic to $V'$ at all finite places outside a finite set $S$ (change the local invariant
    at one finite place inert in $F$ to compensate the change of signature at $\vz$). We take $S$ to contain
    the finite places at which $dh'_p$ does not give volume $1$ to hyperspecial maximal compact subgroups.
  \item[(T2)] Tamagawa measures of inner forms are built from top forms $\omega'$, $\omega''=\psi_*\omega'$ ($\psi$ an
    inner twist; $\psi_*\omega'$ is $F^+$-rational since $\det\mathrm{Ad}=1$ on inner automorphisms) with the same
    convergence factors and the same global constants, so $dh'_v=a_v|\omega'|_v$ with $a_v$ the same for the
    inner form.
  \item[(T3)] The Tamagawa number of $\U(V'')$ is rational (it is $2$).
\end{itemize}

\begin{lemma}\label{ghl:lem:M}
  Let $dh'=\prod_vdh'_v$ be a rational multiple of the Tamagawa measure on $H'(\bfA_{F^+})$ whose local factors at
  the finite places take rational values on compact open subgroups and give volume $1$ to hyperspecial maximal
  compact subgroups at almost all places, for instance the measure of section \ref{ghl:sec:globalmeasures}. Then
  $$
    \mathrm{Vol}_\infty:=\prod_{v\in\Sinf}\volc_{dh'_v}(\U(n-1))\in\QQ^\times.
  $$
\end{lemma}

\begin{proof}
  Put $dh''_v:=a_v|\omega''|_v$ at every place. At a finite place outside $S$ an isomorphism
  $\U(V')_p\to\U(V'')_p$ differs from $\psi$ by an automorphism, which acts on $\omega$ by $\pm1$, so $dh''_p$
  gives volume $1$ to hyperspecial subgroups. At $p\in S$, $dh'_p/|\omega'|_p\in\QQ^\times$, since $dh'_p$ is
  rational on compact open subgroups and $|\omega'|_p(U)\in|a|_p\,q^\ZZ$ for small congruence subgroups $U$,
  so $dh''_p$ is rational. By (T2), $\prod_vdh''_v=r\cdot\mu^{\rm Tam}_{V''}$ with the same $r\in\QQ^\times$ as
  for $V'$. At $v\in\Sinf$ the compact real forms of $\U(V')_v\otimes\CC$ are conjugate, $\det\mathrm{Ad}=1$ and
  $|i^{\dim\liep}|=1$, so $\volc_{dh'_v}(\U(n-1))=\vol_{dh''_v}(\U(V'')_v)$. Since $\U(V'')(F^+_\infty)$ is
  compact, $\U(V'')(F^+)$ is discrete in $\U(V'')(\bfA_{F^+,f})$, and for an open compact $K_f$
  \begin{align*}
    r\,\tau(\U(V''))&=\vol\big(\U(V'')(F^+)\backslash\U(V'')(\bfA_{F^+})\big)\\
    &=\mathrm{Vol}_\infty\cdot\sum_{i=1}^h\frac{dh''_f(K_f)}{|\Gamma_i|}
  \end{align*}
  with finitely many classes, finite groups $\Gamma_i$, and $dh''_f(K_f)\in\QQ^\times$. For the measure of section
  \ref{ghl:sec:globalmeasures} the convention
  $$
    \vol\big(\U(V_1)(F^+)\backslash\U(V_1)(\bfA_{F^+})\big)=1
  $$
  changes $r$ by a rational factor.
\end{proof}

\begin{remark}\label{ghl:rmk:anymeasure}
  Assume (T1)--(T3). Let $dh'=dh'_\infty\,dh'_f$ be any Haar measure on $H'(\bfA_{F^+})$,
  $dh'=c\cdot\mu^{\rm Tam}$, and $K_f\subset H'(\bfA_{F^+,f})$ open compact. The proof of Lemma
  \ref{ghl:lem:M}, in which $r\in\QQ^\times$ is arbitrary, applies to a factorization
  $\mu^{\rm Tam}=\prod_v\mu^0_v$ satisfying the conditions of Lemma \ref{ghl:lem:M}. Comparing $dh'$ with it we get
  $\mathrm{Vol}_\infty\cdot dh'_f(K_f)\in c\cdot\QQ^\times$. All terms of Theorem \ref{ghl:thm:A} and Lemma
  \ref{ghl:lem:definite} other than $I^*_v$, $d(\sigma)$, $\volc$ and $\vol(H'_v)$ are independent of the
  measure, and these are homogeneous of degree $1$, $-1$, $1$, $1$ in $dh'_v$. Hence the identity of Theorem
  \ref{ghl:thm:B} holds for every $dh'_\infty$, and for $f$, $f'$ as there,
  $$
    \prod_{v\in\Sinf}I^*_v\in c\,dh'_f(K_f)^{-1}\cdot E.
  $$
\end{remark}

For totally definite $V'$ the lemma is \cite[Lemma 4.1.5, Corollary 4.1.6, p.~145]{harris13b}. Combined with \eqref{ghl:eq:fdeg}, Lemma \ref{ghl:lem:M} gives
$$
d(\sigma)^{-1}\prod_{v\neq\vz}\vol(H'_v)=\frac{\mathrm{Vol}_\infty}{\dim F_{A'(q)}}\in\QQ^\times
$$
for every discrete series $\sigma$ of $H'_{\vz}$. No $L$-value enters.

\subsection{Proof of Theorem \ref{ghl:thm:B}}\label{ghl:sec:proofB}

In the setting of \cite[Conjecture 4.15]{ghl} all gaps $A_i-A_{i+1}$ are at least $2$, so
\eqref{ghl:eq:gap} holds, and \hypref{ghl:hyp:J}{J} holds by Theorem \ref{ghl:thm:J} $(q\geq1)$ and Proposition
\ref{ghl:prop:Jq0} $(q=0)$. By Theorem \ref{ghl:thm:A} at $\vz$ and Lemma \ref{ghl:lem:definite} at $v\neq\vz$,
$$
  \prod_{v\in\Sinf}I^*_v(f_v,f'_v)=\frac{\volc_{dh'_{\vz}}}{\dim F_{A'(q)}\,s(q)}\prod_{v\neq\vz}
  \frac{\vol(H'_v)}{\dim F_{\Lambda'_v}}\prod_{v\in\Sinf}R^\circ_v
  =\frac{\mathrm{Vol}_\infty}{s(q)}\prod_{v\in\Sinf}\frac{1}{\dim F_{\Lambda'_v}}\prod_{v\in\Sinf}R^\circ_v.
$$
By Lemma \ref{ghl:lem:M}, \eqref{ghl:eq:P} and Lemma \ref{ghl:lem:Rrational} the right side lies in
$\QQ^\times\cdot E$. \qed

If $f$, $f'$ are as in \cite[Theorem 4.14]{ghl}, with $I^{\rm can}(f,f')\neq0$ and rational archimedean
components in the minimal types, decompose $f_v$ and $f'_v$ for $v\neq\vz$ according to the eigenvalues of
$\U(V_-)(F^+_v)$ (Lemma \ref{ghl:lem:DRv}). Since $I^{\rm can}$ is multilinear in the archimedean components,
one of the resulting factorizable pairs, all of whose components are rational eigenvectors, has $I^{\rm can}\neq0$.

\begin{corollary}\label{ghl:cor:perplace}
  Under the assumptions of Theorem \ref{ghl:thm:B}, for every splitting $dh'=\prod_vdh'_v$ as in Lemma \ref{ghl:lem:M} with
  $\vol_{dh'_v}(H'_v)\in\QQ^\times$ at the definite places, $I^*_v(f_v,f'_v)\in E$ for every $v\in\Sinf$. At
  $\vz$,
  $$
  I^*_{\vz}=C\cdot R^\circ_{\vz}\quad\text{with}\quad C=\frac{\mathrm{Vol}_\infty}{\prod_{v\neq\vz}\vol(H'_v)\cdot\dim F_{A'(q)}\cdot s(q)}\in\QQ^\times.
  $$
  For $F^+=\QQ$, \cite[Conjecture 4.15]{ghl} holds as stated, for the field $E$.
\end{corollary}

\begin{proof}
  The first two assertions follow from Lemmas \ref{ghl:lem:definite}, \ref{ghl:lem:M}, \ref{ghl:lem:Rrational} and
  Theorem \ref{ghl:thm:A}. For $F^+=\QQ$ there is no definite place and $dh'_{\vz}$ is fixed up to $\QQ^\times$.
\end{proof}

\begin{remark}\label{ghl:rmk:n2}
  For $n=2$, $H'_{\vz}=\U(0,1)$ is compact and equal to the non-compact $\U(1)$-factor of $K$. The minimal
  $K$-types are lines and $I^*_{\vz}=\vol(H'_{\vz})$. The space $V'$ is totally definite, so
  $\prod_{v\in\Sinf}\vol(H'_v)\in\QQ^\times$ by \cite[Corollary 4.1.6, p.~145]{harris13b}, and with Lemma
  \ref{ghl:lem:definite} the product form holds for $n=2$ with the measures of section \ref{ghl:sec:globalmeasures}, for $f_v$ eigenvectors of
  $\U(V_-)(F^+_v)$ at $v\neq\vz$.
\end{remark}

\begin{remark}\label{ghl:rmk:scope}
  The orientation $\pi'(q)_{\vz}=\pi^\vee_{A'(q)}$ is essential: with $\pi_{A'(q)}$ in place of its dual the
  integral vanishes, already for $n=2$ (if $A_1\neq-\frac12$).
\end{remark}
 \section{Global consequences}\label{ghl:sec:global}

We determine what Theorems \ref{ghl:thm:A}, \ref{ghl:thm:J} and \ref{ghl:thm:B} give for
\cite[Theorem 6.1]{ghl} when $F^+=\QQ$, and, under a further condition, when $[F^+:\QQ]\geq2$ (section
\ref{ghl:sec:Fplus}). All number fields are subfields of $\CC$ \cite[\S1.1]{ghl}, and for a
number field $L$ we write $a\in b\cdot L^\times$ for the relation at the given embedding. The relations $\sim$ of
\cite[Definition 1.5]{ghl} imply it \cite[Definition 1.5 (ii)(1), p.~13]{ghl}.

Section \ref{ghl:sec:globalinputs} reviews \cite[\S6.2]{ghl}, proves Lemma \ref{ghl:lem:pattern} and states the
hypotheses. Sections \ref{ghl:sec:thmC} and \ref{ghl:sec:Fplus} prove Theorems \ref{ghl:thm:C} and
\ref{ghl:thm:Cprime}.
Theorem \ref{ghl:thm:C} is the precise form of Theorem \ref{ghl:thm:Cintro}.

\subsection{The inputs of the proof of \cite[Theorem 6.1]{ghl}}\label{ghl:sec:globalinputs}

The proof \cite[\S6.2, pp.~47--50]{ghl} is an induction on $n$. For $q=0$ it is \cite[Observation 6.6]{ghl},
which does not use Conjecture 4.15 of loc.\ cit. For $q\geq1$ the archimedean
input \cite[Conjecture 4.15]{ghl} is used only in the passage \cite[(6.11)$\to$(6.13)]{ghl}, and the induction
applies \cite[Theorem 6.1]{ghl} to $\Pi'$ and $q'=n-q-2$. Besides \cite[Conjecture 4.15]{ghl} the proof uses
\cite[Conjecture 2.10]{ghl}, through \cite[Proposition 5.12]{ghl} in the proof of \cite[Theorem 5.19]{ghl}, at
every step of the induction with $q\geq1$, for every $F^+$.

At every step of the induction the measures are those of \cite[p.~28]{ghl} for the pair of that step, normalized as
in section \ref{ghl:sec:globalmeasures}. Hence the measure on the smaller group which enters \cite[(6.11)]{ghl}
satisfies the conditions of Lemma \ref{ghl:lem:M} at every step, which is what is needed to insert Theorem
\ref{ghl:thm:B} into \cite[(6.11)]{ghl}. The identity of Theorem \ref{ghl:thm:B} holds for every measure (Remark
\ref{ghl:rmk:anymeasure}). The two Petersson norms on $\pi'(q)$ which enter the induction, for $H'$ as the smaller
group at rank $n$ and as the larger group at rank $n-1$, agree up to $\QQ^\times$, both measures being rational
multiples of Tamagawa measure.

The relation \cite[(6.11)]{ghl} rests on the automorphic periods of \cite[\S4.3, p.~33]{ghl}: $Q(\pi)$ is the
Petersson norm of a de~Rham rational form, independent of the choice up to $E(\pi)^\times$ by
\cite[Proposition 3.19, Lemma 3.20, p.~2033]{harris13a}. We do not use these statements. We use $Q(\pi)$ only for forms whose archimedean component lies in the
minimal type. These are the forms which enter \cite[(6.11)]{ghl}, \cite[Observation 6.6]{ghl} and the induction.
For them the dependence on the choice is determined in Lemma \ref{ghl:lem:pattern} below: for $F^+=\QQ$ there is
none, and for $[F^+:\QQ]\geq2$ it is a product of periods of characters of $T$, which we expect to be
transcendental. We treat $F^+=\QQ$ first (Theorem \ref{ghl:thm:C}) and the general case in section
\ref{ghl:sec:Fplus}.

Let $\pi$ be $\pi(q)$ or $\pi'(q)$, a representation of $H$, resp.\ $H'$, and $W$ the tensor product of the
minimal $K$-type of $\pi_{\vz}$ and of the $\pi_v$, $v\neq\vz$, resp.\ the same for $K'$. Then
$\pi_\infty\otimes\pi_f$ contains $W\otimes\pi_f$. For $m=(m_v)_{v\neq\vz}\in\ZZ^{\Sinf\smallsetminus\{\vz\}}$ let $W[m]$ be the subspace of $W$
on which $\U(V_-)(F^+_v)$ acts by $\iota_v^{m_v}$ for every $v\neq\vz$. A form $\varphi\in\pi$ is of {\em pattern} $m$
if it lies in $W[m]\otimes\pi_f$. The patterns of $\pi$ are the $m$ with $W[m]\neq0$. For $F^+=\QQ$ there is one
pattern, and the condition is that the archimedean component lies in the minimal $K$-type. Fix a line $V_j\subset V$
with $j\neq-$, and put $\gamma_v=e_{-,v}-e_{j,v}$ for $v\in\Sinf$. For $j'\neq-$, $e_{j,v}-e_{j',v}$ lies in $Q$, resp.\
$Q'$, and has all $e_{-,w}$-coordinates $0$. By Lemma \ref{ghl:lem:Tx} and Proposition \ref{ghl:prop:DR} its value at
$t$, resp.\ $t'$, is $1$, so $\gamma_v(t)$, resp.\ $\gamma_v(t')$, does not depend on $j$.

\begin{lemma}\label{ghl:lem:pattern}
  Let $\pi$ be $\pi(q)$, resp.\ $\pi'(q)$, and let $\varphi$, $\tilde\varphi\in\pi$ be non-zero
  de~Rham rational forms of patterns $m$ and $\tilde m$. Then
  $$
    \langle\varphi,\varphi\rangle\in\prod_{v\neq\vz}|\gamma_v(t)|^{2(m_v-\tilde m_v)}\cdot
    \langle\tilde\varphi,\tilde\varphi\rangle\cdot E^\times,
  $$
  with $t'$ in place of $t$ for $\pi'(q)$. In particular, for $F^+=\QQ$ the class of $\langle\varphi,\varphi\rangle$
  modulo $E^\times$ is the same for all non-zero de~Rham rational forms $\varphi\in\pi$ whose archimedean component
  lies in the minimal $K$-type.
\end{lemma}

\begin{proof}
  We write the proof for $\pi(q)$. The space $W\otimes\pi_f$ is irreducible under $K_\infty\times H(\bfA_{F^+,f})$,
  and the de~Rham structure of $\pi$ is stable under $H(\bfA_{F^+,f})$, so its part in $W\otimes\pi_f$ is
  $W'\otimes_E\pi_{f,E}$ for an $E$-structure $W'$ of $W$, and $W'=c\cdot W_E$, $c\in\CC^\times$, with $W_E$ the
  de~Rham structure of $W$ at $x$. Here $\pi_{f,E}$ is an $H(\bfA_{F^+,f})$-stable $E$-structure on $\pi_f$. By
  Schur's lemma the Petersson inner product on
  $W\otimes\pi_f$ is $h_\infty\otimes h_f$ for invariant hermitian forms $h_\infty$ on $W$ and $h_f$ on $\pi_f$.

  At the finite places, let $K_f$ be open compact with $\pi_f^{K_f}\neq0$, $e_1,\dots,e_d$ an $E$-basis of
  $\pi_{f,E}^{K_f}$, and $\mathcal H_\QQ$ the algebra of $\QQ$-valued bi-$K_f$-invariant functions with compact
  support, for the Haar measure with $\vol(K_f)=1$. The algebra $\mathcal H_\QQ$ acts on $\pi_{f,E}^{K_f}$ by
  matrices $A(\phi)$ with entries in $E$, it is stable under $\phi\mapsto\phi^*$, $\phi^*(g)=\phi(g^{-1})$, and
  the invariance of $h_f$ means
  $$
    A(\phi)^tG=G\,\overline{A(\phi^*)}
  $$
  for its Gram matrix $G$. As $\bar E=E$, these are linear equations for
  $G$ with coefficients in $E$. The space $\pi_f^{K_f}$ is an irreducible $\mathcal H$-module, so their solutions
  form a line, which therefore contains a matrix with entries in $E$.
  Hence $h_f$ is a multiple of a form with values in $E$ on
  $\pi_{f,E}^{K_f}$. For $K'_f\subset K_f$ such a form for $K'_f$ restricts to one for $K_f$, so the multiples are
  compatible modulo $E^\times$.

  At the archimedean places, by Proposition \ref{ghl:prop:DR},
  $$
  W_E=\bigoplus_\beta\beta(t)W_{B,\beta},
  $$
  the sum ranging over the weights $\beta$ of $W$. Two weights of $W[m]$ differ by an element of $Q$
  with all $e_{-,w}$-coordinates $0$ ($W$ has a single eigenvalue of $\U(V_-)(F^+_{\vz})$,
  Lemma \ref{ghl:lem:DRv}), which is trivial on $T_x$ by Lemma \ref{ghl:lem:Tx},
  so its value at $t$ is $1$. Hence $W_E\cap W[m]=c_m(W_B\cap W[m])$ with $c_m=\beta(t)$ for any weight
  $\beta$ of $W[m]$. By the proof of Lemma \ref{ghl:lem:Rrational}, at each $v\in\Sinf$, $h_\infty$ is a positive
  multiple of a form with values in $E$ on $W_B$. Since this form and the form on $\pi_{f,E}$ are positive definite,
  every non-zero de~Rham rational $\varphi$ of pattern $m$ has
  $$
  \langle\varphi,\varphi\rangle\in C|c|^2|c_m|^2E^\times,
  $$
  with $C>0$ independent of $\varphi$ and $m$. If $\beta$, $\tilde\beta$ are weights of $W[m]$, $W[\tilde m]$, then
  $\beta-\tilde\beta-\sum_{v\neq\vz}(m_v-\tilde m_v)\gamma_v$ lies in $Q$ and has all $e_{-,w}$-coordinates $0$, so its
  value at $t$ is $1$. Therefore
  $$
  \frac{c_m}{c_{\tilde m}}
  \in\prod_{v\neq\vz}\gamma_v(t)^{m_v-\tilde m_v}E^\times,
  $$
  and, as $\bar E=E$,
  $$
  \left|\frac{c_m}{c_{\tilde m}}\right|^2\in\prod_{v\neq\vz}|\gamma_v(t)|^{2(m_v-\tilde m_v)}E^\times.
  $$
\end{proof}

We formulate the following conditions. None is proved here.

\begin{hypothesis}[2.10]\label{ghl:hyp:nonvanishing}
  \cite[Conjecture 2.10, p.~20]{ghl} holds for the cohomological conjugate self-dual cuspidal automorphic
  representations to which it is applied in the proof of \cite[Theorem 6.1]{ghl}: for a finite set $S$ of
  non-archimedean places of $F^+$, characters $\alpha_v$ of $\GL_1(\mathcal O_{F^+,v})$, $v\in S$, and a
  conjugate self-dual algebraic character $\alpha_\infty$ of $\GL_1(F\otimes_\QQ\RR)$, there is a Hecke character
  $\chi$ of $\GL_1(\bfA_F)$ with $\chi_\infty=\alpha_\infty$, $\chi|_{\GL_1(\mathcal O_{F^+,v})}=\alpha_v$ for
  $v\in S$, and $L^S(\frac12,\Pi\otimes\chi)\neq0$.
\end{hypothesis}

\begin{hypothesis}[dR]\label{ghl:hyp:dR}
  The hermitian spaces of the induction form a flag $V=V^{(0)}\supset V^{(1)}\supset\cdots$, $\dim
  V^{(k)}=n-k$, and there is one decomposition of $V$ into lines adapted to it, the negative line $V_-$
  lying in every $V^{(k)}$ which occurs. For each pair $\U(V^{(k)})\supset\U(V^{(k+1)})$, with the CM points
  attached to this decomposition as in section \ref{ghl:sec:DR}, the de~Rham structures of
  \cite[\S3.2.2]{ghl} on the minimal types are those of \cite{harris13a} at these points, they are
  defined by models in the sense of section \ref{ghl:sec:DR}, and the
  constructions of \cite{harris13a} and \cite{harris14} for similitude groups transfer to the unitary groups
  \cite[pp.~2045--2046]{harris13a}. If no pair occurs ($q=0$), the same is assumed for $\U(V)$ at the point attached to a
  decomposition of $V$ into lines.
\end{hypothesis}

\begin{hypothesis}[Gal]\label{ghl:hyp:Gal}
  For each pair occurring in the induction and every $\sigma\in\Aut(\CC)$ the $\sigma$-conjugate data satisfy
  the hypotheses of Theorem \ref{ghl:thm:B} (in particular $0\leq q_\sigma\leq n-2$ for the conjugate index) and
  \hypref{ghl:hyp:dR}{dR}, with $E({}^\sigma\pi)=\sigma(E(\pi))$ for each representation $\pi$ of a unitary group occurring; the
  field $E$ of section \ref{ghl:sec:DR} is CM or totally real; every $\rho\in\Aut(\CC/F^{\rm Gal})$
  carries the de~Rham structures, the rational elements of $\Hom_{K'}(\tau_\sigma,\tau)$ and the $E$-valued
  invariant forms for $\sigma$ to those for $\rho\sigma$; and $I^{\rm can}({}^\sigma f,{}^\sigma f')\neq0$ for the
  chosen $f$, $f'$.
\end{hypothesis}

\hypref{ghl:hyp:dR}{dR} is \cite[\S2.3, p.~2011; \S4.2, pp.~2035--2036; pp.~2045--2046]{harris13a} and \cite[\S3.1, p.~87]{harris14}, as we understand them. A point adapted to the whole flag exists by \cite[Lemma 4.1]{harris13a}. \hypref{ghl:hyp:Gal}{Gal} is used only in Theorem \ref{ghl:thm:C} (ii).

\subsection{The theorem}\label{ghl:sec:thmC}

\begin{theorem}\label{ghl:thm:C}
  Let $F^+=\QQ$, $F$ imaginary quadratic, $n\geq2$, and let $\Pi$ be a cohomological conjugate self-dual cuspidal
  automorphic representation of $\GL_n(\bfA_F)$ which satisfies \cite[Hypothesis 2.4]{ghl}, with $\Pi_\infty$
  $(n-1)$-regular. For $0\leq q\leq n-2$ let $\pi(q)$ be as in \cite[\S6.1]{ghl}, and let the fields $E(\pi)$ of
  all representations $\pi$ of unitary groups occurring in the induction be enlarged, as in section
  \ref{ghl:sec:DR}, so that for each pair of the induction $E(\pi(q))=E(\pi'(q))$ is the field $E$ defined
  there. In particular they are stable under complex conjugation. Assume \hypref{ghl:hyp:nonvanishing}{2.10} and
  \hypref{ghl:hyp:dR}{dR}. Here, and in (ii), $Q(\pi(q))$ is the class modulo $E^\times$ of $\langle\varphi,\varphi\rangle$ for a non-zero de~Rham
  rational form $\varphi\in\pi(q)$ whose archimedean component lies in the minimal $K$-type. By Lemma
  \ref{ghl:lem:pattern} it does not depend on $\varphi$.
  \begin{itemize}
    \item[(i)] For $0\leq q\leq n-2$,
      $$
        Q(\pi(q))\in p(\check\xi_\Pi,\Sigma)^{-1}\,\frac{P^{(q+1)}(\Pi,\imath_{\vz})}{P^{(q)}(\Pi,\imath_{\vz})}
        \cdot E^\natural(\Pi,q)^\times,
      $$
      where $E^\natural(\Pi,0)=E_0(\Pi)$, and for $q\geq1$ the field $E^\natural(\Pi,q)$ is the compositum of
      $E_q(\Pi)$, $E(\Pi')$, $E(\pi'(q))$, $E_F(\xi)$, $E_F(\check\xi_\Pi)$, $E_F(\check\xi_{\Pi'})$ and
      $E^\natural(\Pi',n-q-2)$, with $\Pi'$, $\pi'(q)$, $\xi$ as in \cite[\S6.2]{ghl}. Consequently (6.4) of
      \cite[Theorem 6.3]{ghl} holds at the given embedding, with $E^{(i)}(\Pi)$ replaced by the compositum of
      $E_0(\Pi)$ and the $E^\natural(\Pi,q)$, $q\leq i-1$.
    \item[(ii)] If in addition \hypref{ghl:hyp:Gal}{Gal} holds, then (i) holds with $E_q(\Pi)$ in place of $E^\natural(\Pi,q)$, as
      relations equivariant under $\Aut(\CC/F^{\rm Gal})$. That is, \cite[Theorem 6.1]{ghl} and (6.4)
      of \cite[Theorem 6.3]{ghl} hold for $F^+=\QQ$ without \cite[Conjecture 4.15]{ghl}, with $Q(\pi(q))$ taken on
      forms whose archimedean component lies in the minimal $K$-type.
  \end{itemize}
\end{theorem}

\begin{proof}
  Since $\Pi_\infty$, the archimedean component of $\Pi$, is $(n-1)$-regular, the Harish-Chandra parameter of $\Pi$ at $\vz$ satisfies
  $A_i-A_{i+1}\geq n$ \cite[\S1.3.1]{ghl}, so $\lambda_{\vz}$ is regular and \eqref{ghl:eq:gap} holds. The representation $\Pi'$ is
  $(n-2)$-regular \cite[p.~48]{ghl}, and the same holds at every step of the induction, where $F^+=\QQ$
  throughout. For $q=0$, (i) is \cite[Observation 6.6, (6.7), p.~47]{ghl} at the given embedding. The field of
  \cite[(6.7)]{ghl} is $E_0(\Pi)E_F(\check\xi_\Pi)=E_0(\Pi)$, so the reduction which follows it in \cite[p.~47]{ghl}
  is not needed: $E_F(\check\xi_\Pi)$ is generated by $F^{\rm Gal}$, the values of $\check\xi_{\Pi,f}$ and the reflex
  fields of the data $(T_F,h_\Psi)$ \cite[\S2.1.1]{ghl}, which lie in $F^{\rm Gal}$, and since the central character
  of ${}^\sigma\Pi_f$ is $\sigma\circ\xi_{\Pi,f}$, the values of $\check\xi_{\Pi,f}$ lie in $\QQ(\Pi_f)$. Hence
  $E_F(\check\xi_\Pi)\subset\QQ(\Pi_f)F^{\rm Gal}\subset E(\Pi)$ \cite[Theorem 2.6]{ghl}. This does not use
  $F^+=\QQ$. Let $q\geq1$.
  We follow \cite[\S6.2]{ghl}. Step 1, \cite[(6.8)]{ghl} and \cite[Observation 6.9]{ghl}, is unchanged. It uses
  \hypref{ghl:hyp:nonvanishing}{2.10} through \cite[Theorem 5.19]{ghl}. In Step 2 let $f$, $f'$ be as in \cite[p.~32]{ghl}. Then
  \cite[(6.11)]{ghl} holds for them, with the measure of \cite[p.~28]{ghl}, and with $Q(\pi(q))$, $Q(\pi''(q))$
  computed on $f$ and $f''=f'\otimes f_0$ by Lemma \ref{ghl:lem:pattern}, where $\langle f'',f''\rangle=\langle
  f',f'\rangle\langle f_0,f_0\rangle$ for the measure on $\U(V_1)$ of \cite[p.~28]{ghl}. The passage from
  the period $P(f,f'')$ of \cite[\S4.1]{ghl} to \cite[(6.11)]{ghl} at the given embedding uses $|I^{\rm can}(f,f')|^2\in E\cdot\bar E=E$. Since
  $\Sinf=\{\vz\}$, the eigenvector condition of Theorem \ref{ghl:thm:B} is empty, and by \hypref{ghl:hyp:dR}{dR} Proposition
  \ref{ghl:prop:DR} applies to the pair. The measure on $H'$ satisfies the conditions of Lemma \ref{ghl:lem:M} (section \ref{ghl:sec:globalinputs}), so Theorem \ref{ghl:thm:B} and \cite[Proposition 6.12]{ghl} give
  $I^*_{\vz}(f_{\vz},f'_{\vz})\in E^\times$, which is \cite[(6.13)]{ghl} at the given embedding.
  Here \cite[Proposition 6.12]{ghl} rests on $I^{\rm can}(f,f')\neq0$ \cite[Theorem 4.14 (2)]{ghl} and on the formula
  of Ichino--Ikeda and N.~Harris \cite[Theorem 4.5]{ghl}. In the present case this formula is proved in \cite[Theorem
  1.10]{blzz}. Indeed, $\pi(q)$ and $\pi'(q)$ are tempered, their base changes $\Pi$ and $\Pi'$ are cuspidal, the
  latter by \cite[Theorem 4.14 (1)]{ghl}, and they are hermitian in the sense of \cite[Definition 1.5]{blzz},
  i.\,e.\ $L(s,\Pi,{\rm As}^{(-1)^n})$ and $L(s,\Pi',{\rm As}^{(-1)^{n-1}})$ are regular at $s=1$, since $\Pi$ and
  $\Pi'$ are base changes from unitary groups \cite[\S1, p.~2]{ghlap}. These are the values in the denominator of
  \cite[(4.4)]{ghl}. The condition of \cite[Theorem 1.8 (2)]{blzz} that the orthogonal complement of $V'$ in $V$ is
  spanned by a vector of norm $1$ is arranged by scaling the hermitian form of $V$, which changes neither $H$ nor $H'$. The scalar is totally positive, since the complement is positive at every $v\in\Sinf$, $V_-$ lying in $V'$, so the signatures are unchanged.
  Step 3, \cite[(6.14), (6.15)]{ghl}, uses \hypref{ghl:hyp:nonvanishing}{2.10} through \cite[Theorems 2.6, 5.14]{ghl} and is unchanged, and
  so are \cite[Lemma 6.16, (6.17)]{ghl}. At the given embedding they give \cite[(6.18)]{ghl} with the field
  $$
    E_q(\Pi)E_{q'}(\Pi')E_F(\xi)E_F(\check\xi_\Pi)E_F(\check\xi_{\Pi'}),
  $$
  where $q'=n-q-2$ and $E_{q'}(\Pi')=E(\Pi')E(\pi'(q))$. The reduction of this field in
  \cite[p.~50]{ghl} by \cite[Lemma 1.34]{grobnerlin} is a statement about relations
  equivariant under $\Aut(\CC/F^{\rm Gal})$ and is not used. The representations $\Pi'$ and $\pi'(q)$ satisfy the assumptions of the theorem for $n-1$ and
  $q'=n-q-2$ \cite[p.~50]{ghl}, by \cite[Theorem 4.14 (1)]{ghl} and \cite[Remark 2.5]{ghl} for
  \cite[Hypothesis 2.4]{ghl}.
  Hypothesis \hypref{ghl:hyp:dR}{dR} is assumed at every level, and for $n-1=2$ or $q'=0$
  the case $q=0$ applies. By \hypref{ghl:hyp:dR}{dR} the point $x'$
  of the pair at rank $n$ is the point $x$ of the pair at rank $n-1$, so $\pi'(q)$ carries the same de~Rham
  structure in both pairs. Its two Petersson norms agree up to $\QQ^\times$ (section \ref{ghl:sec:globalinputs}), and by Lemma
  \ref{ghl:lem:pattern} the class $Q(\pi'(q))$ in (i) for $\Pi'$ is that of $\langle f',f'\rangle$. Inserting (i)
  for $\Pi'$ into \cite[(6.18)]{ghl} gives (i) for $\Pi$. The relation (6.4) follows from
  \cite[Definition 6.2]{ghl} as in \cite[p.~46]{ghl}.

  For (ii), under \hypref{ghl:hyp:Gal}{Gal} and with $\sigma\in\Aut(\CC)$ as there, Theorem \ref{ghl:thm:B} and Lemma
  \ref{ghl:lem:pattern} apply to the $\sigma$-conjugate data,
  and by \hypref{ghl:hyp:Gal}{Gal} the $E$-valued forms of the proof of Lemma \ref{ghl:lem:pattern} are transported by $\sigma$, so the
  classes $Q({}^\sigma\pi(q))$, $\sigma\in\Aut(\CC)$, are well defined up to $\sigma(E)^\times$,
  which gives condition (1) of \cite[Definition 1.5 (ii)]{ghl}. For $\rho\in\Aut(\CC/F^{\rm Gal})$ the factors
  $\mathrm{Vol}_\infty$ of Lemma \ref{ghl:lem:M}, $s(q)$ and $\dim F_{A'(q)}$ are rational and unchanged, and
  $\rho(R^\circ_{\vz}(\sigma))=R^\circ_{\vz}(\rho\sigma)$, since $\rho$ commutes with complex conjugation on $E$ and
  transports the rational data. This is condition (2), so \cite[(6.13)]{ghl} holds as a relation equivariant under $\Aut(\CC/F^{\rm Gal})$, and
  the rest of \cite[\S6.2]{ghl}, including \cite[Lemma 1.34]{grobnerlin}, applies without change.
\end{proof}

\begin{remark}\label{ghl:rmk:routeII}
  Consider \cite[Theorem 5.6]{ghl2}, which is \cite[Theorem 6.1]{ghl} without \cite[Conjecture 2.10]{ghl}, under
  \cite[Hypothesis 2.3]{ghl2} and $(n+4)$-regularity. Its proof \cite[\S5.3, pp.~33--37]{ghl2} follows that of
  \cite[Theorem 6.1]{ghl} step by step, with \cite[Theorem 1]{ghl2}, which is to appear in \cite{ghlr}, in place of \cite[Theorem 5.19]{ghl} for
  \cite[(5.14)]{ghl2}. The argument above would apply to it for $F^+=\QQ$, under \hypref{ghl:hyp:dR}{dR} and with the
  archimedean input replaced as above, provided the relation (2) of \cite[Theorem 1]{ghl2} at $s_0=\frac12$, i.\,e.\
  \cite[(5.14)]{ghl2}, holds for the pairs $(\Pi,\Pi')$ of the induction. We do not treat \cite[Theorem 5.6]{ghl2} here.
\end{remark}

\subsection{The case $F^+\neq\QQ$}\label{ghl:sec:Fplus}

We keep the flag $V=V^{(0)}\supset\cdots\supset V^{(N)}$ and the decomposition of \hypref{ghl:hyp:dR}{dR}. The space $V^{(k)}$ carries the
representation $\pi^{(k)}$ of the $k$-th step of the induction, with $\pi^{(0)}=\pi(q)$. Here $\Pi^{(k)}$ is the
corresponding representation of $\GL_{n-k}(\bfA_F)$, $q_k$ is its index, and $N=N_q$ is the number of steps, so
that $\pi^{(N)}$ is holomorphic. One checks from the recursion $q_0=q$, $q_{k+1}=n-k-q_k-2$ of \cite[p.~50]{ghl}, which gives $q_{k+2}=q_k-1$,
that $N=0$ for $q=0$ and $N=\min(2q,2n-2q-3)$ for $q\geq1$, that $\dim V^{(N)}\geq2$, and that the rank $n-k$ is at
least $3$ whenever $q_k\geq1$. Let $T^{(k)}$ be the product of the $\U(V_i)$, $V_i\subset V^{(k)}$, and $t^{(k)}$
the element $t$ of Proposition \ref{ghl:prop:DR} for $\U(V^{(k)})$. Fix a line $V_j\subset V^{(N)}$ with $j\neq-$ and put
$\gamma_v=e_{-,v}-e_{j,v}$ for $v\in\Sinf$, at every level. Let $E$ be a number field, stable under complex
conjugation (not necessarily CM, Remark \ref{ghl:rmk:nonCM}), which contains the fields $E$ of section \ref{ghl:sec:DR} for all pairs of the inductions
for all $q$ and
for which $\pi^{(k)}_{f,E}=\bigotimes'_p\pi^{(k)}_{p,E}$ for every $k$, with $E$-structures $\pi^{(k)}_{p,E}$ on the
local components. Such a field exists, after a finite extension of the former fields: the local components
have $E$-structures after a finite extension, containing a spherical vector for almost all $p$, and as $\pi^{(k)}_f$
is irreducible, $\pi^{(k)}_{f,E}$ is a multiple of their restricted tensor product. From now on $E$ denotes this field. Lemma \ref{ghl:lem:pattern}
holds at every level, with $t^{(k)}$. We write $Q(\pi^{(k)};m)$ for the class of the Petersson norm of a non-zero
de~Rham rational form of pattern $m$. Since $V_-$ lies in every $V^{(k)}$, the patterns of all the $\pi^{(k)}$ lie in
the same lattice.

For $[F^+:\QQ]\geq2$ the $\gamma_v$, $v\neq\vz$, are linearly independent modulo the characters trivial on $T_x$
(Lemma \ref{ghl:lem:Tx}), and one expects $\prod_{v\neq\vz}|\gamma_v(t)|^{2k_v}\notin\bar\QQ^\times$ for $k\neq0$. We do not
prove this. If so, the automorphic $Q$-period of \cite[\S4.3]{ghl} depends on the pattern (Remark \ref{ghl:rmk:GHL43}). The proof of \cite[Theorem 6.1]{ghl} uses \cite[\S4.3]{ghl} in
\cite[(6.11)]{ghl}, where only the product $\langle f,f\rangle\langle f'',f''\rangle$ for the chosen pair enters, in
the induction \cite[p.~50]{ghl}, where $Q(\pi'(q))$ occurs at two levels, and in \cite[Observation 6.6]{ghl}, the
base of the induction. Step 3 uses $Q(\pi''(q))\sim Q(\pi'(q))Q(\xi)$, which for the chosen $f''=f'\otimes f_0$ is
$\langle f'',f''\rangle=\langle f',f'\rangle\langle f_0,f_0\rangle$. We choose the forms from the bottom of the
induction upwards.

\begin{lemma}\label{ghl:lem:bottomup}
  Let $\pi=\pi^{(k)}$, $\pi'=\pi^{(k+1)}$, $k<N$, and let $f'\in\pi'$ be factorizable, de~Rham rational, of pattern
  $m'$. Then there is a factorizable de~Rham rational $f\in\pi$ of pattern $-m'$ with $I^{\rm can}(f,f')\neq0$.
\end{lemma}

\begin{proof}
  Let $\pi^{\rm min}$, resp.\ $\pi'^{\rm min}$, be the space of forms whose archimedean components lie in the
  minimal types. The forms $f'\in\pi'^{\rm min}$ with $I^{\rm can}(\cdot,f')=0$ on $\pi^{\rm min}$ form a subspace
  stable under $K'_\infty\times H'(\bfA_{F^+,f})$, which preserves $\pi^{\rm min}$. The space $\pi'^{\rm min}$ is
  irreducible under this group, and the subspace is not all of $\pi'^{\rm min}$ by \cite[Theorem 4.14 (2)]{ghl},
  so it is $0$.
  Since $I^{\rm can}$ is $\CC$-bilinear in $(f,f')$ \cite[p.~29]{ghl} and invariant under
  $\U(V_-)(F^+_v)\subset H'_v$, the functional $I^{\rm can}(\cdot,f')$ vanishes on the forms of pattern $m\neq-m'$.
  By Lemma \ref{ghl:lem:DRv} at the archimedean
  places and the choice of $E$ at the finite places, the factorizable de~Rham rational forms of pattern $-m'$ span
  the forms of pattern $-m'$ in $\pi^{\rm min}$, and one of them has $I^{\rm can}(f,f')\neq0$.
\end{proof}

We formulate the following condition (cf.\ \cite[Observation 6.6, p.~47]{ghl} and the definition of
$P^{(I_0)}(\Pi)$ in \cite[p.~19]{ghl}).

\begin{hypothesis}[B0]\label{ghl:hyp:B0}
  There are $m^{\rm hol}=m^{\rm hol}(\Pi,q)\in\ZZ^{\Sinf\smallsetminus\{\vz\}}$ and a character
  $\varepsilon=\varepsilon_{\Pi,q}$ of $T^{(N)}$ such that for every pattern $m$ of $\pi^{(N)}$
  $$
    Q(\pi^{(N)};m)\in\varepsilon(t^{(N)})\prod_{v\neq\vz}|\gamma_v(t^{(N)})|^{2(m_v-m^{\rm hol}_v)}\cdot
    P^{(I_0)}(\Pi^{(N)})\cdot E^\times .
  $$
\end{hypothesis}

By Lemma \ref{ghl:lem:pattern} it suffices to require this for one pattern. Hypothesis \hypref{ghl:hyp:B0}{B0} says that the class of
$Q(\pi^{(N)};m)/P^{(I_0)}(\Pi^{(N)})$ is a period of $T^{(N)}$ of this form. Here $m^{\rm hol}$ and $\varepsilon$ are not
determined by it, and only the factor $b(\Pi,q)$ below enters. With $b(\Pi,q)$ defined as the class of
$Q(\pi^{(N)};m^{(N)})/P^{(I_0)}(\Pi^{(N)})$, Theorem \ref{ghl:thm:Cprime} does not use \hypref{ghl:hyp:B0}{B0}. In \hypref{ghl:hyp:B0}{B0},
$m^{\rm hol}$ need not be a pattern of $\pi^{(N)}$, and $\varepsilon$ need not lie in $Q$.
Both depend on the chain, hence on $q$. If the arithmetic
holomorphic forms which define $P^{(I_0)}(\Pi^{(N)})$ \cite[p.~19]{ghl} are scalar valued and of a single pattern,
\hypref{ghl:hyp:B0}{B0} holds with $\varepsilon=1$ and $m^{\rm hol}$ their pattern. If they are vector valued, with a metric normalized
on another rational structure, a non-trivial $\varepsilon$ and an $m^{\rm hol}$ which is not a pattern can occur. For
$F^+=\QQ$ there is one pattern, and \hypref{ghl:hyp:B0}{B0} with $\varepsilon=1$ is \cite[Observation 6.6]{ghl}.

If $\varepsilon$ is not trivial on the centre of $\U(V^{(N)})$, the class $\varepsilon(t^{(N)})E^\times$ depends on the
model of the standard principal bundle, through the factor $z$ of Proposition \ref{ghl:prop:DR}, and so does
$Q(\pi^{(N)};m)$. Here \hypref{ghl:hyp:B0}{B0} refers to the fixed model. For the model of Proposition \ref{ghl:prop:DR} (ii),
$t^{(N)}\in T^{(N)}_x(\CC)$, where $T^{(N)}_x$ is the Mumford--Tate group of the
point. Only then is $\varepsilon(t^{(N)})$ a
period of the CM point.

For a pattern $m^{(N)}$ of $\pi^{(N)}$ put
$$
  b(\Pi,q)=\varepsilon(t^{(N)})\prod_{v\neq\vz}|\gamma_v(t^{(N)})|^{2(m^{(N)}_v-m^{\rm hol}_v)},\qquad
  m_q=(-1)^Nm^{(N)} .
$$

\begin{theorem}\label{ghl:thm:Cprime}
  Let $F^+$ be arbitrary, and let $n$, $\Pi$, $\pi(q)$ and the fields $E^\natural(\Pi,q)$ be as in Theorem
  \ref{ghl:thm:C}, and $E$ as above. Assume \hypref{ghl:hyp:nonvanishing}{2.10}, \hypref{ghl:hyp:dR}{dR} and \hypref{ghl:hyp:B0}{B0}. Let $0\leq q\leq n-2$, let $N=N_q$ be as
  above, $m^{(N)}$ a pattern of $\pi^{(N)}$, and $b(\Pi,q)$, $m_q$ as above.
  \begin{itemize}
    \item[(i)]
      $$
        \displaystyle Q(\pi(q);m_q)\in b(\Pi,q)^{(-1)^N}\,p(\check\xi_\Pi,\Sigma)^{-1}\,
        \frac{P^{(q+1)}(\Pi,\imath_{\vz})}{P^{(q)}(\Pi,\imath_{\vz})}\cdot\big(E^\natural(\Pi,q)E\big)^\times.
      $$
    \item[(ii)] For every pattern $m$ of $\pi(q)$,
      $$
        Q(\pi(q);m)\in\prod_{v\neq\vz}|\gamma_v(t^{(0)})|^{2(m_v-m_{q,v})}\cdot Q(\pi(q);m_q)\cdot E^\times.
      $$
    \item[(iii)] (6.4) of \cite[Theorem 6.3]{ghl} holds at the given embedding, with
      $$
        P_i(\Pi,\imath_{\vz})=p(\check\xi_\Pi,\Sigma)\,b(\Pi,i-1)^{-(-1)^{N_{i-1}}}Q(\pi(i-1);m_{i-1})
      $$
      for $1\leq i\leq n-1$ in \cite[Definition 6.2]{ghl}, each $\pi(i-1)$ with its own chain, and $E^{(i)}(\Pi)$
      replaced by the compositum of $E_0(\Pi)$, $E$ and the $E^\natural(\Pi,q)$, $q\leq i-1$.
  \end{itemize}
  For $F^+=\QQ$, with $\varepsilon=1$, this is Theorem \ref{ghl:thm:C} (i), with the field enlarged by $E$.
\end{theorem}

We note that \cite[Theorem 6.1]{ghl} as stated is not expected for every de~Rham rational form when
$[F^+:\QQ]\geq2$ (Remark \ref{ghl:rmk:GHL43}).

\begin{proof}[Proof of Theorem \ref{ghl:thm:Cprime}]
  Choose $f^{(N)}\in\pi^{(N)}$ factorizable, de~Rham rational, of pattern $m^{(N)}$, and for $k=N-1,\dots,0$ let
  $f^{(k)}\in\pi^{(k)}$ be the form of Lemma \ref{ghl:lem:bottomup} for $f'=f^{(k+1)}$. It has pattern
  $(-1)^{N-k}m^{(N)}$. For $k<N$ we have $q_k\geq1$ and rank at least $3$, and we apply \cite[\S6.2]{ghl}, as in the
  proof of Theorem \ref{ghl:thm:C}, to $\Pi^{(k)}$, $\pi^{(k)}$, $\pi^{(k+1)}$ and the pair $f^{(k)}$, $f^{(k+1)}$:
  Steps 1 and 3 and \cite[Lemma 6.16]{ghl} are unchanged, $\Pi^{(k)}_\infty$ being $(n-k-1)$-regular at every
  $v\in\Sinf$ \cite[\S1.3.1 and p.~48]{ghl}. In Step 2, \cite[Lemma 4.2, Theorems 4.5, 4.11, Lemma
  4.8]{ghl} give \cite[(6.11)]{ghl} with
  $$
    \langle f^{(k)},f^{(k)}\rangle\langle f^{(k+1)},f^{(k+1)}\rangle\langle f_0,f_0\rangle
  $$
  in place of $Q(\pi(q))Q(\pi''(q))$, for the measure on $\U(V_1)$ of \cite[p.~28]{ghl}. Theorem \ref{ghl:thm:B}
  (whose eigenvector condition holds for $f^{(k)}$ and $f^{(k+1)}$), Lemma \ref{ghl:lem:M} and
  \cite[Proposition 6.12]{ghl} give \cite[(6.13)]{ghl}. With $r_k$ the right side of (i) for $\Pi^{(k)}$ and $q_k$, without the factor
  $b$, we obtain
  $$
    \langle f^{(k)},f^{(k)}\rangle\langle f^{(k+1)},f^{(k+1)}\rangle\in r_kr_{k+1}(E^\natural(\Pi,q)E)^\times,
  $$
  by \cite[(6.18)]{ghl}. By \hypref{ghl:hyp:dR}{dR} the form $f^{(k+1)}$ has the same de~Rham structure, and by section
  \ref{ghl:sec:globalinputs} its norm is the same up to $\QQ^\times$, in the two pairs in which it occurs. By \hypref{ghl:hyp:B0}{B0} and \cite[Observation 6.6, (6.7)]{ghl}, applied to $\Pi^{(N)}$, with the field as in
  the proof of Theorem \ref{ghl:thm:C},
  $$
    \langle f^{(N)},f^{(N)}\rangle\in b(\Pi,q)\,r_N(E^\natural(\Pi,q)E)^\times .
  $$
  Hence
  $$
    \langle f^{(0)},f^{(0)}\rangle\in b(\Pi,q)^{(-1)^N}r_0(E^\natural(\Pi,q)E)^\times,
  $$
  which is (i) by Lemma
  \ref{ghl:lem:pattern}. Part (ii) is Lemma \ref{ghl:lem:pattern}, and (iii) follows from (i) and \cite[Definition
  6.2]{ghl} as in \cite[p.~46]{ghl}.
\end{proof}

\begin{remark}\label{ghl:rmk:Fplus}\leavevmode\par
  \noindent {\rm(1)} The pattern $m_q$ and the factor $b(\Pi,q)$ are determined by the induction and by $P^{(I_0)}(\Pi^{(N)})$, not
  by $\pi(q)$ alone. That the induction forces opposite patterns at successive levels is Lemma
  \ref{ghl:lem:bottomup}, and
  we know no rule intrinsic to $\pi(q)$ which selects $m_q$. The class
  $$
    Q^\flat(\pi)=Q(\pi;m)\prod_{v\neq\vz}|\gamma_v(t)|^{-2m_v}
  $$
  does not depend on $m$. For $\pi=\pi(q)$ we have $t=t^{(0)}$. If
  $$
    \prod_{v\neq\vz}\big(|\gamma_v(t^{(N)})|^2|\gamma_v(t)|^{-2}\big)^{m^{(N)}_v-m^{\rm hol}_v}\in(E^\natural(\Pi,q)E)^\times,
  $$
  for instance if $\gamma_v(t^{(N)})\in\gamma_v(t)E^\times$ for all $v\neq\vz$, (i) is equivalent to
  $$
    Q^\flat(\pi(q))\in r_0\Big(\varepsilon(t^{(N)})\prod_{v\neq\vz}|\gamma_v(t)|^{-2m^{\rm hol}_v}\Big)^{(-1)^N}
    (E^\natural(\Pi,q)E)^\times .
  $$
  Here $r_0$ is as in the proof of Theorem \ref{ghl:thm:Cprime}. Since (i) holds, the condition is also necessary for this relation. Only the levels $0$ and $N$ enter the condition.

  \smallskip
  \noindent {\rm(2)} The proof of (6.5) of \cite[Theorem 6.3]{ghl} in \cite[p.~47]{ghl} uses a relation between the motivic
  period $Q_i(M(\Pi),\imath)$ and $Q(\pi(i-1))$. By Lemma \ref{ghl:lem:pattern} it can hold for at most one
  pattern, and we do not know whether this is $m_{i-1}$. For $[F^+:\QQ]\geq2$, Theorem \ref{ghl:thm:Cprime} gives
  (6.4), not (6.5).

  \smallskip
  \noindent {\rm(3)} For $[F^+:\QQ]\geq2$ the Galois-equivariant form of Theorem \ref{ghl:thm:C} (ii) would need \hypref{ghl:hyp:Gal}{Gal} and the
  transport of patterns to the conjugate data, whose chains differ. Theorem \ref{ghl:thm:Cprime} is stated at the
  given embedding.

  \smallskip
  \noindent {\rm(4)} If $P^{(I)}(\Pi)$, defined in \cite[p.~19]{ghl} by Petersson norms of arithmetic holomorphic forms,
  depends on the pattern at the definite places of $H_I$ as $Q(\pi)$ does (Lemma \ref{ghl:lem:pattern}), this affects
  \cite[Theorems 2.6, 5.14, 5.19]{ghl}, which we use as stated.
\end{remark}

\begin{remark}\label{ghl:rmk:GHL43}
  In \cite[\S4.3, p.~33]{ghl}, the authors define $Q(\pi)$ by the Petersson norm of a de~Rham rational form which
  ``generates'' the de~Rham rational structure. Every non-zero rational vector of an absolutely irreducible module
  generates it, so no vector is singled out. Assume that the de~Rham structures of \cite{ghl} are those of
  \cite{harris13a} at the points of \hypref{ghl:hyp:dR}{dR} and are defined by models in the sense of section \ref{ghl:sec:DR},
  so that Proposition \ref{ghl:prop:DR} applies, and that $\prod_{v\neq\vz}|\gamma_v(t)|^{2k_v}\notin\bar\QQ^\times$ for
  $k\neq0$, as predicted by the period conjecture. Then by Lemma \ref{ghl:lem:pattern}, for $[F^+:\QQ]\geq2$ and $\pi$
  with at least two
  patterns, as is the case for $\pi(q)$ by the regularity of $\Pi_\infty$ \cite[p.~12]{ghl}, the independence
  asserted in
  \cite[\S4.3]{ghl}, and in \cite[Proposition 3.19, Lemma 3.20]{harris13a}, does not hold, and $Q(\pi)$ is not
  well defined as stated. In this case \cite[Theorem 6.1]{ghl} is not meaningful as stated.
  If \cite[Conjecture 4.15]{ghl}
  is read for all rational vectors, it is not expected at the definite places (Remark
  \ref{ghl:rmk:definiteperiods} (2)), and the theorem, which assumes it, is vacuous. If it is read for suitably
  chosen vectors, then \cite[Theorem 6.1]{ghl}, in the form proved in \cite[\S6.2]{ghl}, holds for the forms of at
  most one pattern, and Theorem \ref{ghl:thm:Cprime} is the corresponding statement.
  The same question arises for \cite[Theorem 5.6]{ghl2},
  whose proof follows that of \cite[Theorem 6.1]{ghl} (Remark \ref{ghl:rmk:routeII}). The transcendence is not proved. For special $F$, for instance with CM types induced from a subfield, a
  relation between CM periods might make $|\gamma_v(t)|^2$ algebraic.
\end{remark}

For $q\geq1$ Theorems \ref{ghl:thm:C} and \ref{ghl:thm:Cprime} use \hypref{ghl:hyp:nonvanishing}{2.10}, a global input of the proof of
\cite[Theorem 6.1]{ghl} which is independent of the archimedean integrals. For $[F^+:\QQ]\geq2$, Theorem
\ref{ghl:thm:Cprime} needs in addition \hypref{ghl:hyp:B0}{B0} and concerns $Q(\pi(q))$ on a pattern fixed by the
induction, up to the factor $b(\Pi,q)^{(-1)^N}$. The relation (6.5) of \cite[Theorem 6.3]{ghl} uses
further inputs \cite[pp.~46--47]{ghl} and is not treated here.

\end{document}